\documentclass[10pt]{amsart}
\usepackage{amsmath,amsfonts,amssymb}
\usepackage{amssymb,url,xspace}
\usepackage{mathtools}
\usepackage[margin=2.1cm]{geometry}
\usepackage[english]{babel}
\usepackage[T1]{fontenc}
\usepackage[sort&compress,numbers]{natbib}

\numberwithin{equation}{section}
\usepackage{amssymb,url,xspace}
\usepackage{amsmath,amsfonts,amssymb}
\usepackage{mathtools}
\usepackage{mathrsfs}
\usepackage{textcomp}
\usepackage{fancyhdr}
\usepackage{eso-pic}
\usepackage{bbold}
\usepackage{comment}
\usepackage{color}
\usepackage{xcolor}
\usepackage{stmaryrd}
\usepackage[artemisia]{textgreek}

\def\XXint#1#2#3{{\setbox0=\hbox{$#1{#2#3}{\int}$}
     \vcenter{\hbox{$#2#3$}}\kern-.5\wd0}}

\let\dem=\proof
\let\enddem=\endproof
\let\vph=\varphi

\DeclarePairedDelimiter\abs\lvert\rvert
\DeclarePairedDelimiter\norm\lVert\rVert
\DeclarePairedDelimiter\scalar\langle\rangle
\DeclareMathOperator\dvg{div}

\DeclareMathOperator\adj{Adj}
\newcommand\R{\mathbb{R}}
\newcommand\N{\mathbb{N}}
\newcommand\D{\mathbb{D}}

\newcommand\X{\mathscr{X}}
\newcommand\dpt{\partial_t}

\newcommand\cC{\mathscr{C}}

\newcommand\loc{\text{loc}}
\newtheorem{theo}{Theorem}[section]
\newtheorem{lemm}[theo]{Lemma}

\theoremstyle{definition}

\newtheorem{prop}[theo]{Proposition}
\newtheorem{nota}[theo]{Notation}
\newtheorem{rema}[theo]{Remark}
\newtheorem{coro}[theo]{Corollary}

\definecolor{mycolor}{HTML}{D35400}

\definecolor{hide}{HTML}{A0AC81}
\definecolor{myred}{HTML}{8A1538}

\usepackage[colorlinks=true,citecolor=blue,urlcolor=cyan]{hyperref}
\usepackage{cleveref}
\crefformat{prop}{#2 Proposition~#1#3}{}
\crefformat{lemma}{#2Lemma~#1#3}{}
\crefformat{rema}{#2Remark~#1#3}{}
\crefformat{theo}{#2Theorem~#1#3}{}
\crefformat{appendix}{#2Appendix~#1#3}{}
\crefformat{section}{#2Section~#1#3}{}
\crefformat{coro}{#2Corollary~#1#3}{}
\crefformat{nota}{#2Notation~#1#3}{}
\crefformat{def}{#2Definition~#1#3}{}
\begin{document}

\title[]{A well-posedness result for the compressible two-fluid model with density-dependent viscosity}

\author{}
\address{}
\curraddr{}
\email{}
\thanks{}

\author[]{Sagbo Marcel ZODJI}
\address{Université Paris Cité and Sorbonne Université, CNRS, IMJ-PRG.}
\email{sagbo-marcel.zodji@u-paris.fr}
\thanks{}


\date{\today}

\dedicatory{}


\begin{abstract}
In this paper, we study a system of PDEs describing the motion of two compressible viscous fluids occupying  the whole space $\R^d\;(d\in \{2,3\}$). The two phases of the mixture are separated by a $\cC^{1+\alpha}$-regular sharp interface $\mathcal{C}$ across which the density can experience jumps.  We prove the existence of a unique local-in-time  solution assuming that 
 the initial density is $\alpha$-H\"older continuous on both sides of $\mathcal{C}$. The initial velocity belongs to the Sobolev space $H^1(\R^d)$, and the divergence of the initial stress tensor  belongs to $L^2(\R^d)$. The later assumption expresses somehow the continuity of the stress tensor.  This result is more general than the one by Tani \cite{tani1984two}, as it allows for less regular initial data and furthermore it can serve as  a building block for the construction of global-in-time solutions.
\end{abstract}

\maketitle

\section{Introduction and main result}

\subsection{Presentation of the model}
In this paper, we are concerned with the following system of PDE describing the evolution of two compressible fluids with density-dependent viscosity in the whole space $\R^d\; (d\in \{2,3\})$,
\begin{gather}\label{epNS}
\begin{cases}
\dpt \rho_\pm +\dvg (\rho_\pm u)=0,\\
\dpt(\rho u^j)+\dvg (\rho u^j u)=\dvg \Pi^j.
\end{cases}
\end{gather}
Above, $\rho_\pm=\rho_\pm(t,x)\geqslant 0$ is the density of the phase fluid labelled $\pm$, $\rho:=\rho_++\rho_-$ is the total 
density of the mixture and $u=u(t,x)\in \R^d$ is the single velocity of the mixture. The stress tensor $\Pi$ is given by 
\begin{gather}\label{ep1}
\Pi^{jk}= 2\mu \D^{jk} u+ (\lambda\dvg u- P+\widetilde P)\delta^{jk}\quad\text{with}\quad \D^{jk} u=\dfrac{1}{2}\big(\partial_j u^k+ \partial_k u^j\big),
\end{gather}
where $P=P(\rho_+,\rho_-)$, $\mu= \mu(\rho_+,\rho_-)$, $\lambda=\lambda(\rho_+,\rho_-)$ are respectively pressure, dynamic and kinematic viscosity of the mixture and are given $\cC^2$-regular functions of $\rho_\pm$, and finally $\widetilde P$ is a positive equilibrium pressure.  We refer to \cite{Ishi75} for the derivation of \eqref{epNS} for the more general case of  heat conducting fluids mixture. 

It is assumed that the two pure components of the mixture are immiscible, so that each fluid is implicitly assumed to occupy its own domain and that $\rho_+\rho_-=0$. In particular they are separated by a sharp interface  $\mathcal{C}$. The relevant physical situation for multi-phase flows is exactly where the total density $\rho$ experiences  discontinuity across the interface $\mathcal{C}$. Since the two densities $\rho_\pm$ of the mixture are assumed to be different, we  see that the class of \emph{piecewise regular} functions is suitable for the study of this system. This motivates the following definitions.

\subsubsection{Piecewise regularity}\label[section]{notabene}
\noindent 
\begin{enumerate}
\item \label{c3.24} Throughout this paper $\alpha\in (0,1)$, and we assume that $\mathcal{C}$ is  a compact $\cC^{1+\alpha}$-regular manifold of dimension $d-1$, which is the boundary of an open, bounded and simply connected domain $D\subset \R^d$. It is defined as follows.

There exists a finite number $J\in \N^*$ of open bounded sets $V_j\subset \R^{d-1}$, $j\in \llbracket 1,J\rrbracket$ and maps
\begin{equation*}
 \gamma_j\colon  V_j\longrightarrow \mathcal C\subset\R^d, \quad
 \mathcal{C} = \bigcup_{j=1}^J \gamma_j (V_j),
\end{equation*}
such that $\gamma_j\in \cC^{1+\alpha}(V_j)$ are homeomorphisms on their images. For $(j,k)\in \llbracket 1, J\rrbracket\times \llbracket 1, J\rrbracket$ such that  $\gamma_j(V_j)\cap\gamma_k( V_k)\neq \emptyset$, the map 
\[
\gamma_j^{-1} \circ \gamma_{k}\colon \gamma_k^{-1}\big(\gamma_j(V_j)\cap\gamma_k( V_k)\big) \longrightarrow \gamma_j^{-1}\big(\gamma_j(V_j)\cap\gamma_k( V_k)\big)
\]
is a $\cC^{1+\alpha}$-diffeomorphism. Furthermore, these mappings provide an analytical expression for the outward normal vector field $n_x$ on $\mathcal{C}$. We introduce the following notations:
\[
\norm*{\mathcal{C}}_{\text{Lip}}:=\sup_j\norm{\nabla \gamma_j}_{L^\infty(V_j)},\,\text{ and }\,  \norm*{\mathcal{C}}_{\dot \cC^{1+\alpha}}:=  \sup_j\norm{\nabla \gamma_j}_{\dot \cC^\alpha(V_j)}.
\]
 
 Following \cite[page 6]{gancedo2021quantitative},  we define the quantity:
 \begin{gather}\label{c3.14}
 \abs*{\mathcal{C}}_*:=\Bigg(\inf_j \bigg(\inf_{\substack{s\neq s'\\
s,s'\in V_j}
}\dfrac{\abs{\gamma_j(s)-\gamma_j(s')}}{\abs{s-s'}}\bigg)\Bigg)^{-1}\in (0,\infty).
 \end{gather}

For $(j,k)\in \llbracket 1, d\rrbracket\times \llbracket 1, d\rrbracket$, we define the second-order 
Riesz operator as $\mathcal R_{jk}=-\partial_j (-\Delta)^{-1}\partial_k$. The piecewise H\"older estimate in \cref{lemA2} for the operator $\mathcal{T} = \mathcal{R}_{jk}$ (respectively, $\mathcal{T} = \mathcal{R}_{jk} \mathcal{R}_{lm}$) involves a polynomial denoted by $\mathfrak{P}^{jk}$ (respectively, $\mathfrak{P}^{jklm}$). Next, we consider a polynomial $\mathfrak{P}$ that satisfies:
\[
\mathfrak{P}^{jk}(a),\, \mathfrak{P}^{jklm}(a)\leqslant \mathfrak{P}(a),\; \text{ for all }\; a>0.
\]
Finally, we introduce the following quantity: 
\begin{gather}\label{c3.34}
 \mathfrak{P}_{\mathcal{C}}:=\big(1+\abs*{\mathcal{C}}\big)\mathfrak{P}\big( \norm*{\mathcal{C}}_{\text{Lip}}+\abs*{\mathcal{C}}_*\big)\norm*{\mathcal{C}}_{\dot \cC^{1+\alpha}}.
\end{gather}

 As noted in \cite[Section 3.1]{kiselev2016finite} (see also \cite[Remark 2, page 35]{finn1957asymptotic}),  there exists  a function $\vph\colon \R^d\mapsto \R\in \cC^{1+\alpha}$ such that
\[
D= \{x\in \R^d\colon \vph(x)>0\},\quad\text{and} \quad \abs{\nabla \vph}_{\text{inf}}:=\inf_{x\in\partial D}\abs{\nabla \vph(x)}>0. 
\]
We then define:
\begin{gather}\label{c3.55}
\ell_{\vph}=\min\left\{1,\left(\dfrac{\abs{\nabla \vph}_{\text{inf}}}{\norm{\nabla\vph}_{\dot \cC^\alpha}}\right)^{1/\alpha}\right\}.
\end{gather}
\item Since
\[
\R^d= D\cup \mathcal{C}\cup (\R^d\setminus \overline{D}),
\]
we   define the space of piecewise $\alpha$-H\"older continuous  functions with respect to $\mathcal{C}$ as follows: 
\[
\dot \cC^\alpha_{pw,\mathcal{C}}(\R^d):=\dot \cC^\alpha(\overline{D})\cap \dot \cC^\alpha(\R^d\setminus D),\;\,\text{with}\;\,\norm{g}_{\dot \cC^\alpha_{pw,\mathcal{C}}(\R^d)}:= \norm{g}_{\dot \cC^\alpha(\overline{D})}+\norm{g}_{\dot \cC^\alpha(\R^d\setminus D)}
\]
and the non-homogeneous space
\[
 \cC^\alpha_{pw,\mathcal{C}}(\R^d):= L^\infty(\R^d)\cap \dot \cC^\alpha_{pw,\mathcal{C}}(\R^d)\;\text{ with }\; \norm{g}_{\cC^\alpha_{pw,\mathcal{C}}(\R^d)}:=\norm{g}_{L^\infty(\R^d)}+\norm{g}_{\dot \cC^\alpha_{pw,\mathcal{C}}(\R^d)}.
\]
This space strictly contains the H\"older space $\cC^\alpha(\R^d)$. 
\item Given a function $g\in \cC^{\alpha}_{pw,\mathcal{C}}(\R^d)$, the jump $\llbracket g\rrbracket$, and the average $\scalar{g}$ of $g$ on $\mathcal{C}$ are defined as follows: for all $\sigma\in \mathcal{C}$,
\[
\llbracket g\rrbracket(\sigma):=\lim_{r\to 0}\big[g(\sigma+r n_x(\sigma))-g(\sigma-r n_x(\sigma))\big]\text{ and }\scalar{g}(\sigma):=\dfrac{1}{2}\lim_{r\to 0}\big[g(\sigma+r n_x(\sigma))+g(\sigma-r n_x(\sigma))\big].
\]
Above $n_x$ denotes the outward normal vector of $\mathcal{C}$.
\item   Consider a time-dependent interface $\mathcal{C}=\mathcal{C}(t)$, with  local parametrizations $\gamma_j=\gamma_j(t,s)$,\, $j\in \llbracket 1, J\rrbracket$ defined on $I\times V_j\subset \R\times \R^{d-1}$. We assume that $\gamma_j\in \cC(I,\cC^{1+\alpha}(V_j))$, and for all $t\in I$, $\mathcal{C}(t)$ is a $\cC^{1+\alpha}$-regular manifold of dimension $d-1$, that forms the boundary of an open, bounded and simply connected domain $D(t)$. We introduce the following space:
\[
L^p( I, \;\; \cC^\alpha_{pw,\mathcal{C}}(\R^d)):=\left\{ g=g(t,x)\colon 
\begin{cases}
    \displaystyle\int_I\norm{g(t)}_{ \cC^\alpha_{pw,\mathcal{C}(t)}(\R^d)}^pdt<\infty\quad \text{ if }\quad 1\leqslant p<\infty,\\
    \displaystyle\sup_{t\in I}\text{ess} \norm{g(t)}_{ \cC^\alpha_{pw,\mathcal{C}(t)}(\R^d)}<\infty \quad \text{ if }\quad p=\infty
\end{cases} \right\}.
\]
 For $I=(0,T)$, the corresponding space will sometimes be shortly  denoted by $L^p_T\cC^\alpha_{pw,\mathcal{C}}(\R^d)$.
\end{enumerate}
\subsubsection{Reformulation of the system}
A suitable way to rewrite the system \eqref{epNS} is to introduce the so-called fractional density or concentration $c$:  
\[
c:=\dfrac{\rho_+}{\rho}.
\]
From the first two equations of \eqref{epNS}, one checks easily, assuming some regularity, that the concentration solves 
the transport equation 
\begin{gather}\label{c4.4}
    \dpt c+ u\cdot \nabla c=0.
\end{gather}
Let us notice that the domain occupied by the phase $\pm$ can be defined by:
\[
D_\pm(t):=\{x\in \R^d\colon  \rho_\pm (t,x)>0\}.
\]
As well, one can  write the domain occupied by the one phase fluid with the help of the concentration as follows:
\begin{gather}\label{c4.5}
D_+(t)=\{x\in \R^d\colon  c (t,x)=1\},\;\quad D_-(t)=\{x\in \R^d\colon  c (t,x)=0\}
\end{gather}
and consequently the interface is:
\[
\mathcal{C} (t):=\partial D_+ (t)= \partial D_-(t).
\]
Assume that the velocity is regular, for instance Lipschitz, then it admits a continuous flow map $X(t)$ that helps together with \eqref{c4.4} and \eqref{c4.5} to obtain that the initial interface is just transported by the fluid
\[
\mathcal{C}(t)=X(t) \mathcal{C}(0).
\]
On the other hand, summing  up the mass equation on the density of the one phase $\rho_\pm$ $\eqref{epNS}_1$, one obtains that the total density solves a mass equation
with velocity $u$. Besides, the density $\rho_\pm$ of the one phase fluid can be written in terms of the total density and the concentration as follows:
\[
\rho_+= c\rho \quad \text{ and } \quad \rho_-=(1-c)\rho, 
\]
thus, the pressure can be written as function of the total density and the concentration, $P=P(\rho,c)$, such that 
$P_+(\rho)= P(\rho,1)$ and $P_-(\rho)=P(\rho,0)$ as well as the viscosity. Hence, formally $(c,\rho,u)$ solves:
\begin{gather}\label{ep1.1}
    \begin{cases}
        \dpt c+ u\cdot\nabla c=0,\\
        \dpt \rho + \dvg (\rho u)=0,\\
        \dpt(\rho u)+ \dvg(\rho u\otimes u)+\nabla P(\rho,c)=\dvg(2\mu(\rho,c)\D u)+\nabla (\lambda(\rho,c)\dvg u).
    \end{cases}
\end{gather}
Conversely, if $(c,\rho, u)$ solves \eqref{ep1.1}, one can set 
\[
\rho_+= c\rho \quad \text{and}\quad \rho_- =(1-c)\rho
\]
and checks easily that $\rho_\pm$ solve the equations $\eqref{epNS}_1$. Just by replacing 
\[
\rho=\rho_++\rho_- \quad\text{and} \quad c=\dfrac{\rho_+}{\rho_++\rho_-}
\]
in the pressure and the viscosity, one infers that $(\rho_+,\rho_-, u)$ solves \eqref{epNS}. In fact, we obtain that 
$(\rho_+,\rho_-,u)$ solves \eqref{epNS} if and only if $(c,\rho,u)$ solves \eqref{ep1.1} owing to some change of pressure and viscosity law, at least for classical solutions. The equivalence also holds  true for weak solutions  \emph{à la Lions} (we mean weak solutions which satisfy the classical energy inequality associated with \eqref{epNS}, see \cite{vasseur2019global}) due to 
the  result in \cite[Lemma 2.5]{vasseur2019global}, which  is part of the  DiPerna-Lions theory.  Therefore, we will only work with the system \eqref{ep1.1}.

The equations \eqref{ep1.1} are supplemented with initial data 
\begin{gather}\label{ep3}
c_{|t=0}=c_0,\quad\;{\rho}_{|t=0}=\rho_{0}\quad \text{and}\quad u_{|t=0}=u_0.
\end{gather}
Given the equilibrium state of the pressure $\widetilde P$ introduced in \eqref{ep1},  the equilibrium state of the density, $\widetilde\rho$, is defined by the following relation:
\begin{gather}\label{c4.36}
\widetilde P=P(\widetilde\rho,1).
\end{gather}
Next, we define the equilibrium states of the viscosities:
\begin{gather*}
 \widetilde \mu:=\mu(\widetilde\rho,1)\quad \text{and}
\quad \widetilde \lambda:=\lambda(\widetilde\rho,1).
\end{gather*}

We recall that $\mu$, $\lambda$ and $P$ are $\cC^2$-regular, with $P$ and $\lambda$ being non-negative, while $\mu(\rho,c)>0$ for $\rho>0$.

The purpose of this paper is to  construct a  unique local-in-time solution $(c,\rho,u)$ of the equations \eqref{ep1.1} with discontinuous initial data and with density-dependent viscosity: the initial  density $\rho_0$ is assumed to be H\"older continuous on both sides of a suitable interface 
$\mathcal{C}_0$ through which it is discontinuous. The initial velocity belongs to the Sobolev space $H^1(\R^d)$
with the following compatibility condition:
\begin{gather}\label{c4.1}
\dvg(\Pi)_{|t=0}\in L^2(\R^d).
\end{gather}
This condition expresses the continuity of the normal component of the stress tensor and does not require smoothness of the density.
The parabolic effect of the momentum equations ensures that \eqref{c4.1} holds true at positive times even for less regular initial data (see Hoff \cite{hoff1995global}). In particular it does not  prevent to extend solutions until blow-up.

Our result is more general than the one by Tani \cite{tani1984two} (see \cref{review} below for more details), as it allows for less regular initial data, interface, and it can serve as building blocks for global solutions.

We do not have enough information in order  to prove existence of global weak solution for the two-fluid model, even in the case of small initial data. However, we are able to prove 
    a global well-posedness result in the spirit of \cite{hoff2002dynamics} for a simplified situation which is the Navier-Stokes equations with density-dependent viscosity. To keep the paper a reasonable size, we plan to address this issue in the near future.
\subsection{Review of known results}\label[section]{review}

The study of the two-phase model \eqref{epNS} goes back to Tani’s work \cite{tani1984two} where the author showed existence of a local-in-time solution for the heat conductive variant of \eqref{epNS}.  As in the pioneering work of Nash \cite{nash1962probleme}, Tani's proof relies on a change  of variables into Lagrangian coordinates, followed by the construction of an approximate solution using the kernel associated to the linearized equations.
 The solution is constructed in piecewise H\"older spaces\,: the density is $\cC^{1+\alpha}$ and the velocity is $\cC^{2+\alpha}$,
 on both sides of an interface $\mathcal{C}$ which is $\cC^{2+\alpha}$.  To the best of our knowledge, Tani's result is the only one addressing the two-fluid model with H\"older regularity for interfaces.

 In the literature, other results mostly focus on interfaces with Sobolev regularity.  As an example, Denisova in \cite{denisova2000problem}  has proved the local well-posedness of the system \eqref{epNS} including surface tension with different constant  viscosity for each component of the mixture and with Sobolev-Slobodetskii interface $\mathcal{C}\in W^{\frac{5}{2}+l,2}$, 
 with $l\in (1/2,1)$. Her proof also relies on a change in  Lagrangian coordinates and the density and the velocity are piecewise $W^{1+l,2}$ regular on both sides of the interface $\mathcal{C}$.
 
 In contrast to the aforementioned results based on a change in Lagrangian coordinates, there are other works \cite{jang2016compressible,jang2016compressible-stability,jang2016compressible-taylor} by  Jang \emph{et al.} that rely on a change into flattened coordinates. 
The authors considered the two-phase problem in a cylinder, taking into account scenarios with and without surface tension. This problem involves two fluids, with one fluid on top of another. In this context, there are two interfaces to consider: the interface between the lower and upper fluids, and the interface between the upper fluid and the atmosphere. Both interfaces are assumed to be represented by graphs. Additionally, the viscosity of the two fluids are assumed to be constant and different from each other.  In \cite{jang2016compressible}, they successfully constructed local-in-time solution in Sobolev 
 framework, in \cite{jang2016compressible-taylor}, they studied the Rayleigh-Taylor instability which occurs when the heavier fluid is on top of a lighter fluid. They successfully proved that when the surface tension is below a critical value, then the system is 
 nonlinearly unstable. Finally in \cite{jang2016compressible-stability}, they proved that if the heavier fluid is on top of the lighter one  and the surface tension is sufficiently large, then, the Rayleigh-Taylor instability is prevented and in the converse 
 configuration, the system is stable around the equilibrium state with non-negative surface tension. Moreover, they 
 studied the  zero surface tension limit problem.
 
 On the other hand, there are also results \cite{shibata2015local,yan2020global} that address the two-phase problem involving a compressible fluid and an incompressible fluid.
In \cite{shibata2015local}, Shibata has proved the local well-posedness of this  two-phase model, including temperature and surface tension, with an almost flat interface. The author worked on an exterior domain and with a viscosity that depends on both density and temperature. The initial density of the compressible fluid has a Sobolev $W^{1,p}$ regularity and the initial velocity and temperature have the  Besov regularity  $B^{2(1-1/p)}_{q,p}$, on both sides of the interface, for some $1<p,q<\infty$ that verify:
 \[
\dfrac{2}{p}+\dfrac{d}{q}<1.
 \]
 The same Compressible-Incompressible model is considered in \cite{yan2020global} by Yan  and Zhao with surface tension in horizontal periodic domain. The viscosity of the fluids involved are different from each other and assumed to be constant. 
 They successfully achieved a global well-posedness result by making the assumption that the initial density of the compressible fluid is close to a constant state in $H^4$  with respect to the domain occupied by  the compressible fluid, and the initial velocity is small in $H^4$ on both sides of the interface.  To our knowledge, this result stands as the only one to address the existence of strong global solution of the two-phase problem.  

 Finally, we may notice that Kubo  \emph{et al.} in \cite{kubo2016some} have addressed the local well-posedness of the Compressible-Compressible model with different approach. Their methods rely on proving the existence of a $\mathcal{R}$-bounded solution operator and maximal regularity of the linearised system. The interface is uniformly $W^{r,2-1/r}$ for $r>d$, the density has a $W^{1,r}$ regularity on both sides of the interface and the initial velocity belongs to the Besov space  $B^{2(1-1/p)}_{r,p}$ for some $p>2$  on both sides of the initial interface. More recently, in \cite{shibatakuo2023}, Kuo and Shibata have constructed local-in-time solution of the one phase model in half-space with non-slip boundary conditions. As the result above, the solution is also formulated, in some maximal regularity Besov space framework.

As explained, for instance in \cite{gancedo2021global}, a different approach to tracking the interface is to construct weak 
solutions for the full model in a class that allows the study of its dynamic. To the best of our knowledge, the only 
available results for weak solutions pertain to fluid mixture with constant viscosity that are the same for both components. 
For instance, we refer to the work by Vasseur, Wen, and Yu \cite{vasseur2019global} in a three-dimensional bounded domain where they obtained a weak solution with finite initial energy. The main difficulty in their proof lies in the compactness of an approximate sequence of the pressure that depends on two variables, in contrast to the single-variable case \cite{lions1996mathematical, feireisl2001existence}. They successfully overcome the difficulty by considering the following pressure law
\[
P(\rho_+,\rho_-)=A_+{\rho_+}^{\gamma_+}+A_-{\rho_-}^{\gamma_-}
\]
 with the following closeness constraints on $\gamma_\pm$: 
\[
\gamma_\pm>9/5\quad \text{and}\quad \max\left\{\dfrac{3\gamma_+}{4}; \gamma_+-1; \dfrac{3(\gamma_++1)}{5}\right\}< \gamma_- < \min \left\{\dfrac{4\gamma_+}{3}; \gamma_++1; \dfrac{5\gamma_+}{3}-1\right\}.
\]
This result has been improved by Novotn\'y  \cite{novotny2020weak1}, Wen \cite{wen2021global} by allowing more general pressure law,  Novotn\'y and Pokortn\'y \cite{novotny2020weak} and Novotn\'y \emph{et al.} \cite{kracmar2022weak} for multi-components fluid models. It is important to notice that, although the pressure law of the components of the mixture are different, the viscosity of the components  are supposed to be the same constant. The solutions that the authors constructed are too weak in order to track down discontinuities in the density. In one dimension, Bresch, Burtea and Lagoutière in \cite{bresch2022mathematical} successfully solved  the global well-posedness for non common constant viscosity,  with numerical considerations in \cite{bresch:hal-03738116}. 

 A class of weak solutions with more regularity (for the one phase model) has been proposed in the pioneering work of  Hoff \cite{hoff1995global}, see also \cite{perepelitsa2006global,vaigant1995existence} for bulk viscosity coefficient $\lambda$ depending on the density and Perepelitsa \cite{perepelitsa2014weak} for the case of bounded domains. However, these papers do not address the density patch problem.  
This issue is explored in Hoff \cite{hoff2002dynamics} and Hoff and Santos \cite{hoff2008lagrangean}, which examine the dynamics of sharp interfaces for the single-phase model, namely the Navier-Stokes equations with constant viscosity.
Persistence of $\cC^{1+\alpha}$-regular interface  is achieved in \cite{hoff2002dynamics} for $d=2$, under the condition that the density is  piecewise $\alpha$-H\"older continuous. However, 
 the analysis is restricted to a linear pressure law and small viscosity $\lambda$. In \cite{hoff2008lagrangean} the density is only bounded, the gradient of the velocity is only $BMO$ in space and for interfaces that are initially $\cC^\alpha$ one can only ensure $\cC^{\alpha(t)}$-regularity for $t>0$ where $\alpha(t)$ is decaying exponentially to zero. Although not stated, a similar result i.e. propagation of H\"older regularity with exponential loss  also holds in the context of \cite{perepelitsa2006global,vaigant1995existence} in which the viscosity $\lambda$ is allowed to depend on the density, the analysis being very similar to \cite{hoff2008lagrangean}. However, when the viscosity $\mu$ depends on the density, there is no clear notion of effective flux, the analysis complicates and it is not even clear how can one propagate the $L^\infty$-norm of the density, we refer to an interesting work by Bresch and Jabin \cite{bresch2018global} for a discussion regarding the complications arising in this case.

In this paper, we consider density-dependent viscosity for both components of the mixture, and we achieve the construction of unique weak solution in $\R^d,\, (d\in \{2,3\})$ for the system \eqref{ep1.1} which allows the study of the dynamics of the interface. 
\subsection{Main result}
We consider the Cauchy problem associated with equations \eqref{ep1.1} and initial data \eqref{ep3}. We assume the
existence of an open, bounded, and simply connected domain $D_0$ such that:
\begin{gather}\label{epq2}
c_0=\mathbb{1}_{D_0^c}. 
\end{gather}
Also, we assume that the boundary $\mathcal{C}_0:=\partial D_0$ fulfills the conditions outlined  in \cref{notabene} \cref{c3.24}.
Additionally, the initial density $\rho_0$ and velocity $u_0$  satisfy:
\begin{gather}\label{epq1}
    u_0\in H^1(\R^d), \quad P(\rho_0,c_0)-\widetilde P\in L^2(\R^d)\cap \cC^\alpha_{pw,\mathcal{C}_0}(\R^d),\quad\text{and}\quad \dvg(\Pi)_{|t=0}\in L^2(\R^d).
\end{gather}
A technical step requires us to obtain estimate for the $L^{8/d}((0,T),L^4(\R^d))$ norm of $\sigma^{\tfrac{d}{4}}\nabla \dot u$, where 
\[
\sigma(t)=\min \{1,t\}\quad\text{and}\quad\dot u:= \dpt u+ (u\cdot\nabla)u
\]
is the material acceleration of fluid particles. To achieve this,  we need to look at  the second material derivative of the velocity, which is 
\[
\ddot u:= \dpt \dot u+ (u\cdot\nabla)\dot u.
\]
 This variable appears while applying the material derivative $\dpt\,\cdot+\dvg(\; \cdot\; u)$ to the momentum equation $\eqref{ep1.1}_3$. 
 Our main result reads as follows:

\begin{theo}\label[theo]{thlocal}
        Let $(c_0, \rho_0, u_0)$ be initial data associated with the system \eqref{ep1.1}, satisfying the conditions \eqref{c4.36},  \eqref{epq2} and  \eqref{epq1}.  Additionally, assume:
    \begin{gather}\label{c3.37}
    \rho_{*,0}:= \inf_{x\in \R^d}\rho_0(x)>0 \quad \text{ and } \quad \mu_{*,0}:= \inf_{x\in \R^d}\mu(\rho_0(x),c_0(x))>0.
    \end{gather}
    There exists $[\mu]>0$ depending only on the dimension $d\in \{2,3\}$, $\alpha$ and $\widetilde\mu$ such that if \footnote{We refer to \eqref{c3.34}-\eqref{c3.55} for the definition of $\ell_{\vph_0}$ and $\mathfrak{P}_{\mathcal{C}_0}$.} 
    \begin{multline}
    \left[1+\norm{\lambda(\rho_0,c_0)}_{\dot \cC^\alpha_{pw,\mathcal{C}_0}(\R^d)}+\left( \mathfrak{P}_{\mathcal{C}_0} +\ell^{-\alpha}_{\vph_0}\right)\norm*{\llbracket\lambda(\rho_0,c_0)\rrbracket}_{L^\infty(\mathcal{C}_0)}\right]\norm{\mu(\rho_0,c_0)-\widetilde\mu}_{\cC^\alpha_{pw,\mathcal{C}_0}(\R^d)}\\
    +\left( \mathfrak{P}_{\mathcal{C}_0} +\ell^{-\alpha}_{\vph_0}\right)\left[\norm{\llbracket\mu(\rho_0,c_0)\rrbracket}_{L^\infty(\mathcal{C}_0)}+\norm{\llbracket\mu(\rho_0,c_0)\rrbracket,\, \llbracket \lambda(\rho_0,c_0)\rrbracket}_{L^\infty(\mathcal{C}_0)}\left\|1-\dfrac{\widetilde\mu}{\scalar{\mu(\rho_0,c_0)}}\right\|_{L^\infty(\mathcal{C}_0)}\right]\leqslant \dfrac{[\mu]}{4},\label{ep3.2}
    \end{multline}
    then, there exist a  time $T>0$ and a unique solution $(c,\rho, u)$ of the system \eqref{ep1.1} verifying the following:
    \begin{enumerate}
        \item $\forall\; 0< t\leqslant T,\;\;c(t)=\mathbb{1}_{D(t)^c}$\;\; with\;\; $\displaystyle\partial D(t)=:\mathcal{C}(t)$\, and with local parametrization  $\gamma(t)\in \cC^{1+\alpha}$;
        \item\label{epq10}  $P(\rho,c)-\widetilde P\in \cC([0,T], L^2(\R^d)\cap \cC^{\alpha}_{pw,\mathcal{C}}(\R^d))$; 
        \item\label{epq7} $u\in \cC([0,T], H^1(\R^d))\cap L^\infty((0,T), \dot W^{1,6}(\R^d))\cap L^{16}((0,T), \dot W^{1,8}(\R^d))$,\, $\sigma^{r/4}\nabla u\in L^4((0,T), \cC^\alpha_{pw,\mathcal{C}}(\R^d))$ for 
        \begin{gather}\label{c4.38}
            r=\begin{cases}
                \max\left(1/3,\,2\alpha\right)\quad & \text{if}\quad d=2,\\
                \max\left(\dfrac{9}{14},\,2\alpha\right) \quad & \text{if}\quad d=3,\; \alpha\in (0,1/2),\\
                \dfrac{2}{3}(5\alpha-1) \quad &  \text{if }\quad d=3,\; \alpha\in [1/2,1).
            \end{cases}
        \end{gather}
        \item \label{epq8}  $\dot u\in \cC([0,T], L^2(\R^d))\cap L^2((0,T),\dot H^1(\R^d))$,\, $\sqrt{\sigma}\nabla \dot u\in L^\infty((0,T), L^2(\R^d))$,\; $\sigma^{\tfrac{d}{4}}\nabla \dot u\in L^{8/d}((0,T), L^4(\R^d))$;
        \item  \label{epq11}$\sqrt{\sigma} \ddot u\in L^2((0,T)\times\R^d)$, $\sigma\ddot u\in L^\infty((0,T), L^2(\R^d))\cap L^2((0,T), \dot H^1(\R^d))$.
    \end{enumerate}
\end{theo}

To prove  \cref{thlocal}, we transform  the system \eqref{ep1.1} into Lagrangian coordinates, see \eqref{ep2.7} below, to fix the domain. The existence of a unique solution $(\overline{c}, \overline{\rho}, \overline{u})$ of \eqref{ep2.7} is obtained in  \cref{remaL} below. The flow of $\overline{u}$ is regular enough to ensure that, back in Eulerian coordinates, $(c,\rho,u)$ possesses all the properties stated  in \cref{thlocal}. 

\begin{rema}
\noindent
\begin{enumerate}
    \item  The smallness assumption \eqref{ep3.2} arises from the study of the linearized equations associated with the equations resulting from the transformation of  \eqref{ep1.1} into Lagrangian coordinates:
    \[
    \rho_0 v_t-\dvg (\mu \D v+ \lambda \dvg v I_d)= \dvg F.
    \]
    To obtain $L^p(\R^d)$ and piecewise H\"older estimates for $\nabla u$, we express (see \eqref{ep1.9}):
        \begin{align}
           \nabla v  &=\frac{1}{\widetilde\mu}(-\Delta)^{-1}\nabla\mathcal{P}
\dvg F   -\nabla^2 (-\Delta)^{-1}\left[\dfrac{1}{2\widetilde\mu+\lambda}(-\Delta)^{-1}\dvg \dvg F\right]\nonumber\\
&-\frac{1}{\widetilde\mu}(-\Delta)^{-1}\nabla\mathcal{P}\left(  \rho_{0}v_t\right) +\nabla^2 (-\Delta)^{-1}\left[\dfrac{1}{2\widetilde\mu+\lambda}(-\Delta)^{-1}\dvg(\rho_0 v_t)\right]\nonumber\\
&+\frac{1}{\widetilde\mu}(-\Delta)^{-1}\nabla \mathcal{P}
\dvg\{2(\mu-\widetilde\mu)\D  v\} -\nabla^2 (-\Delta)^{-1}\left[\dfrac{1}{2\widetilde\mu+\lambda}(-\Delta)^{-1}\dvg \dvg \{2(\mu-\widetilde\mu)\D v\}\right].\label{ep3.3}
\end{align}
where $\mathcal{P}$ is the Leray projector onto the space of divergence-free vector fields. Under  the smallness condition \eqref{ep3.2}, 
the last two terms of  \eqref{ep3.3} are small in comparison to the left-hand side.

\item   One of the main difficulties is to establish an estimate for the $L^\infty(\R^d)$-norm for $\nabla v$, (see page \pageref{pagecCalpha}). Except for the third term in the expression \eqref{ep3.3}, which is regular in space, all other terms are second-order or fourth-order Riesz transforms of discontinuous functions (see the regularity of $\mu$, $\lambda$ and $F$ in \eqref{ep1.2}-\eqref{c7.6}). It is well known that Riesz operators are not continuous on  $L^\infty(\R^d)$, which makes
the estimation of the $L^\infty(\R^d)$-norm of $\nabla u$ more challenging.  We employ a quantitative piecewise H\"older estimate for even-order Riesz operators obtained by Gancedo and Garc\'ia-Ju\'arez (see \cref{lemA2} or \cite{gancedo2021global,gancedo2021quantitative}) to derive a piecewise H\"older estimate for $\nabla v$. Roughly speaking, the authors proved that  these operators are continuous on $L^q(\R^d)\cap \cC^{\alpha}_{pw,\mathcal{C}_0}(\R^d)$, $q\in [1,\infty)$, with a norm that may depend on the regularity of $\mathcal{C}_0$ (see \cref{notabene}-item 1).

\end{enumerate}
\end{rema}
\begin{rema}\label[rema]{compa}
    Compared to Tani's work \cite{tani1984two}, we only require the interface $\mathcal{C}$ to be $\cC^{1+\alpha}$-regular   
  and the density to be piecewise $\cC^\alpha$-H\"older continuous.  Furthermore, initially non-Lipschitz velocities are allowed.
  Indeed, for constant viscosity coefficients, the compatibility condition \eqref{c4.1} reads:  
    \[
    \widetilde\mu\Delta u_0+(\widetilde\mu+\widetilde\lambda)\nabla\dvg u_0=\nabla P(\rho_0)+\dvg (\Pi)_{|t=0},
    \]
    which leads to (similar to the derivation of \eqref{ep3.3}): 
\begin{align*}
\nabla u_0&=-\dfrac{1}{\widetilde \mu}\nabla(-\Delta)^{-1}\mathcal{P} \dvg (\Pi)_{|t=0}-\dfrac{1}{2\widetilde\mu+\widetilde\lambda}\nabla (-\Delta)^{-1}\big(\mathcal{I}-\mathcal{P}\big) \dvg (\Pi)_{|t=0}\\
&-\dfrac{1}{2\widetilde\mu+\widetilde\lambda}\nabla (-\Delta)^{-1}\nabla (P(\rho_0,c_0)-\widetilde P).
\end{align*}
  We observe that, in general, the first two terms are  not bounded for $\dvg (\Pi)_{|t=0}\in L^2(\R^d)$.
\end{rema}
 \begin{rema}\label[rema]{remalocal}
 The solution $(c,\rho, u)$ constructed in the above theorem exhibits the following properties.
     \begin{enumerate}
         \item There exists a constant $R$ depending on the lifespan $T$, on the norms of the initial data  and 
         on the initial interface $\mathcal{C}_0$ such that: 
         \[
         \sup_{[0,T]} \left\{\norm{u}_{H^1(\R^d)}+ \norm{\dot u}_{L^2(\R^d)}\right\}\leqslant R.
         \]
         \item  From \cref{epq8}, we have $\rho \dot u\in \cC([0,T], L^2(\R^d))$
         and from the momentum equation \eqref{ep1.1} we infer that  $\dvg (\Pi)\in \cC([0,T], L^2(\R^d))$:
          the compatibility condition \eqref{epq1} is propagated in time. 
         \item\label{epq44} From \cref{epq7} and interpolation we infer that the velocity $u$ is continuous in time and space. Additionally, from  \cref{epq8}, its material derivative, $\dot u$ is continuous in space. Furthermore, from \cref{epq11} we deduce that $\ddot u$ is at least continuous across the interface 
         $\mathcal{C}(t)$. Roughly, for all $j\in \llbracket 1,J\rrbracket$, $s\in V_j$, we have:
         \[
         \llbracket \dot u(t,\gamma_j(t,s))\rrbracket=0\implies \dfrac{d}{dt} \llbracket \dot u(t,\gamma_j(t,s))\rrbracket=\llbracket \ddot u(t,\gamma_j(t,s)\rrbracket=0.
         \]
         \item From \cref{epq7}, we infer that $\nabla u\in L^1((0,T), \cC^\alpha_{pw,\mathcal{C}}(\R^d))$, and hence we have the  persistence of the initial interface's regularity. Indeed  $\mathcal C(t)$ is $\cC^{1+\alpha}$-regular manifold of dimension $d-1$, with the following estimates:
          \begin{gather}\label{epq12}
          \begin{cases}
           \displaystyle   \abs*{\mathcal{C}(t)}_*&\displaystyle\leqslant \abs*{\mathcal{C}_0}_* \exp\left(\int_0^t\norm{\nabla u(\tau)}_{L^\infty(\R^d)}d\tau \right),\\
          \displaystyle\norm{\nabla\gamma_j(t)}_{\cC^\alpha}&\displaystyle\leqslant \left(\norm{\nabla\gamma_{0,j}}_{\cC^\alpha} +\norm{\nabla\gamma_{0,j}}_{L^\infty}^{1+\alpha}\int_0^t \norm{\nabla u(\tau)}_{\dot\cC^\alpha_{pw,\mathcal{C}(\tau)}(\R^d)}d\tau\right)\\
          &\times \exp\left((2+\alpha)\int_0^t \norm{\nabla u(\tau)}_{L^\infty(\R^d)}d\tau\right).
          \end{cases}
          \end{gather}  
          \item  Applying $\dpt\cdot\,+ \dvg (\, \cdot\; u)$ to the momentum equations $\eqref{ep1.1}_3$, we  find:
          \begin{gather}\label{epq13}
          \rho \ddot u=\dvg (2\mu(\rho,c)\D \dot u)+\nabla(\lambda(\rho,c)\dvg \dot u)+ \dvg \mathcal{F},
          \end{gather}
          where $\mathcal{F}$ is a lower order term that belongs, at least, to $L^2 ((\tau,T),\;\cC^\alpha_{pw,\mathcal{C}}(\R^d))$, for all $\tau\in (0,T)$.  With the help of \cref{epq11}, small fluctuation assumption on $\mu(\rho,c)$, quantitative H\"older estimate  \cref{lemA2}, we infer from  \eqref{epq13} that:
          \[
          \nabla \dot u\in L^2 ((\tau,T),\; \cC^\beta_{pw,\mathcal{C}}(\R^d))\quad\text{with}\quad \beta=\begin{cases}
              \alpha,\quad\text{in}\quad d=2,\\
              \min\{\alpha,\,1/2\},\quad\text{in}\quad d=3,
          \end{cases}
          \]
          for all $\tau\in (0,T)$.
       In particular, $\nabla \dot u$ is also continuous on both sides of $\mathcal{C}$.
     \end{enumerate}
 \end{rema}
 \cref{thlocal} immediately implies:
\begin{coro}[Blow-up criterion]\label[theo]{thblowup}
    Let $(c, \rho, u)$ be the solution constructed in \cref{thlocal} defined up to a maximal time $T^*$. If  $T^*<\infty$, then we have either:
    \begin{align*}
     \limsup_{t\to T^*}& \left\{\abs*{\mathcal{C}(t)}_*+\norm{\mathcal{C}(t)}_{\text{Lip}}+\norm{\mathcal{C}(t)}_{\dot \cC^{1+\alpha}}+\bigg\|\dfrac{1}{\rho(t)},\;\;\dfrac{1}{\mu(\rho,c)(t)}\bigg\|_{L^\infty(\R^d)}\right\}\\
    +\limsup_{t\to T^*}&\left\{\norm{u(t)}_{H^1(\R^d)}+\norm{(\rho \dot u)(t)}_{L^2(\R^d)}+
    \norm{P(\rho,c)(t)-\widetilde P}_{\cC^\alpha_{pw,\mathcal{C}(t)}(\R^d)}\right\}=\infty,
    \end{align*}
    or
    \begin{multline}
    \limsup_{t\to T^*}\left[1+\norm{\lambda (\rho,c)(t)}_{\dot \cC^\alpha_{pw,\mathcal{C}(t)}(\R^d)}+\left( \mathfrak{P}_{\mathcal{C}(t)} +\ell^{-\alpha}_{\vph(t)}\right)\norm*{\llbracket\lambda(\rho,c)(t)\rrbracket}_{L^\infty(\mathcal{C}(t))}\right]\norm{\mu(\rho,c)(t)-\widetilde\mu}_{\cC^\alpha_{pw,\mathcal{C}(t)}(\R^d)}\\
    +\limsup_{t\to T^*}\left( \mathfrak{P}_{\mathcal{C}(t)} +\ell^{-\alpha}_{\vph(t)}\right)\left[\norm*{\llbracket\mu(\rho,c)(t)\rrbracket}_{L^\infty(\mathcal{C}(t))}+\norm*{\llbracket\mu(\rho,c)(t)\rrbracket, \llbracket \lambda(\rho,c)(t)\rrbracket}_{L^\infty(\mathcal{C}(t))}\left\|1-\dfrac{\widetilde\mu}{\scalar{\mu(\rho,c)(t)}}\right\|_{L^\infty(\mathcal{C}(t))}\right]\geqslant [\mu].\label{c3.38}
    \end{multline}
\end{coro}
\begin{rema}
    Our result can be seen as a first step towards the existence of global solutions with small data which we plan to attack in a forthcoming paper. The construction of discontinuous weak-solutions for compressible models with variable viscosity cannot be achieved in the classical way by regularizing the initial data and applying some classical results \emph{à la Matsmura-Nishida} since this procedure \emph{destroys} the discontinuities in the density. Some difficult constructions  are then needed in order to obtain discontinuous weak-solutions, see \cite[Section 3]{hoff2002dynamics}. We believe that since we dispose of powerful functional analysis tools, see \cref{app1}, an approach where local-well posedness is first proved to be quite natural. With \cref{thlocal} in hand, in order to attack the global existence issue, we do not have to worry anymore about constructing appropriate approximate systems that preserve the main feature of the system \eqref{ep1.1}.
\end{rema}
\begin{rema}
     The difficulty in proving global existence for the two-fluid model lies in the fact that we are not able to estimate the $L^1((0,T), \cC^\alpha_{pw,\mathcal{C}}(\R^d))$
    norm  of $\nabla u$ uniformly with respect to time. In particular, see \eqref{epq12},  the characteristic $\abs*{\mathcal{C}(t)}_*$ and the $\cC^{1+\alpha}$-norm of the interface $\gamma_j(t)$ are expected to grow exponentially with respect to time which prevents the validity of the condition \eqref{c3.38} globally in time.

    In the study of the Navier-Stokes equations\textit{ with constant viscosity for $d=2$}, Hoff \cite{hoff2002dynamics} has observed that the exponential-in-time growth of these quantities can be compensated by the exponential decay of the pressure jump across the interface (in our case we would be needing the exponential-in-time decay of the viscosity jump, see the explosion criterion \eqref{c3.38}).  In order to obtain exponential decay, we use the continuity of the normal component of the stress tensor, we follow the arguments starting from \eqref{c3.32} below to establish:
    \[
    \llbracket (2\mu(\rho,c)+\lambda(\rho,c))\dvg u- P(\rho,c)\rrbracket=  2\llbracket \mu(\rho,c)\rrbracket \left(\scalar{\dvg u}-\scalar{\D^{jk} u}n^{j}_xn^k_x\right).
    \]
    Next, we observe that $f=f(\rho,c)$ defined  up to a constant as 
    \[
   f(\rho,c)=\int^\rho \dfrac{2\mu(s,c)+\lambda(s,c)}{s}ds
   \]
   satisfies: 
   \begin{gather*}
       \dpt f(\rho,c)+ u\cdot \nabla f(\rho,c)+ P(\rho,c)-\widetilde P=-\left[(2\mu(\rho,c)+\lambda(\rho,c))\dvg u- P(\rho,c)\right].
   \end{gather*}
    We write the above equation along the flow $X$ of the velocity and then take jumps to obtain:
   \begin{gather}
   \dfrac{d}{dt}\llbracket f(\rho,c)\rrbracket(t,X(t,s))+  \llbracket P(\rho,c)\rrbracket(t,X(t,s)) =2\llbracket \mu(\rho,c)\rrbracket \left(\scalar{\dvg u}-\scalar{\D^{jk} u}n^{j}_xn^k_x\right)(t,X(t,s)),\; s\in \mathcal{C}_0.
   \end{gather}
   Observe that in the constant viscosity case, the RHS of the above equation vanishes and using the fact that the pressure is strictly increasing one infers in a straightforward manner the exponential-in-time decay of the jump of $f$. 
    \emph{For the two-fluid model, when the pressure and viscosity laws of the two phases are different, it is less obvious how to obtain exponential decay of the jump or even if this property holds in general}. First of all the RHS of the above identity is no longer $0$ and second of all in order to obtain exponential-in-time decay, the jump of $f,P$ and $\mu$ need to be "comparable", (this is the case for constant $\mu$ and increasing $P=P(\rho)$). For multi-phase flows the minimal requirements seem to be that:
    \begin{gather}\label{c3.54}
    0<\text{\textnu}_*\leqslant \dfrac{P(\rho,1)-P_-(\rho',0)}{f(\rho,1)-f(\rho',0)}\leqslant \text{\textnu}^* \quad \text{ and }\quad \dfrac{\abs{\mu(\rho,1)-\mu(\rho',0)}}{\abs{f(\rho,1)-f(\rho',0)}}+\dfrac{\abs{\lambda(\rho,1)-\lambda(\rho',0)}}{\abs{f(\rho,1)-f(\rho',0)}}\leqslant \text{\textnu}^*,
    \end{gather}
  at least, for all $\rho$ close to $\widetilde\rho(1)$ and $\rho'$ close to $\widetilde\rho(0)$ where $\widetilde P=P(\widetilde \rho(1),1)=P(\widetilde\rho(0),0)$.  This requirement does not seem to occur in general. For instance, it does not hold in the case of two different constant viscosity coefficients say $\mu(\rho,1)=\widetilde\mu,\; \mu(\rho,0)=\mu_-$, where $\mu_-\neq \widetilde\mu$. 
\end{rema}

\subsection*{\emph{Outline of the paper}}
The remainder of this paper is  divided into  three sections. In the first one, \cref{lagrangian}, we rewrite system \eqref{ep1.1} in Lagrangian coordinates, see \eqref{ep2.7} below, and in the second one, \cref{linearised}, we study the linearized system associated to \eqref{ep2.7}. Finally, the study of the full nonlinear system in Lagrangian coordinates \eqref{ep2.7} is presented in \cref{nonlinearised} from which \cref{thlocal} follows.
\section{Change of variables in Lagrangian coordinates}\label[section]{lagrangian}
The aim of this section is to write the system \eqref{ep1.1}
in Lagrangian coordinates. We know that for a Lipschitz vector field $u$, one can define 
a flow map $\mathscr{X}$ given by the Duhamel formula:
\[
\mathscr{X}(t,y)=y+\int_0^t u(s, \X(s,y))ds 
\]
and adopting the following notation: for all function $g=g(t,x)$ we denote by $\overline g$ the function:
\[
\overline g(t,y):= g(t,\X(t,y))
\]
then, the flow of the vector field $u$ can be written in term of $\overline u$ by :
\[
\X (t,y)=:\X_{\overline{u}}(t,y)= y+ \int_0^t \overline u(s,y)ds.
\]
And for any $g=g(t,x)$, one has that 
\begin{gather}\label{ep2.6}
\partial_t \overline g(t,y)= (\dpt g + u\cdot\nabla g)(t,\X_{\overline{u}}(t,y)).
\end{gather}
We introduce some notations.
\begin{nota}
Let us fix a square matrix $B=(b_{ij})_{ij}\in \mathcal{M}_{n}(\R)$, $n\in \N^*$.
    \begin{itemize}
    \item  The trace of $B$ is the sum of the diagonal elements of $B$:
    \[
    \text{Tr}(B)= b_{ii}.
    \]
     The transposed matrix of $B$, denoted  $B^T$,  is:
    \[
    B^T= (b_{ji})_{ij}
    \]
    and for  $B'=(b'_{ij})_{ij}\in \mathcal{M}_{n}(\R)$, we have:
    \[
    \text{Tr}(B\cdot B')= b_{ij}b'_{ji}=B\colon B'= (B')^T\colon B^T.
    \]
    \item For $i,j\in \llbracket 1,n\rrbracket$, $M_{ij}(B)$ denotes the determinant of the $(n-1)\times (n-1)$ matrix 
    obtained by removing the $i$-th row and $j$-th column of $B$.
    \item The adjugate matrix of $B$, denoted as $\adj (B)$, is given by:
    \[
    (\adj (B))_{ij}= (-1)^{i+j} M_{ji}(B).
    \]
    \item  If $\det (B)\neq 0$, then $B$ is invertible and we have:
    \[
    B^{-1}=\dfrac{1}{\det (B)} \adj (B).
    \]
    \end{itemize}
\end{nota}
If $D G$ denotes the Jacobian matrix of $G$, then chain rule yields:
\[
D \overline G = D G \circ \X_{\overline{u}} \cdot D \X_{\overline{u}} ,
\]
but, since  
\begin{gather}\label{ep2.4}
(D \X_{\overline{u}} (t,y))_{ij}= \delta^{ij}+\int_0^t\partial_j \overline u^ i(s,y)ds,
\end{gather}
it holds, at least for small time, that the Jacobian matrix of the change of variable is invertible and hence:
\[
\overline {D G}(t,y)= D \overline G(t,y)\cdot (D\X_{\overline{u}}(t,y))^{-1}= \dfrac{1}{\det D\X_{\overline{u}} (t,y)}D \overline G(t,y)\cdot \adj D \X_{\overline{u}}(t,y).
\]
In particular, one has that:
\begin{gather}\label{ep2.2}
\overline{\partial_j G^ i }=\dfrac{1}{\det D\X_{\overline{u}}}\partial_k \overline G^i\left(\adj D \X_{\overline{u}}\right)_{kj}.
\end{gather}
On the other hand, from  Piola's formula: for any matrix $A\in \mathcal{M}_{n\times n}(\R)$, 
\[
\partial_k (\adj (A))_{kj}=0
\]
we infer that:
\begin{gather}\label{ep2.3}
   \overline{\partial_j G^ i }=\dfrac{1}{\det D\X_{\overline{u}} }\partial_k \left[\overline G^i\left(\adj D \X_{\overline{u}}\right)_{kj}\right].
\end{gather}
The first form \eqref{ep2.2} will be referred as the non-conservative form and the second \eqref{ep2.3} is the conservative one and both will
be used. From the non-conservative form, we find  
\[
\overline{ \dvg G}= \dfrac{1}{\det D\X_{\overline{u}}}\text{Tr}(D \overline G\cdot \adj D \X_{\overline{u}}), 
\]
and by setting 
\[
A_{\overline{u}}(t,y)= (D \X_{\overline{u}}(t,y))^{-1}
\]
we have
\begin{gather}\label{ep2.5}
\dvg_{A_{\overline{u}}} G =\overline{ \dvg G}= D\overline G : A_{\overline{u}} = A^T_{\overline{u}}:\nabla \overline G.
\end{gather}

On the other hand,  for all matrix $B \in \cC^1( [0,T], \mathcal{M}_{ n}(\R))$, one has: 
\begin{align*}
\dfrac{d \det B}{ds}(s)&= \dfrac{\partial \det B}{\partial b_{ij}} \dfrac{\partial b_{ij}}{\partial s} (s)\\
&=(-1)^{i+j} M_{ij}(B(s)) \dfrac{\partial b_{ij}}{\partial s} (s)\\
&= \adj(B(s))_{ji}\left(\dfrac{\partial B}{\partial s} (s)\right)_{ij}\\
&=   \dfrac{\partial B(s)}{\partial s}:\adj(B(s)).
\end{align*}
In particular, for $B=D\X_{\overline{u}}$, and remembering from \eqref{ep2.4} that 
\[
\dfrac{\partial D\X_{\overline{u}}}{\partial t}= D \overline u
\]
we find 
\[
\dfrac{\partial \det D\X_{\overline{u}}}{\partial t}(t,y)=D \overline u (t,y): \adj D\X_{\overline{u}} (t,y)).
\]
Hence, the identity \eqref{ep2.5} yields:
\begin{gather}\label{eq3.4}
\dfrac{\partial \det D \X_{\overline{u}}}{\partial t} (t,y)=\det D \X_{\overline{u}}(t,y) \dvg_{A_{\overline{u}}} u(t,y)  .
\end{gather}
Finally, combining \eqref{ep2.6} and the previous identity, we notice that the continuity equation on the density 
reads:
\[
\dfrac{\partial (\overline \rho \det D \X_{\overline{u}})}{\partial t} (t,y)=0
\]
and consequently 
\[
\overline \rho(t,y)=\dfrac{\rho_0(y)}{\det D\X_{\overline{u}} (t,y)}.
\]
Now, let us proceed to rewrite the momentum equations  in Lagrangian coordinates. First, from what precedes, it is clear that 
\[
\left[\dpt(\rho u^j)+\dvg(\rho u^j u)\right](t,\X_{\overline{u}}(t,y))=\dfrac{1}{\det D\X_{\overline{u}}(t,y)}\rho_0(y)\dpt \overline u(t,y)
\]
so it only remains to compute the term involving the stress tensor. From the conservative form \eqref{ep2.3},
we have:
\begin{align*}
\overline{\dvg \Pi^j}&=\dfrac{1}{\det D\X_{\overline{u}}(t,y)}\partial_k\left\{(\adj D \X_{\overline{u}})_{kl}  \overline {\Pi}^j_l \right\} \\
&=\dfrac{1}{\det D\X_{\overline{u}}(t,y)}\partial_k\left\{(\adj D \X_{\overline{u}})_{kl} \left\{ 2\mu(\overline\rho,c_0)\D_{A_{\overline{u}}}^{jl}\overline u +\left(\lambda(\overline\rho,c_0)\dvg_{A_{\overline{u}}} \overline u-P(\overline\rho,c_0)+\widetilde P\right)\delta_{jl}\right\}\right\}
\end{align*}
where, we use the following notation: 
\[
2\D_B v= D v\cdot B+B^T\cdot \nabla v.
\]
This ends the transformation of the momentum equation in Lagrangian coordinates.
Summing up, we find that $(\overline c, \overline \rho, \overline u)$ solves:
\begin{gather}\label{ep2.7}
    \begin{cases}
    \dpt \overline{c}=0,\\
    \dpt \left\{\overline \rho J_{\overline{u}} \right\}=0,\\
    \rho_0\dpt \overline u=\dvg \left\{ \adj(D\X_{\overline{u}})\left\{2\mu(\overline\rho,c_0)\D_{A_{\overline{u}}} \overline u+\left(\lambda(\overline\rho,c_0)\dvg_{A_{\overline{u}}} \overline u-P(\overline\rho,c_0)+\widetilde P\right) I_d\right\}\right\},
    \end{cases}
\end{gather}
where $J_{\overline{u}}:=\abs{\det D \X_{\overline{u}}}$. We establish the well-posedness of \eqref{ep2.7} in  the following theorem. 
\begin{theo}\label[theo]{remaL}
    Assume the hypotheses in \cref{thlocal} hold true for $(c_0,\rho_0,u_0)$, the pressure law,  and  the viscosity laws $\mu$, $\lambda$. Additionally, suppose that the small fluctuation assumption \eqref{ep3.2} also holds true, and consider $r$ as given in \eqref{c4.38}. 
    
    Then, there exist $T>0$ and a unique solution $\overline{u}$ of the Cauchy problem 
    \begin{gather}\label{c4.39}
    \begin{cases}
        \rho_0\dpt \overline u=\dvg \left\{ \adj(D\X_{\overline{u}})\left\{2\mu(\rho_0J_{\overline{u}}^{-1} ,c_0)\D_{A_{\overline{u}}} \overline u+\left(\lambda(\rho_0J_{\overline{u}}^{-1},c_0)\dvg_{A_{\overline{u}}} \overline u-P(\rho_0J_{\overline{u}}^{-1},c_0)+\widetilde P\right) I_d\right\}\right\},\\
        \overline{u}_{|t=0}= u_0,
    \end{cases}
    \end{gather}
    verifying the following:
    \begin{enumerate}
        \item $\overline u\in \cC([0,T], H^1(\R^d))\cap L^\infty((0,T), \dot W^{1,6}(\R^d))\cap L^{16}((0,T), \dot W^{1,8}(\R^d))$,\, $\sigma^{r/4}\nabla \overline  u\in L^4((0,T), \cC^\alpha_{pw,\mathcal{C}_0}(\R^d))$;
        \item   $\dpt \overline  u\in \cC([0,T], L^2(\R^d))\cap L^2((0,T),\dot H^1(\R^d))$,\, $\sqrt{\sigma}\nabla \dpt \overline u\in L^\infty((0,T), L^2(\R^d))$,\, $\sigma^{\tfrac{d}{4}}\nabla \dpt \overline  u\in L^{8/d}((0,T), L^4(\R^d))$;
        \item  $\sqrt{\sigma} \partial_{tt} \overline  u\in L^2((0,T)\times\R^d)$, $\sigma\partial_{tt}\overline   u\in L^\infty((0,T), L^2(\R^d))\cap L^2((0,T), \dot H^1(\R^d))$.
    \end{enumerate}
\end{theo}

In \cref{linearised} below, we study the linearized equations associated with \eqref{c4.39}. The proof of the existence of $\overline{u}$ is the purpose of \cref{nonlinearised}, and it relies on a fixed point argument. The proof of uniqueness follows a similar argument to that in \cite{danchin2020well} and is therefore omitted here.  It only requires viscosity far from vacuum and the velocity gradient to meet the following conditions:
\[
\nabla \overline u\in L^1((0,T), L^\infty(\R^d)),\quad \text{and}\quad \sqrt\sigma\nabla \overline u\in L^2((0,T), L^\infty(\R^d)).
\]
In particular, no small fluctuation assumption on the viscosity $\mu$ is required. Once we have constructed $\overline{u}$, 
we express $(\overline{c},\overline{\rho})$ from $\eqref{ep2.7}_{1,2}$, and we finally go back into Eulerian coordinates to obtain that $(c,\rho,u)$ solves \eqref{ep1.1} with initial data $(c_0,\rho_0,u_0)$. 
We now proceed to the study of the linearized equations associated with \eqref{c4.39}.

\section{Study of the linearised system}\label[section]{linearised}
\subsection{Statement of the main result}
The aim of this section is to establish the well-posedness of the linear non-homogeneous parabolic equation 
\begin{gather}\label{ep2}
    \rho_0 u_t -\dvg(2\mu\D  u+\lambda\dvg  u I_d) = \dvg F
\end{gather}
associated to \eqref{ep2.7}. The initial density $\rho_0=\rho_0(x)$ and the viscosity $\mu=\mu(t,x), \;\lambda=\lambda(t,x)$, are supposed to verify: 
\begin{gather}\label{ep1.2}
\begin{cases}
    0<\rho_*\leqslant \rho_0\leqslant \rho^*,\;0<\mu_*\leqslant \mu\leqslant \mu^*,\;\; 0\leqslant \lambda\leqslant \mu^*;\\
     \mu,\;\lambda\in L^\infty((0,T), \cC^\alpha_{pw,\mathcal{C}_0}(\R^d));\\
      \mu_t,\;\lambda_t\in L^\infty((0,T), L^4(\R^d))\cap L^1((0,T), L^\infty(\R^d))\cap L^{16}((0,T), L^8(\R^d));\\
      \sqrt\sigma (\mu_t,\, \lambda_t)\in L^2((0,T), L^\infty(\R^d)),\;
       \sigma(\mu_{tt},\lambda_{tt})\in L^{2}((0,T), L^4(\R^d)),
\end{cases}
\end{gather}
for positive constants $\rho_*, \rho^*, \mu_*, \mu^*$ and time $0<T<\infty$. 
The source term $F$ in \eqref{ep2} and the initial datum $u_0$ are such that  $(F,u_0)\in Y_T(\mu_0,\lambda_0)$ where:
\begin{align}
 Y_T(\mu_{0},\lambda_{0}):=\left\{(F,u_0)\right.&\left.\in \cC([0,T],L^2(\R^d))\times H^1(\R^d)
         \colon F\in  L^\infty((0,T),  L^6(\R^d)),\right.\notag\\
         &\left.F\in  L^{16}((0,T), L^8(\R^d)),\;\sigma^{r/4} F\in L^4((0,T),\, \cC^\alpha_{pw,\mathcal{C}_0}(\R^d)),\right.\notag\\
         &\left.
         F_t\in L^2((0,T)\times\R^d),\;\sigma^{\tfrac{d}{4}} F_t\in L^{8/d}((0,T), L^4( \R^d)),\right.\notag\\
         &\left.\sigma F_{tt}\in L^2((0,T)\times\R^d),\; \dvg (2\mu_{0}\D u_0+\lambda_{0}\dvg u_0 I_d+ F_{|t=0})\in L^2(\R^d)\right\}.\label{c7.6}
\end{align}
Above, $\mu_0:=\mu_{|t=0}$, $\lambda_0:=\lambda_{|t=0}$ and $0<r<8/3$  depends on $0<\alpha<1$ and the dimension, see \eqref{c5.1} below. The existence and uniqueness of  solution to \eqref{ep2} is established in the following space:
\begin{align}
X_T:=&\left\{u\in \cC([0,T], H^1(\R^d))\colon \right.\nabla u\in L^\infty((0,T),L^6(\R^d))\cap  L^{16}((0,T), L^8(\R^d)),\,\sigma^{r/4} \nabla u\in L^4((0,T),\, \cC^\alpha_{pw,\mathcal{C}_0}(\R^d)),\nonumber\\
&\left.u_t\in \cC([0,T], L^2(\R^d))\cap L^2((0,T), \dot H^1(\R^d)),\;  \sqrt{\sigma}\nabla u_t \in L^\infty((0,T), L^2(\R^d)),\; \sigma^{\tfrac{d}{4}}\nabla u_t\in L^{8/d}((0,T), L^4(\R^d)),\right.\nonumber\\
&\left.\sqrt{\sigma}  u_{tt}\in L^2((0,T)\times\R^d),\;\sigma u_{tt} \in L^\infty((0,T), L^2(\R^d))\cap L^2((0,T), \dot H^1(\R^d))\right\}.\label{ep1.3}
\end{align}
 The time weights are inspired by the pioneering work of Hoff \cite{hoff1995global}.  One of the main difficulties is to obtain estimates for  $\nabla u$, the solution of \eqref{ep2}, that allow us to propagate the regularity of the interface $\mathcal{C}_0$. The difficulty arises from the fact that the density and viscosity do not have any (global-in-space) Sobolev regularity, and therefore we are unable to differentiate \eqref{ep2} with respect to space.

The main result of this section is stated in the following theorem:
 
\begin{theo}\label[theo]{Theorem 3.1}
    Let $T\in (0,\infty)$, the viscosity $\mu,\, \lambda$ and the initial density $\rho_0$ verify \eqref{ep1.2}, $0<\alpha<1$.  There exists a constant $[\mu]>0$ depending only on the dimension $d$, $\alpha$ and $\widetilde\mu$  such that if:
       \begin{align}
    \sup_{[0,T]}\left[1+\norm{\lambda}_{\dot \cC^\alpha_{pw,\mathcal{C}_0}(\R^d)}\right.&\left.+\left( \mathfrak{P}_{\mathcal{C}_0} +\ell^{-\alpha}_{\vph_0}\right)\norm{\llbracket\lambda\rrbracket}_{L^\infty(\mathcal{C}_0)}\right]\sup_{[0,T]}\norm{\mu-\widetilde\mu}_{\cC^\alpha_{pw,\mathcal{C}_0}(\R^d)}\nonumber\\
    &+\left( \mathfrak{P}_{\mathcal{C}_0} +\ell^{-\alpha}_{\vph_0}\right)\sup_{[0,T]}\left[\norm{\llbracket\mu\rrbracket}_{L^\infty(\mathcal{C}_0)}+\norm{\llbracket\mu\rrbracket,\, \llbracket \lambda\rrbracket}_{L^\infty(\mathcal{C}_0)}\left\|\dfrac{\scalar{\mu-\widetilde\mu}}{\scalar{\mu}}\right\|_{L^\infty(\mathcal{C}_0)}\right]\leqslant [\mu],\label{eq3.64}
    \end{align}
    then, for all $(F,u_0)\in Y_T(\mu_0,\lambda_0)$,  there exists a unique solution   of the Cauchy problem: 
    \begin{gather}\label{ep1.6}
        \begin{cases}
            \rho_0 u_t-\dvg(\mu \D u)-\nabla (\lambda \dvg u)=\dvg F,\\
            u_{|t=0}= u_0,
        \end{cases}
    \end{gather}
    in $X_T$. Moreover,  the following estimates hold true with a constant $C_*$ that depends polynomially on the upper and lower bounds of the density and viscosity (see first line of \eqref{ep1.2}):
  \begin{itemize}
      \item  \textbf{Estimates for $u$ and $u_t$}:
  \begin{align}
     \sup_{[0,T]}\norm{u,\,u_t,\, \nabla u}_{L^2(\R^d)}^{2}&+\int_{0}^{T}
\norm{\nabla u,\, \nabla u_t,\, u_t}_{L^2(\R^d)}^{2}\leqslant C_* \left(\norm{(F,u_0)}_{Y_T(\mu_0,\lambda_0)}^2+\int_0^T\norm{\mu_t,\,\lambda_t}_{L^4(\R^d)}^2\norm{F}_{L^2\cap L^6(\R^d)}^2\right)\nonumber\\
&\times \exp\left( C_*\int_0^T\left[\norm{\mu_t,\,\lambda_t}_{L^\infty(\R^d)}
+\norm{\mu_t,\,\lambda_t}_{L^4(\R^d)}^2\right] \right).\label{ep1.5}
\end{align}
\item  \textbf{Higher integrability estimates for  $\nabla u$}:
\begin{align}\label{eq3.36}
\begin{cases}\vspace{0,1cm} \displaystyle
    \norm{\nabla u}_{L^\infty((0,T),L^6(\R^d))}\leqslant C_*\norm{F}_{L^\infty((0,T),L^6(\R^d))}+ C_*\sup_{[0,T]}\norm{F, u_t, \nabla u}_{L^2(\R^d)},\\
    \norm{\nabla u}_{L^{16}((0,T), L^8(\R^d))} \leqslant C_* \left(\norm{ F}_{L^{16}((0,T), L^8(\R^d))}+  \norm{F, u_t, \nabla u}_{L^2\cap L^\infty((0,T), L^2(\R^d))}+\norm{\nabla u_t}_{L^2((0,T)\times \R^d)}\right).
\end{cases}
\end{align}
\item   \textbf{Estimates for $u_{tt}$}:
\begin{align}
    \int_0^T\norm{ \sqrt\sigma u_{tt},\; \sigma \nabla u_{tt}}_{L^2(\R^d)}^2&+\norm{\sqrt\sigma \nabla u_t,\; \sigma u_{tt}}^2_{L^\infty((0,T),L^2(\R^d))}\leqslant C_* \left[ \int_0^{\sigma(T)}\norm{\nabla u_t}^2_{L^2(\R^d)}+\int_0^{T} \norm{F_t,\, \sigma F_{tt}}_{L^2(\R^d)}^2\right.\nonumber\\
    &\left.+  \int_0^T\norm{\mu_t,\,\lambda_t,\, \sigma(\mu_{tt},\,\lambda_{tt})}_{L^4(\R^d)}^2\norm{\nabla u}_{L^4(\R^d)}^2\right]\nonumber\\
    &\times\exp\left(C_*\int_0^T\left[\norm{\mu_t,\,\lambda_t}_{L^\infty(\R^d)}
    + \sigma \norm{\mu_t,\,\lambda_t}_{L^\infty(\R^d)}^2\right]\right).\label{eq3.34}
\end{align}
The norms of $u$ involved in this estimate are controlled thanks to \eqref{ep1.5}-\eqref{eq3.36}.
\item  \textbf{$L^{8/d}((0,T), L^4(\R^d))$ norm estimate for $\sigma^{\tfrac{d}{4}}\nabla u_t$}:
       \begin{align}
       \int_0^T\sigma^{2}\norm{\nabla u_t}_{L^4(\R^d)}^{8/d}&\leqslant C_*\int_0^T\sigma^{2}\norm{F_t}_{L^4(\R^d)}^{8/d}+C_*\int_0^T\sigma^{2}\norm{\nabla u}_{L^8(\R^d)}^{8/d}\norm{\mu_t,\,\lambda_t}_{L^8(\R^d)}^{8/d}\nonumber\\
       &+C_*\norm{\sigma u_{tt}}_{L^\infty((0,T),L^2(\R^d))}^2\int_0^T\norm{\nabla u_t,\, F_t,\,\lambda_t \dvg u,\, \mu_t\D u}_{L^2(\R^d)}^{\tfrac{8}{d}-2}.\label{eq3.48}
\end{align}
\item \textbf{H\"older estimate for $\nabla u$}: 
\begin{align}
 \int_0^T\sigma^{r}\norm{\nabla u}_{\cC^\alpha_{pw,\mathcal{C}_0}(\R^d)}^4 
         &\leqslant  C_{**}\int_0^T\sigma^{r}\norm{F}_{\cC^\alpha_{pw,\mathcal{C}_0}(\R^d)}^4 +C_{**}\norm{F,\,  \nabla u,\, u_t,\,\sqrt\sigma\nabla u_t}_{L^\infty((0,T),L^2(\R^d))}^4\nonumber\\
         &+C_{**}\left[\int_0^T \norm{F,\, \nabla u,\,u_t,\,  \nabla u_t}_{L^2(\R^d)}^2\right]^2+C_{**}\left[\int_0^T \sigma^{2}\norm{\nabla u_t}_{L^4(\R^d)}^{8/d}\right]^{d/2},\label{eq3.49}
\end{align}
with a constant $C_{**}$ that depends polynomially on the same parameters as $C_*$, and  on $\alpha$, $\ell_{\vph_0}^{-\tfrac{\alpha}{2}}$, $\mathfrak{P}_{\mathcal{C}_0}$ and $\norm{\mu,\;\lambda}_{\dot \cC^{\alpha}_{pw,\mathcal{C}_0}(\R^d)}$.
\end{itemize}
\end{theo}
\begin{rema}
The first estimate \eqref{ep1.5} follows from the classical energy estimate and the so-called Hoff estimates 
introduced in \cite{hoff1995global,hoff2002dynamics}. These estimates provide  $H^1$ regularity for the time derivative of the velocity $u$. The other estimate \eqref{eq3.34} can be seen as  a higher-order Hoff estimate 
and has be used in \cite{raphaelglobamuniquesolution} in the study of the incompressible Navier-Stokes model.
\end{rema}

The proof of \cref{Theorem 3.1} is the subject of \cref{prop..2} below and it is based on the classical  method of a priori estimates combined with  a homotopy argument.

\subsection{\texorpdfstring{Proof of \protect{\cref{Theorem 3.1}}}{}}\label{prop..2}
The proof of \cref{Theorem 3.1}  relies on a priori estimates and a homotopy argument and is structured into three sections: in  \cref{sect1}, we derive a priori estimates \eqref{ep1.5},  \eqref{eq3.36}, \eqref{eq3.34}, \eqref{eq3.48} and \eqref{eq3.49}, whereas \cref{proof1}-\cref{proof2} focus on the homotopy argument.
\subsubsection{A priori estimates}\label[section]{sect1}
The aim of  this section is to derive the a priori estimates \eqref{ep1.5}, \eqref{eq3.36}, \eqref{eq3.34}, \eqref{eq3.48} and \eqref{eq3.49}.
We observe that the regularity of $u$ is sufficient to use  use $u$ or $u_t$ as test functions in subsequent calculations. We use 
 $C$ to denote any constant that depends only on the dimension, $\widetilde\mu$, and $\alpha$. Conversely,  $C_*$ denotes a constant that may additionally depend polynomially on the lower and upper bounds of $\rho_0$, $\mu$ and $\lambda$.
\proof[\textbf{Proof of \eqref{ep1.5}}] Here we derive energy estimates for $u$ and $u_t$.
We take  the scalar product of the equation \eqref{ep1.6} with the velocity $u$. By doing so, we obtain:
\begin{gather*}
  \dfrac{1}{2}\dfrac{d}{dt}\norm{\sqrt{\rho_0} u}_{L^2(\R^d)}^{2}
+\int_{\mathbb{R}^{d}} 2\mu\abs{\D u}^2+\int_{\R^d}\lambda\abs{\dvg u}^{2} 
=-\int_{\R^d} F^{jk}\partial_k u^j.
\end{gather*}
We then integrate the above in time and we have:
\begin{gather}\label{ep10}
\dfrac{1}{2}\norm{\sqrt {\rho_{0}}u(t)}_{L^2(\R^d)}^2
+\int_0^t\int_{\mathbb{R}^{d}} 2\mu\abs{\D u}^2+\int_0^t\int_{\R^d}\lambda\abs{\dvg u}^{2}= \dfrac{1}{2}\norm{\sqrt{\rho_{0}} u_0}_{L^2(\R^d)}^2
-\int_0^t \int_{\R^d} F^{jk} \partial_ku^j.
\end{gather}
Using the identity:
\begin{gather}\label{c7.1}
\norm{\D w}_{L^2(\R^d)}^2=\dfrac{1}{2}\norm{\nabla w}^2_{L^2(\R^d)}+\dfrac{1}{2}\norm{\dvg w}_{L^2(\R^d)}^2
\end{gather}
we deduce that
\begin{gather}\label{c3.25}
\mu_*\left(\norm{\nabla w}_{L^2(\R^d)}^2+\norm{\dvg w}_{L^2(\R^d)}^2\right)\leqslant \int_{\R^d}2\mu(x) \abs{\D w(x)}^2dx.
\end{gather}
This inequality will be freely used in the forthcoming computations to translate the $L^2(\R^d)$-norm estimate of $\sqrt{\mu}\D w$ into that of $\nabla w$ with $w\in \{u,u_t,u_{tt}\}$.  For instance, by applying \eqref{c3.25} with $w=u$ and 
using H\"older's and Young's inequalities, we  obtain from \eqref{ep10}:
\begin{gather}\label{c3.41}
\sup_{[0,T]}\norm{\sqrt{\rho_0}u}_{L^2(\R^d)}^2+\int_0^T \norm{\nabla u}_{L^2(\R^d)}^2\leqslant C_* \norm{\sqrt{\rho_0}u_0}_{L^2(\R^d)}^2 + C_*\int_0^T\norm{F}_{L^2(\R^d)}^2.
\end{gather}

Next, we use $u_t$ as a test  function in the weak formulation of \eqref{ep2} in order to obtain:
\begin{align*}
\int_{0}^{t}\norm{\sqrt{\rho_{0}}  u_t}_{L^2(\R^d)}^{2}
&+\int_{\R^d}\mu(t)\abs{\D u(t)}^2+\dfrac{1}{2}\int_{\R^d}\lambda(t)(\dvg u(t))^2
=  \int_{\R^d}\mu_0\abs{\D u_0}^2+\dfrac{1}{2}\int_{\R^d}\lambda_0(\dvg u_0)^2\\
&-\int_{\R^d}F^{jk}(s) \partial_ku^j(s)\bigg|_{s=0}^{s=t}+\int_{0}^{t}\int_{\R^d} F^{jk}_t \partial_k u^j
+\int_{0}^{t}\int_{\mathbb{R}^{d}}\mu_t\abs{\D u}^2+\dfrac{1}{2}\int_{0}^{t}\int_{\mathbb{R}^{d}}\lambda_t (\dvg u)^{2}.
\end{align*}
In the following lines, we will estimate the right hand side in order to use Gr\"onwall's lemma. The third and the fourth terms are:
\begin{gather*}
    \Bigg|\int_{\R^d}F^{jk}(s) \partial_ku^j(s)\bigg|_{s=0}^{s=t}\Bigg| \leqslant 
    \eta \norm{\nabla u(t)}_{L^2(\R^d)}^2+ \dfrac{1}{4\eta} \norm{F(t)}_{L^{2}(\R^d)}^2+ \norm{\nabla u(0)}_{L^2(\R^d)}^2+ \norm{F(0)}_{L^{2}(\R^d)}^2,
\end{gather*}
for some positive $\eta>0$ small in order to absorb the norm of the gradient in the left hand side of the above equality. Next, the fifth term can be bounded as follows:
\[
\bigg|\int_{0}^{t}\int_{\R^d} F^{jk}_t \partial_k u^j\bigg|\leqslant \int_0^t\norm{\nabla u}_{L^2(\R^d)}^2 +\int_0^t\norm{ F_t}_{L^2(\R^d)}^2.
\]
Finally, the last two terms are controlled by:
\begin{gather*}
    \int_{0}^{t}\int_{\mathbb{R}^{d}}\mu_t\abs{\D u}^2+\dfrac{1}{2}\int_{0}^{t}\int_{\mathbb{R}^{d}}\lambda_t (\dvg u)^{2}
    \leqslant C_* \int_0^t \left(\norm{  \mu_t}_{L^\infty(\R^d)}+\norm{\lambda_t}_{L^\infty(\R^d)}\right)\norm{\nabla u}_{L^2(\R^d)}^2.
\end{gather*}
Gathering all of these estimates, and choosing $\eta$ small, one  has:
\begin{align}
\int_{0}^{t}\norm{\sqrt{\rho_{0}}  u_t}_{L^2(\R^d)}^{2}
+\sup_{[0,t]}\norm{\nabla u}^2_{L^2(\R^d)}&\leqslant C_*\norm{\nabla u_0}^2_{L^2(\R^d)}+C_*\sup_{[0,t]}\norm{F}_{L^{2}(\R^d)}^2
+C_*\int_0^t\norm{ F_t}_{L^2(\R^d)}^2\nonumber\\
&+C_*\int_0^t\norm{\mu_t,\, \lambda_t}_{L^\infty(\R^d)}\norm{\nabla u}_{L^2(\R^d)}^2.\label{c3.39}
\end{align}

We  differentiate \eqref{ep1.6} in time and we use $u_t$ as a test function in the weak formulation of the resulting equation. By doing so, we obtain:
\begin{align}
   \dfrac{1}{2}\int_{\R^d}\rho_{0}\abs{ u_t(t)}^{2}+\int_{0}^{t}
\int_{\R^d} 2\mu \abs{\D  u_t}_{L^2(\R^d)}^{2}&+\int_{0}^{t}\int_{\R^d}\lambda (\dvg u_t)^{2}
=\dfrac{1}{2}\int_{\R^d}\rho_{0} \abs{{u_t}_{|t=0}}^{2}\nonumber\\
&  -\int_0^t \int_{\R^d}2 \mu_t \D^{jk} u\partial_{k}u^j_t- \int_0^t \int_{\R^d} \lambda_t \dvg u \dvg  u_t+ \int_0^t\int_{\R^d}  F^{jk}_t\partial_ku^j_t\nonumber\\
&\leqslant \dfrac{1}{2}\int_{\R^d}\rho_{0} \abs{{u_t}_{|t=0}}^{2}+\eta \int_0^t\norm{\nabla  u_t}_{L^2(\R^d)}^2 + \dfrac{C_*}{\eta} \int_0^t \norm{ F_t}_{L^{2}(\R^d)}^2 \nonumber\\
&+ \dfrac{C_*}{\eta} \int_0^t\norm{\mu_t,\,\lambda_t}_{L^4(\R^d)}^2\norm{\nabla u}_{L^4(\R^d)}^2,\label{eq4.1}
\end{align}
with some positive $\eta $ small related to the lower bound of the viscosity $\mu$. We now turn to the estimate of the $L^4(\R^d)$ norm of $\nabla u$ that appears in the above inequality. We rewrite \eqref{ep1.6} as: 
\begin{gather}\label{c3.26}
-\widetilde\mu \Delta u-\nabla\{(\widetilde\mu+\lambda)\dvg u\}=\dvg \{2(\mu-\widetilde\mu)\D u\}-\rho_0 u_t +\dvg F.
\end{gather}
We denote by $\mathcal{P}=I-\mathcal{Q}$ and $\mathcal{Q}=-\nabla (-\Delta)^{-1}\dvg $  the usual projectors onto the spaces of divergence-free and curl-free vector fields, respectively.
First, we apply $\mathcal{P}$ to \eqref{c3.26} and we have:
\begin{gather}\label{c3.27}
\widetilde\mu\mathcal{P}u=(-\Delta)^{-1}\mathcal{P}\dvg \{2(\mu-\widetilde\mu)\D u+F\}-(-\Delta)^{-1}\mathcal{P}(\rho_0 u_t).
\end{gather}
Second, we apply the divergence operator to obtain:
\begin{gather}\label{c3.28}
(2\widetilde\mu+\lambda)\dvg u=(-\Delta)^{-1}\dvg \dvg \{2(\mu-\widetilde\mu)\D u+F\}-(-\Delta)^{-1}\dvg(\rho_0 u_t)
\end{gather}
and given that $\mathcal{Q} w=-\nabla (-\Delta)^{-1}\dvg w$, we find:
\[
\mathcal{Q} u= -\nabla (-\Delta)^{-1}\left[\dfrac{1}{2\widetilde\mu+\lambda}(-\Delta)^{-1}\dvg \dvg \{2(\mu-\widetilde\mu)\D u+F\}\right]+\nabla (-\Delta)^{-1}\left[\dfrac{1}{2\widetilde\mu+\lambda}(-\Delta)^{-1}\dvg(\rho_0 u_t)\right]
\]
Summing up, we express:
\begin{align}
    \nabla u&= \nabla \mathcal{P} u+\nabla\mathcal{Q} u\nonumber\\
            &=\frac{1}{\widetilde\mu}(-\Delta)^{-1}\nabla\left[  \mathcal{P}
\dvg\{2(\mu-\widetilde\mu)\D  u+F\}  -\mathcal{P}\left(  \rho_{0} u_t\right) \right] +\nabla^2 (-\Delta)^{-1}\left[\dfrac{1}{2\widetilde\mu+\lambda}(-\Delta)^{-1}\dvg(\rho_0 u_t)\right] \nonumber\\
&  -\nabla^2 (-\Delta)^{-1}\left[\dfrac{1}{2\widetilde\mu+\lambda}(-\Delta)^{-1}\dvg \dvg \{2(\mu-\widetilde\mu)\D u+F\}\right]. \label{ep1.9}
\end{align}
Furthermore, from \eqref{c3.26}, we have: 
\[
(-\Delta)^{-1}\nabla(\rho_0u_t)=\nabla(-\Delta)^{-1} \dvg (2\mu\D u)+\nabla(-\Delta)^{-1}\nabla(\lambda\dvg u)+\nabla (-\Delta)^{-1}\dvg F
\]
and owing to the fact that $\nabla(-\Delta)^{-1}\nabla$ is continuous over $L^p(\R^d)$, $p\in (1,\infty)$, we have:
\begin{gather}\label{c3.29}
\norm{(-\Delta)^{-1}\nabla(\rho_0u_t)}_{L^2(\R^d)}\leqslant C_*\norm{\nabla u, F}_{L^2(\R^d)}.
\end{gather}
Additionally, singular integral theory helps us to obtain from \eqref{ep1.9}:
\begin{align}
\norm{\nabla u}_{L^4(\R^d))}&\leqslant \dfrac{C}{\widetilde\mu}\norm{\mu-\widetilde\mu}_{L^\infty(\R^d)}\norm{\nabla u}_{L^4(\R^d))} +C\norm{F}_{L^4(\R^d))}+C\norm{(-\Delta)^{-1}\nabla(\rho_0 u_t)}_{ L^4(\R^d)}.\label{ep1.7}
\end{align}
Therefore, under the hypothesis \eqref{eq3.64}, we  absorb the first term 
of the right hand side in the left one.  With the help of  Gagliardo-Nirenberg's inequality and \eqref{c3.29}, the last terms above can be bounded as (recall $d\in \{2,3\})$:
\[
\norm{(-\Delta)^{-1}\nabla (\rho_0 u_t)}_{L^4(\R^d))}
\leqslant C_*\norm{\rho_0 u_t}_{L^2(\R^d)}^{d/4}\norm{\nabla u, F}_{L^2(\R^d)}^{1-d/4}.
\]
In consequence, the last term of \eqref{eq4.1} can be estimated as:
\begin{align}
\int_0^t\norm{\mu_t,\,\lambda_t}_{L^4(\R^d))}^2\norm{\nabla u}_{L^4(\R^d))}^2&\leqslant C_*\int_0^t\norm{\mu_t,\,\lambda_t}_{L^4(\R^d))}^2\norm{F}_{L^4(\R^d))}^2\nonumber\\
&+C_*\int_0^t\norm{\mu_t,\,\lambda_t}_{L^4(\R^d))}^2\norm{F,\; \nabla u,\;  u_t}_{L^2(\R^d))}^2.\label{eq4.2}
\end{align}
We now use \eqref{eq4.2} in \eqref{eq4.1},  choose $\eta$ related to the lower bound of the viscosity and density in order to absorb the second term of the right hand side of \eqref{eq4.1} in the left hand side. Applying Gronwall's Lemma and adding \eqref{c3.39} we find:
\begin{multline}\label{c3.40}
     \sup_{[0,T]}\norm{u_t,\, \nabla u}_{L^2(\R^d)}^{2}+\int_{0}^{T}
\norm{\nabla u_t,\, u_t}_{L^2(\R^d)}^{2}
\leqslant C_* \left(\norm{(F,u_0)}_{Y_T(\mu_0,\lambda_0)}^2+\int_0^T\norm{\mu_t,\,\lambda_t}_{L^4(\R^d))}^2\norm{F}_{L^2\cap L^6(\R^d))}^2\right)\\
\times \exp\left( C_*\int_0^T\left[\norm{\mu_t,\,\lambda_t}_{L^\infty(\R^d)}+\norm{\mu_t,\,\lambda_t}_{L^4(\R^d)}^2\right] \right).
\end{multline}
Finally, \eqref{ep1.5} follows from combining \eqref{c3.40} and  \eqref{c3.41}. We now turn to the proof of \eqref{eq3.36}.
\endproof
\dem[\textbf{Proof of \eqref{eq3.36}}] Here we derive $L^q(\R^d)$-norm estimate for the velocity gradient.
According to \eqref{ep1.9}, for $q\in \{6,8\}$ we have:
\begin{gather*}
\norm{\nabla u}_{L^q(\R^d))}\leqslant \dfrac{C}{\widetilde\mu}\norm{\mu-\widetilde\mu}_{L^\infty(\R^d)}\norm{\nabla u}_{L^q(\R^d)} +C\norm{F}_{ L^q(\R^d)}
+C\norm{(-\Delta)^{-1}\nabla (\rho_0 u_t)}_{L^q(\R^d)}.
\end{gather*}
 Therefore, under the hypothesis \eqref{eq3.64}, we absorb the first term 
of the right hand side in the left one. We only have to estimate the last term. Using Gagliardo-Nirenberg's inequality and \eqref{c3.29} we obtain: 
\[
\norm{(-\Delta)^{-1}\nabla(\rho_0 u_t)}_{L^6(\R^d))}
\leqslant C_*\norm{\rho_0 u_t}_{L^2(\R^d)}^{\beta}\norm{\nabla u, F}_{L^2(\R^d)}^{1-\beta}
\]
with $\beta= d/2-d/6\in \{2/3,1\}$. In consequence, we have:
\begin{gather}
\norm{\nabla u}_{L^\infty((0,T),L^6(\R^d))} \leqslant  C_*\norm{F}_{L^\infty((0,T),L^6(\R^d))}+ C_*\sup_{[0,T]}\norm{F, u_t, \nabla u}_{L^2(\R^d)}.
\end{gather}
For $q=8$ and $d=2$,  we have:
\[
\norm{(-\Delta)^{-1}\nabla(\rho_0 u_t)}_{L^8(\R^2))}
\leqslant C_* \norm{\rho_0 u_t}_{L^2(\R^2)}^{3/4}\norm{\nabla u, F}_{L^2(\R^2)}^{1/4}
\]
whereas for $d=3$, Sobolev's  and interpolation inequalities yield:
\[
\norm{(-\Delta)^{-1}\nabla \mathcal{P}(\rho_0 u_t)}_{L^8(\R^3))}
\leqslant C_*\norm{\rho_0 u_t}_{L^{24/11}(\R^3)}\leqslant C_* \norm{\nabla u_t}_{L^2(\R^3)}^{1/8}\norm{ u_t}_{L^2(\R^3)}^{7/8},
\]
such that,
\begin{gather}
\norm{\nabla u}_{L^{16}((0,T), L^8(\R^d))} \leqslant C_* \norm{ F}_{L^{16}((0,T), L^8(\R^d))}+ C_*\norm{F, u_t, \nabla u}_{L^2\cap L^\infty((0,T), L^2(\R^d))}+C_*\norm{\nabla u_t}_{L^2((0,T)\times\R^d)}.
\end{gather}
This concludes \eqref{eq3.36} and we  turn to the proof of \eqref{eq3.34}.
\enddem
\proof[\textbf{Proof of \eqref{eq3.34}}]
Here we establish estimates for $u_{tt}$. We first  differentiate equation \eqref{ep1.6} in time:
\begin{gather}\label{eq4.3}
\rho_0 u_{tt}-\dvg (2\mu\D u_t +\lambda \dvg u_t I_d)=\dvg F_t+\dvg (2\mu_t \D u+\lambda_t \dvg u I_d).
\end{gather}
Testing \eqref{eq4.3} with $\sigma u_{tt}$  (recall $\sigma(t)=\min \{1,t\}$ and hence $\sigma(0)=0$) leads to:
\begin{align}
\dfrac{1}{2}\int_0^t\sigma\norm{\sqrt{\rho_0} u_{tt}}_{L^2(\R^d)}^2&+\sigma(t)\int_{\R^d} \mu \abs{\D u_t}^2+
\dfrac{\sigma(t)}{2} \int_{\R^d} \lambda\abs{\dvg u_t}^2= \int_0^{t}\sigma'(t)\int_{\R^d} \mu \abs{\D u_t}^2+
\dfrac{1}{2} \int_0^{t}\sigma'(t)\int_{\R^d} \lambda\abs{\dvg u_t}^2\nonumber\\
&+\int_0^t\sigma\int_{\R^d} \mu_t \abs{\D u_t}^2+\dfrac{1}{2} \int_0^t\sigma\int_{\R^d} \lambda_t\abs{\dvg u_t}^2
-2\int_0^t\sigma\int_{\R^d}\mu_t \D^{jk} u \partial_k u^j_{tt}\nonumber\\
&-\int_0^t\sigma\int_{\R^d}\lambda_t \dvg u\dvg u_{tt}-\int_0^t\sigma\int_{\R^d} F_t^{jk}\partial_k u^j_{tt}.\label{eq3.39}
\end{align}
Using H\"older's inequality, the last term can be estimated as:
\[
\bigg|\int_0^t\sigma\int_{\R^d} F_t^{jk}\partial_k u^j_{tt}\bigg| \leqslant \eta \int_0^t\sigma^2 \norm{\nabla u_{tt}}_{L^2(\R^d)}^2
+\dfrac{1}{4\eta}\int_0^{t} \norm{F_t}_{L^2(\R^d)}^2
\]
and  the remaining terms are bounded by: 
\begin{align*}
\bigg|\dfrac{1}{2} \int_0^t\sigma\int_{\R^d} \lambda_t\abs{\dvg u_t}^2&-\int_0^t\sigma\int_{\R^d}\lambda_t \dvg u\dvg u_{tt}
+\int_0^t\sigma\int_{\R^d} \mu_t\abs{\D u_t}^2-2\int_0^t\sigma\int_{\R^d}\mu_t \D^{jk} u \partial_k u^j_{tt}\bigg|\\
&\leqslant \eta \int_0^t\sigma^2\norm{\nabla u_{tt}}_{L^2(\R^d)}^2+ \dfrac{C_*}{\eta}\int_0^t \norm{\mu_t,\,\lambda_t}_{L^4(\R^d)}^2\norm{\nabla u}_{L^4(\R^d)}^2\\
&+C_*\int_0^t\sigma \norm{\nabla u_t}_{L^2(\R^d)}^2\norm{\mu_t,\,\lambda_t}_{L^\infty(\R^d)}.
\end{align*}
Gathering all of these estimates, we have:
\begin{align}
    \int_0^t\sigma\norm{\sqrt{\rho_0} u_{tt}}_{L^2(\R^d)}^2&+\sigma(t)\int_{\R^d} \mu \abs{\D u_t}^2+
\dfrac{\sigma(t)}{2} \int_{\R^d} \lambda\abs{\dvg u_t}^2\nonumber\\
    &\leqslant  C_*\int_0^{t}\sigma'\norm{\nabla u_t}^2_{L^2(\R^d)}+2\eta\int_0^t \sigma^2\norm{\nabla u_{tt}}_{L^2(\R^d)}^2+\dfrac{C_*}{\eta} \int_0^{t} \norm{F_t}_{L^2(\R^d)}^2\nonumber\\
    &   +\dfrac{C_*}{\eta}   \int_0^t \norm{\mu_t,\,\lambda_t}_{L^4(\R^d))}^2\norm{\nabla u}_{L^4(\R^d))}^2+C_*\int_0^t \sigma\norm{\mu_{t},\lambda_{t}}_{L^\infty(\R^d)}\norm{\nabla u_t}_{L^2(\R^d)}^2.  \label{eq4.4}
\end{align}
The second term on the right hand side above will be reduced with the terms coming from the second energy which we turn to.
We  differentiate \eqref{eq4.3} with respect to time before multiplying the resulting equation with $\sigma^2u_{tt}$. By doing so, we obtain:
\begin{align*}
    \dfrac{\sigma^2}{2}\norm{\sqrt{\rho_0} u_{tt}}_{L^2(\R^d)}^2&+2\int_0^t\sigma^2\norm{\sqrt \mu\D u_{tt}}_{L^2(\R^d)}^2+
    \int_0^t\sigma^2\norm{\sqrt \lambda\dvg u_{tt}}_{L^2(\R^d)}^2= \int_0^{t}\sigma\sigma'\norm{\sqrt{\rho_0} u_{tt}}_{L^2(\R^d)}^2\\
    &-4\int_0^t \sigma^2 \mu_t \D^{jk} u_t \partial_k u_{tt}^j-2\int_0^t \sigma^2\int_{\R^d}\lambda_t \dvg u_t \dvg u_{tt}-\int_0^t\sigma^2\int_{\R^d} F^{jk}_{tt}
    \partial_k u^j_{tt}\\
    &-2\int_0^t \sigma^2\int_{\R^d} \mu_{tt} \D^{jk} u\partial_k u^j_{tt}-\int_0^t\sigma^2\int_{\R^d} \lambda_{tt} \dvg u \dvg u_{tt}.
\end{align*}
We estimate the second and third terms of the right hand side as follows:
\[
\bigg|4\int_0^t \sigma^2 \mu_t \D^{jk} u_t \partial_k u_{tt}^j+2\int_0^t \sigma^2\int_{\R^d}\lambda_t \dvg u_t \dvg u_{tt}\bigg|
\leqslant \eta \int_0^t \sigma^2 \norm{\nabla u_{tt}}^2_{L^2(\R^d)}+ \dfrac{C_*}{\eta}\int_0^t \sigma^2 \norm{\mu_t,\,\lambda_t}_{L^\infty(\R^d)}^2\norm{\nabla u_t}_{L^2(\R^d)}^2,
\]
whereas the last two terms are bounded by 
\begin{align*}
\bigg|2\int_0^t \sigma^2\int_{\R^d} \mu_{tt} \D^{jk} u\partial_k u^j_{tt}&+\int_0^t\sigma^2\int_{\R^d} \lambda_{tt} \dvg u \dvg u_{tt} \bigg|\\
&\leqslant  \eta \int_0^t \sigma^2 \norm{\nabla u_{tt}}^2_{L^2(\R^d)}+ \dfrac{C_*}{\eta}\int_0^t \sigma^2 \norm{\mu_{tt},\,\lambda_{tt}}_{L^4(\R^d)}^2\norm{\nabla u}_{L^4(\R^d)}^2.
\end{align*}
Finally, the fourth term of the right hand side is:
\[
\bigg|\int_0^t\sigma^2\int_{\R^d} F^{jk}_{tt}
    \partial_k u^j_{tt}\bigg|\leqslant  \eta \int_0^t \sigma^2 \norm{\nabla u_{tt}}^2_{L^2(\R^d)}+ \dfrac{1}{4\eta} \int_0^t\sigma^2\norm{F_{tt}}_{L^2(\R^d)}^2.
\]
Gathering all of these estimates, adding \eqref{eq4.4} and choosing $\eta$ small related to the lower bounds of the viscosity, we finally apply Gr\"onwall's Lemma, leading to  (noticing $\sigma'(t)=0$ for $t>1$): 
\begin{multline}
    \int_0^T\norm{ \sigma u_{tt},\; \sigma^2 \nabla u_{tt}}_{L^2(\R^d)}^2+\norm{\sigma \nabla u_t,\; \sigma^2 u_{tt}}^2_{L^\infty((0,T),L^2(\R^d))}\leqslant C_* \left[ \int_0^{\sigma(T)}\norm{\nabla u_t}^2_{L^2(\R^d)}+\int_0^{T} \norm{F_t,\, \sigma^2F_{tt}}_{L^2(\R^d)}^2\right.\nonumber\\
    \left.+  \int_0^T\norm{\mu_t,\,\lambda_t,\, \sigma\mu_{tt},\,\sigma\lambda_{tt}}_{L^4(\R^d)}^2\norm{\nabla u}_{L^4(\R^d)}^2\right]\exp\left(C_*\int_0^T\left[\norm{\mu_t,\,\lambda_t}_{L^\infty(\R^d)}+\sigma^2 \norm{\mu_t,\,\lambda_t}_{L^\infty(\R^d)}^2\right]\right).
\end{multline}
This concludes \eqref{eq3.34} and we turn to \eqref{eq3.48}.
\endproof
\dem[\textbf{Proof of \eqref{eq3.48}}] Here we derive $L^{8/d}((0,T), L^4( \R^d))$-norm estimate for $\sigma^{\tfrac{d}{4}}\nabla u_t$.
Differentiating \eqref{c3.27}-\eqref{c3.28} with respect to time yields:
\begin{gather*}
\begin{cases}\vspace{0,2cm}
\widetilde\mu\mathcal{P}u_t&=(-\Delta)^{-1}\mathcal{P}\dvg \{2(\mu-\widetilde\mu)\D u_t+2\mu_t\D u+F_t\}-(-\Delta)^{-1}\mathcal{P}(\rho_0 u_{tt}),\\
(2\widetilde\mu+\lambda)\dvg u_t&=-\lambda_t\dvg u+(-\Delta)^{-1}\dvg \dvg \{2(\mu-\widetilde\mu)\D u_t+2\mu_t\D u+F_t\}-(-\Delta)^{-1}\dvg(\rho_0 u_{tt})
\end{cases}
\end{gather*}
and therefore (recall $\mathcal{Q} w=-\nabla (-\Delta)^{-1}\dvg w$):
\begin{align*}
    \nabla u_t&=\nabla \mathcal{P} u_t+\nabla \mathcal{Q} u_t\\
    &=\frac{1}{\widetilde\mu}(-\Delta)^{-1}\nabla\left[  \mathcal{P}
\dvg\{2(\mu-\widetilde\mu)\D  u_t+2\mu_t \D u+F_t\}  -\mathcal{P}\left(  \rho_{0}u_{tt}\right) \right] +\nabla^2 (-\Delta)^{-1}\left[\dfrac{1}{2\widetilde\mu+\lambda}(-\Delta)^{-1}\dvg(\rho_0 u_{tt})\right] \nonumber\\
& +\nabla^2 (-\Delta)^{-1}\left[\dfrac{\lambda_t}{2\widetilde\mu+\lambda}\dvg u\right] -\nabla^2 (-\Delta)^{-1}\left[\dfrac{1}{2\widetilde\mu+\lambda}(-\Delta)^{-1}\dvg \dvg \{2(\mu-\widetilde\mu)\D u_t+2\mu_t\D u+F_t\}\right].
\end{align*}
Whence, under the hypothesis \eqref{eq3.64} we obtain:
\[
\norm{\nabla u_t}_{L^4(\R^d)}\leqslant C_*\norm{F_t}_{L^4(\R^d)}+C_*\norm{\nabla u}_{L^8(\R^d)}\norm{\mu_t,\,\lambda_t}_{L^8(\R^d)}
+C_*\norm{(-\Delta)^{-1}\nabla(\rho_{0}u_{tt})}_{L^4(\R^d)}
\]
and we only need to estimate the last term.  Following the computations leading to \eqref{c3.29}, we have from \eqref{eq4.3}:
\begin{gather}\label{c3.30}
\norm{(-\Delta)^{-1}\nabla(\rho_0 u_{tt})}_{L^2(\R^d)}\leqslant C_*\norm{\nabla u_t,\, F_t,\,\lambda_t \dvg u,\, \mu_t\D u}_{L^2(\R^d)}.
\end{gather}
and  Gagliardo-Nirenberg's inequality implies:
\begin{align*}
\int_0^T\sigma^2 \norm{(-\Delta)^{-1}\nabla  (\rho_{0}u_{tt})}_{L^4(\R^d)}^{8/d} 
&\leqslant C_*\int_0^T\sigma^{2}\norm{\nabla u_t,\, F_t,\,\lambda_t \dvg u,\, \mu_t\D u}_{L^2(\R^d)}^{\tfrac{8}{d}-2}\norm{u_{tt}}_{L^2(\R^d)}^{2}\\
&\leqslant C_*\norm{\sigma u_{tt}}_{L^\infty((0,T),L^2(\R^2))}^{2}\int_0^T\norm{\nabla u_t,\, F_t,\,\lambda_t \dvg u,\, \mu_t\D u}_{L^2(\R^d)}^{\tfrac{8}{d}-2}.
\end{align*}
In summary, we obtain:
\begin{align*}
       \int_0^T\sigma^{2}\norm{\nabla u_t}_{L^4(\R^d)}^{8/d}&\leqslant C_*\int_0^T\sigma^{2}\norm{F_t}_{L^4(\R^d)}^{8/d}+C_*\int_0^T\sigma^{2}\norm{\nabla u}_{L^8(\R^d)}^{8/d}\norm{\mu_t,\,\lambda_t}_{L^8(\R^d)}^{8/d}\\
       &+C_*\norm{\sigma u_{tt}}_{L^\infty((0,T),L^2(\R^d))}^2\int_0^T\norm{\nabla u_t,\, F_t,\,\lambda_t \dvg u,\, \mu_t\D u}_{L^2(\R^d)}^{\tfrac{8}{d}-2}.
\end{align*}
This ends the proof of \eqref{eq3.48} and we turn to the last a priori estimate.
\enddem


\dem[\textbf{Proof of \eqref{eq3.49}}] 
\label{pagecCalpha}
Here, we  translate the bounds of $u_t$ into a piecewise H\"older regularity for the velocity gradient. To achieve this, we will use the results from \cref{app1} to derive a piecewise H\"older estimate for the discontinuous terms in the expression of the velocity gradient \eqref{ep1.9}. We fix $q=2d/\alpha$.
\paragraph*{\textbf{Step 1: Terms involving $F$}} We estimate the terms 
 \[
    \frac{1}{\widetilde\mu}(-\Delta)^{-1}\nabla\mathcal{P}
\dvg F \quad\text{and} \quad\nabla^2 (-\Delta)^{-1}\left[\dfrac{1}{2\widetilde\mu+\lambda}(-\Delta)^{-1}\dvg \dvg F\right].
 \]
    With the help of \eqref{c3.31} and \eqref{ez1.96}, we have:
    \begin{gather}\label{c3.43}
    \begin{cases}\vspace{0,2cm}
    \norm{(-\Delta)^{-1}\nabla^2 F,\; (-\Delta)^{-1}\nabla\mathcal{P}
\dvg F}_{L^\infty(\R^d)}&\leqslant C\left( \ell_{\vph_0}^{-\tfrac{\alpha}{2}}\norm{F}_{L^q(\R^d)}+\norm{F}_{\cC^\alpha_{pw,\mathcal{C}_0}(\R^d)}\right),\\
    \norm{(-\Delta)^{-1}\nabla^2 F,\; (-\Delta)^{-1}\nabla\mathcal{P}
\dvg F}_{\dot \cC^\alpha_{pw,\mathcal{C}_0}(\R^d)}&\leqslant C(1+\mathfrak{P}_{\mathcal{C}_0})\norm{F}_{\cC^\alpha_{pw,\mathcal{C}_0}(\R^d)},
    \end{cases}
    \end{gather}
    where $\mathfrak{P}_{\mathcal{C}_0}$ is defined in \eqref{c3.34}. Similarly, we have:
    \begin{align*}
        \left\|\nabla^2 (-\Delta)^{-1}\right.&\left.\left[\dfrac{1}{2\widetilde\mu+\lambda}(-\Delta)^{-1}\dvg \dvg F\right]\right\|_{L^\infty(\R^d)}\\
        &\leqslant C\left[\left(\norm{\lambda}_{\dot\cC^\alpha_{pw,\mathcal{C}_0}(\R^d)}+\ell_{\vph_0}^{-\tfrac{\alpha}{2}}\right)\norm{F}_{L^q(\R^d)}+\norm{(-\Delta)^{-1}\dvg \dvg  F}_{\cC^\alpha_{pw,\mathcal{C}_0}(\R^d)}\right]\\
        &\leqslant C\left[\left(\norm{\lambda}_{\dot\cC^\alpha_{pw,\mathcal{C}_0}(\R^d)}+\ell_{\vph_0}^{-\tfrac{\alpha}{2}}\right)\norm{F}_{L^q(\R^d)}+(1+\mathfrak{P}_{\mathcal{C}_0})\norm{F}_{\cC^\alpha_{pw,\mathcal{C}_0}(\R^d)}\right],
    \end{align*}
    and 
    \begin{align*}
        \left\|\nabla^2 (-\Delta)^{-1}\left[\dfrac{1}{2\widetilde\mu+\lambda}(-\Delta)^{-1}\dvg \dvg F\right]\right\|_{\dot \cC^\alpha_{pw,\mathcal{C}_0}(\R^d)}&\leqslant C \norm{(-\Delta)^{-1}\dvg \dvg  F}_{\dot \cC^\alpha_{pw,\mathcal{C}_0}(\R^d)}\\
        &+C\left(\norm{\lambda}_{\dot\cC^\alpha_{pw,\mathcal{C}_0}(\R^d)}+\mathfrak{P}_{\mathcal{C}_0}\right)\norm{(-\Delta)^{-1}\dvg \dvg  F}_{L^\infty(\R^d)}\\
        &\leqslant C\left(1+\mathfrak{P}_{\mathcal{C}_0}+\norm{\lambda}_{\dot\cC^\alpha_{pw,\mathcal{C}_0}(\R^d)}\right)\norm{F}_{\cC^\alpha_{pw,\mathcal{C}_0}(\R^d)} \\
        &+C\left(\norm{\lambda}_{\dot\cC^\alpha_{pw,\mathcal{C}_0}(\R^d)}+\mathfrak{P}_{\mathcal{C}_0}\right)\ell_{\vph_0}^{-\tfrac{\alpha}{2}}\norm{F}_{L^q(\R^d)}.
    \end{align*}
    In sum, we have: 
    \begin{align}
        \left\|(-\Delta)^{-1}\nabla\mathcal{P}
\dvg F,\right.&\left.\;\nabla^2 (-\Delta)^{-1}\left[\dfrac{1}{2\widetilde\mu+\lambda}(-\Delta)^{-1}\dvg \dvg F\right]\right\|_{\cC^\alpha_{pw,\mathcal{C}_0}(\R^d)} \nonumber\\
&\leqslant C\left(1+\mathfrak{P}_{\mathcal{C}_0}+\norm{\lambda}_{\dot\cC^\alpha_{pw,\mathcal{C}_0}(\R^d)}\right)\left(\ell_{\vph_0}^{-\tfrac{\alpha}{2}}\norm{F}_{L^q(\R^d)}+\norm{F}_{\cC^\alpha_{pw,\mathcal{C}_0}(\R^d)}\right).\label{c3.53}
    \end{align}
 \paragraph*{\textbf{Step 2: Small terms}}
We  estimate the terms
    \[
    \frac{1}{\widetilde\mu}(-\Delta)^{-1}\nabla\mathcal{P}
\dvg\{2(\mu-\widetilde\mu)\D  u\} \quad\text{and}\quad \nabla^2 (-\Delta)^{-1}\left[\dfrac{1}{2\widetilde\mu+\lambda}(-\Delta)^{-1}\dvg \dvg \{2(\mu-\widetilde\mu)\D u\}\right].
    \]
     Again, with the help of \eqref{c3.31}, we have:
    \begin{align}
    \norm{(-\Delta)^{-1}\nabla^2\{2(\mu-\widetilde\mu)\D  u\}}_{L^\infty(\R^d)}&\leqslant C\left(\norm{\mu}_{\dot \cC^\alpha_{pw,\mathcal{C}_0}(\R^d)}+\ell_{\vph_0}^{-\tfrac{\alpha}{2}}\norm{\mu-\widetilde\mu}_{L^\infty(\R^d)}\right)\norm{\nabla u}_{L^q(\R^d)}\nonumber\\
    &+C \norm{\mu-\widetilde\mu}_{L^\infty(\R^d)}\norm{\nabla u}_{\cC^{\alpha}_{pw,\mathcal{C}_0}(\R^d)}.\label{c3.44}
    \end{align} 
     Next, we gather \cref{exten} and \eqref{ez1.96} to obtain:
   \begin{align*}
       \norm{(-\Delta)^{-1}\nabla^2(w)}_{\dot \cC^\alpha_{pw,\mathcal{C}_0}(\R^d)}&\leqslant \norm{(-\Delta)^{-1}\nabla^2(w^e)}_{\dot \cC^\alpha_{pw,\mathcal{C}_0}(\R^d)}+\norm{(-\Delta)^{-1}\nabla^2(w-w^e)}_{\dot \cC^\alpha_{pw,\mathcal{C}_0}(\R^d)}\\
       &\leqslant C\left[\norm{w^e}_{\dot \cC^\alpha(\R^d)}+ \norm{w-w^e}_{\dot \cC^\alpha_{pw,\mathcal{C}_0}(\R^d)}+\mathfrak{P}_{\mathcal{C}_0} \norm{w-w^e}_{L^\infty(\R^d)}
\right]\\
       &\leqslant C\left[\norm{w}_{\dot \cC^\alpha_{pw,\mathcal{C}_0}(\R^d)}+ \left(\ell^{-\alpha}_{\vph_0}+\mathfrak{P}_{\mathcal{C}_0}\right)\norm{\llbracket w\rrbracket}_{L^\infty(\mathcal{C}_0)}\right].
   \end{align*}
   Therefore:
    \begin{gather}\label{c3.45}
        \norm{(-\Delta)^{-1}\nabla^2\{2(\mu-\widetilde\mu)\D  u\}}_{\dot \cC^\alpha_{pw,\mathcal{C}_0}(\R^d)}
        \leqslant C\left[\norm{(\mu-\widetilde\mu)\D  u}_{\dot\cC^\alpha_{pw,\mathcal{C}_0}(\R^d)}+\left( \ell^{-\alpha}_{\vph_0}+\mathfrak{P}_{\mathcal{C}_0}\right)\norm{\llbracket (\mu-\widetilde\mu)\D  u\rrbracket}_{L^\infty(\mathcal{C}_0)}\right].
    \end{gather}
    Similarly, we have:
    \begin{align} 
        \left\|\nabla^2 (-\Delta)^{-1}\vphantom{\dfrac{1}{2\widetilde\mu+\lambda}}\right.&\left.\left[\dfrac{1}{2\widetilde\mu+\lambda}(-\Delta)^{-1}\dvg \dvg \{2(\mu-\widetilde\mu)\D u\}\right]\right\|_{L^\infty(\R^d)}\nonumber\\ \vspace{0,3cm}
        &\leqslant C\left(\norm{\lambda}_{\dot\cC^\alpha_{pw,\mathcal{C}_0}(\R^d)}+\ell_{\vph_0}^{-\tfrac{\alpha}{2}}\right)\norm{\mu-\widetilde\mu}_{L^\infty(\R^d)}\norm{\nabla u}_{L^q(\R^d)}+ C \norm{(-\Delta)^{-1}\dvg \dvg \{2(\mu-\widetilde\mu)\D u\}}_{\cC^{\alpha}_{pw,\mathcal{C}_0}(\R^d)}\nonumber\\ \vspace{0,3cm}
        &\leqslant C\left[\norm{\mu}_{\dot\cC^\alpha_{pw,\mathcal{C}_0}(\R^d)}+\left(\norm{\lambda}_{\dot\cC^\alpha_{pw,\mathcal{C}_0}(\R^d)}+\ell_{\vph_0}^{-\tfrac{\alpha}{2}}\right)\norm{\mu-\widetilde\mu}_{L^\infty(\R^d)}\right]\norm{\nabla u}_{L^q(\R^d)}\nonumber\\
        &+C\left( \mathfrak{P}_{\mathcal{C}_0} +\ell^{-\alpha}_{\vph_0}\right)\norm{\llbracket (\mu-\widetilde\mu)\D  u\rrbracket}_{L^\infty(\mathcal{C}_0)}+C\norm{\mu-\widetilde\mu}_{ \cC^\alpha_{pw,\mathcal{C}_0}(\R^d)}\norm{\nabla u}_{ \cC^{\alpha}_{pw,\mathcal{C}_0}(\R^d)},\label{c3.46}
    \end{align}
    and 
    \begin{align}
         \left\|\nabla^2 (-\Delta)^{-1}\vphantom{\dfrac{1}{2\widetilde\mu+\lambda}}\right.&\left.\left[\dfrac{1}{2\widetilde\mu+\lambda}(-\Delta)^{-1}\dvg \dvg \{2(\mu-\widetilde\mu)\D u\}\right]\right\|_{\dot \cC^\alpha_{pw,\mathcal{C}_0}(\R^d)}\nonumber\\ 
         &\leqslant C\left\|\dfrac{1}{2\widetilde\mu+\lambda}(-\Delta)^{-1}\dvg \dvg \{2(\mu-\widetilde\mu)\D u\}\right\|_{\dot\cC^\alpha_{pw,\mathcal{C}_0}(\R^d)}\nonumber\\ \vspace{0,3cm}
         &+C\left( \mathfrak{P}_{\mathcal{C}_0} +\ell^{-\alpha}_{\vph_0}\right)\left\|\llbracket \dfrac{1}{2\widetilde\mu+\lambda}(-\Delta)^{-1}\dvg \dvg \{2(\mu-\widetilde\mu)\D u\}\rrbracket\right\|_{L^\infty(\mathcal{C}_0)}\nonumber\\ \vspace{0,3cm}
         &\leqslant C\left(\norm{\lambda}_{\dot \cC^\alpha_{pw,\mathcal{C}_0}(\R^d)}+\left( \mathfrak{P}_{\mathcal{C}_0} +\ell^{-\alpha}_{\vph_0}\right) \norm{\llbracket\lambda\rrbracket}_{L^\infty(\mathcal{C}_0)}\right)\norm{(-\Delta)^{-1}\dvg \dvg \{2(\mu-\widetilde\mu)\D u\}}_{L^\infty(\R^d)} \nonumber\\ \vspace{0,3cm}
         &+C\left( \mathfrak{P}_{\mathcal{C}_0} +\ell^{-\alpha}_{\vph_0}\right)\left\|\llbracket (-\Delta)^{-1}\dvg \dvg \{2(\mu-\widetilde\mu)\D u\}\rrbracket\right\|_{L^\infty(\mathcal{C}_0)} \nonumber\\ \vspace{0,3cm}
         &+C\left\|(-\Delta)^{-1}\dvg \dvg \{2(\mu-\widetilde\mu)\D u\}\right\|_{\dot\cC^\alpha_{pw,\mathcal{C}_0}(\R^d)} \nonumber\\ \vspace{0,3cm}
         &\leqslant C\left(\norm{\lambda}_{\dot \cC^\alpha_{pw,\mathcal{C}_0}(\R^d)}+\left( \mathfrak{P}_{\mathcal{C}_0} +\ell^{-\alpha}_{\vph_0}\right)\norm{ \llbracket\lambda\rrbracket}_{L^\infty(\mathcal{C}_0)}\right)\left(\norm{\mu}_{\dot \cC^\alpha_{pw,\mathcal{C}_0}(\R^d)}+\ell_{\vph_0}^{-\tfrac{\alpha}{2}}\norm{\mu-\widetilde\mu}_{L^\infty(\R^d)}\right)\norm{\nabla u}_{L^q(\R^d)} \nonumber\\ \vspace{0,25cm}
         &+C\left[1+\norm{\lambda}_{\dot \cC^\alpha_{pw,\mathcal{C}_0}(\R^d)}+\left( \mathfrak{P}_{\mathcal{C}_0} +\ell^{-\alpha}_{\vph_0}\right)\norm{\llbracket\lambda\rrbracket}_{L^\infty(\mathcal{C}_0)}\right]\norm{\mu-\widetilde\mu}_{\cC^\alpha_{pw,\mathcal{C}_0}(\R^d)}\norm{\nabla u}_{\cC^{\alpha}_{pw,\mathcal{C}_0}(\R^d)}\nonumber\\
         &+C\left( \mathfrak{P}_{\mathcal{C}_0} +\ell^{-\alpha}_{\vph_0}\right)\left\|\llbracket (-\Delta)^{-1}\dvg \dvg \{2(\mu-\widetilde\mu)\D u\}\rrbracket,\; \llbracket (\mu-\widetilde\mu)\D  u\rrbracket\right\|_{L^\infty(\mathcal{C}_0)}.\label{c3.47}
    \end{align}  
We now turn to estimating the norms of the jumps involved in the above inequalities.  We first notice that the solution $u$ of \eqref{ep2} satisfies:
   \begin{gather}\label{c3.32}
   \llbracket 2\mu \D^{jk} u+ \lambda \dvg u \delta^{jk} + F^{jk}\rrbracket n_x^k=0,
   \end{gather}
   where $n_x$ is the normal vector field of $\mathcal{C}_0$ and $\delta^{jk}$ denotes the Kronecker symbol. We consider an unit vector field $\tau_{x,l}$, $l\in \llbracket 1,d-1\rrbracket$ such that $(\tau_{x,1},\dots,\, \tau_{x,d-1},\,n_x)$ forms a basis of $\R^d$. Given that $u$ is continuous in the whole space and $\nabla u$ is continuous on both sides of $\mathcal{C}_0$,
    then discontinuity in $\nabla u$ may  appear only in the normal direction of $\mathcal{C}_0$. Namely, we have:
    \[
    \forall\, l\in \llbracket 1, d-1\rrbracket,\quad \llbracket \nabla u\rrbracket \cdot \tau_{x,l}=0
    \]
    and there exists a vector field $\Vec{\mathrm a}\in \R^d$ such that:
    \begin{gather}\label{c3.33}
    \llbracket \nabla u\rrbracket= \Vec{\mathrm{a}}\otimes n_x.
    \end{gather}
If $\scalar{g}$ denotes the average of $g$ at the interface, then we rewrite \eqref{c3.32} as
\[
2\scalar{\mu}\llbracket \D^{jk} u\rrbracket  n_x^{k}+ \scalar{\lambda}\llbracket \dvg u\rrbracket n_x^j =-2\llbracket \mu\rrbracket\scalar{\D^{jk} u} n_x^k-\llbracket\lambda\rrbracket \scalar{\dvg u}n_x^j-\llbracket F^{jk}\rrbracket n_x^k,
\]
and scalar products with $n_x$ and $\tau_{x,l}$ together with \eqref{c3.33} lead to:
\begin{gather}\label{c3.42}
\begin{cases}
    \scalar{2\mu+\lambda}\Vec{\mathrm{a}}\cdot n_x= -2\llbracket \mu\rrbracket\scalar{\D^{jk} u} n_x^kn_x^j-\llbracket\lambda\rrbracket \scalar{\dvg u}-\llbracket F^{jk}\rrbracket n_x^kn_x^j,\\
    \scalar{\mu}\Vec{\mathrm{a}}\cdot \tau_{x,l}= -2\llbracket \mu\rrbracket\scalar{\D^{jk} u} n_x^k\tau_{x,l}^j-\llbracket F^{jk}\rrbracket n_x^k\tau_{x,l}^j.
\end{cases}
\end{gather}
Therefore:
\begin{align*}
\llbracket (\mu-\widetilde\mu)\D u\rrbracket&=\llbracket \mu \rrbracket \scalar{\D u}+ \scalar{\mu-\widetilde\mu}\llbracket \D u\rrbracket\\
        &=\llbracket \mu \rrbracket \scalar{\D u}+ \dfrac{\scalar{\mu-\widetilde\mu}}{\scalar{\mu}}\scalar{\mu}\llbracket \D u\rrbracket\\
\end{align*}
and recalling \eqref{c3.33}, \eqref{c3.42} with
\[
\Vec{\mathrm{a}}=\Vec{\mathrm{a}}\cdot n_x n_x+ \sum_l\Vec{\mathrm{a}} \cdot \tau_{x,l} \tau_{x,l},
\]
we infer:
\begin{align}
\norm{\llbracket (\mu-\widetilde\mu)\D  u\rrbracket}_{L^\infty(\mathcal{C}_0)} &\leqslant 
C_*\norm{F}_{L^\infty(\R^d)}+C\norm{\llbracket\lambda\rrbracket}_{L^\infty(\mathcal{C}_0)}\left\|\dfrac{\scalar{\mu-\widetilde\mu}}{\scalar{\mu}}\right\|_{L^\infty(\mathcal{C}_0)}\norm{\nabla u}_{L^\infty(\R^d)}\nonumber\\
&+ C\norm{\llbracket\mu\rrbracket}_{L^\infty(\mathcal{C}_0)}\left(1+\left\|\dfrac{\scalar{\mu-\widetilde\mu}}{\scalar{\mu}}\right\|_{L^\infty(\mathcal{C}_0)}\right)\norm{\nabla u}_{L^\infty(\R^d)}.\label{c3.48}
\end{align}
Moreover, from \eqref{c3.28} and $\eqref{c3.42}_1$, we find (recall $\llbracket\dvg u\rrbracket=a\cdot n_x$):
\begin{align*}
\llbracket (-\Delta)^{-1}\dvg \dvg \{(\mu-\widetilde\mu)\D u\}\rrbracket&=\llbracket (2\widetilde\mu+\lambda)\dvg u\rrbracket - \llbracket (-\Delta)^{-1}\dvg \dvg F\rrbracket\\
&=\dfrac{\scalar{2\widetilde\mu+\lambda}}{\scalar{2\mu+\lambda}}\scalar{2\mu+\lambda} \llbracket \dvg u\rrbracket +\llbracket \lambda\rrbracket \scalar{\dvg u}- \llbracket (-\Delta)^{-1}\dvg \dvg F\rrbracket\\
&=\dfrac{\scalar{2\widetilde\mu+\lambda}}{\scalar{2\mu+\lambda}}\left(-2\llbracket \mu\rrbracket\scalar{\D^{jk} u} n_x^kn_x^j-\llbracket F^{jk}\rrbracket n_x^kn_x^j\right) \\
&+2\dfrac{\scalar{\mu-\widetilde\mu}}{\scalar{2\mu+\lambda}}\llbracket \lambda\rrbracket \scalar{\dvg u}- \llbracket (-\Delta)^{-1}\dvg \dvg F\rrbracket
\end{align*}
and hence (recall \eqref{c3.43}):
\begin{align}
\norm{\llbracket (-\Delta)^{-1}\dvg \dvg &\{(\mu-\widetilde\mu)\D u\}\rrbracket}_{L^\infty(\mathcal{C}_0)}\leqslant  C\norm{\llbracket\mu\rrbracket}_{L^\infty(\mathcal{C}_0)}\left(1+\left\|\dfrac{\scalar{\mu-\widetilde\mu}}{\scalar{\mu}}\right\|_{L^\infty(\mathcal{C}_0)}\right)\norm{\nabla u}_{L^\infty(\R^d)}\nonumber\\
&+C\norm{\llbracket\lambda\rrbracket}_{L^\infty(\mathcal{C}_0)}\left\|\dfrac{\scalar{\mu-\widetilde\mu}}{\scalar{\mu}}\right\|_{L^\infty(\mathcal{C}_0)}\norm{\nabla u}_{L^\infty(\R^d)}+C\left( \ell_{\vph_0}^{-\tfrac{\alpha}{2}}\norm{F}_{L^q(\R^d)}+\norm{F}_{\cC^\alpha_{pw,\mathcal{C}_0}(\R^d)}\right). \label{c3.49}
\end{align}
Summing up \eqref{c3.44}, \eqref{c3.45}, \eqref{c3.46}, \eqref{c3.47}, \eqref{c3.48} and \eqref{c3.49} we find: 
\begin{align}
    \left\|(-\Delta)^{-1}\right.&\left.\nabla\mathcal{P}\vphantom{\dfrac{1}{2\widetilde\mu+\lambda}}
\dvg\{(\mu-\widetilde\mu)\D  u\},\;\nabla^2 (-\Delta)^{-1}\left[\dfrac{1}{2\widetilde\mu+\lambda}(-\Delta)^{-1}\dvg \dvg \{(\mu-\widetilde\mu)\D u\}\right]\right\|_{\cC^\alpha_{pw,\mathcal{C}_0}(\R^d)}\nonumber\\
&\leqslant C_{**}\left( \norm{\nabla u}_{L^q(\R^d)}+\norm{F}_{L^q(\R^d)}+\norm{F}_{\cC^\alpha_{pw,\mathcal{C}_0}(\R^d)}\right)+C\left( \mathfrak{P}_{\mathcal{C}_0} +\ell^{-\alpha}_{\vph_0}\right) \norm{\llbracket\mu\rrbracket}_{L^\infty(\mathcal{C}_0)}\norm{\nabla u}_{L^\infty(\R^d)}\nonumber\\
&+C\left(1+\norm{\lambda}_{\dot \cC^\alpha_{pw,\mathcal{C}_0}(\R^d)}+\left( \mathfrak{P}_{\mathcal{C}_0} +\ell^{-\alpha}_{\vph_0}\right)\norm{\llbracket\lambda\rrbracket}_{L^\infty(\mathcal{C}_0)}\right)\norm{\mu-\widetilde\mu}_{\cC^\alpha_{pw,\mathcal{C}_0}(\R^d)}\norm{\nabla u}_{\cC^{\alpha}_{pw,\mathcal{C}_0}(\R^d)}\nonumber\\
&+C\left( \mathfrak{P}_{\mathcal{C}_0} +\ell^{-\alpha}_{\vph_0}\right) \norm{\llbracket\lambda,\; \mu\rrbracket}_{L^\infty(\mathcal{C}_0)}\left\|\dfrac{\scalar{\mu-\widetilde\mu}}{\scalar{\mu}}\right\|_{L^\infty(\mathcal{C}_0)}\norm{\nabla u}_{L^\infty(\R^d)}\label{c3.52}
\end{align}
Above, the constant $C_{**}$ depends polynomially on the lower and upper bounds of the viscosity and density,  on $\displaystyle\sup_{[0,T]}\norm{\mu,\,\lambda}_{\dot\cC^\alpha_{pw,\mathcal{C}_0}(\R^d)}$, and  on $\ell^{-\tfrac{\alpha}{2}}_{\vph_0}$ and $\mathfrak{P}_{\mathcal{C}_0}$.
\paragraph*{\textbf{Step 3: Remaining terms}}
We estimate the terms
    \[
    \frac{1}{\widetilde\mu}(-\Delta)^{-1}\nabla\mathcal{P}\left(  \rho_{0}u_t\right) \quad\text{and}\quad\nabla^2 (-\Delta)^{-1}\left[\dfrac{1}{2\widetilde\mu+\lambda}(-\Delta)^{-1}\dvg(\rho_0 u_t)\right]. 
    \]
 The first term above is H\"older continuous in the whole space, given the regularity previously obtained  for $u_t$. Indeed, 
  by the interpolation inequality and \eqref{c3.29}:
\begin{align}
\int_0^T \sigma^{\tfrac{d\beta}{4}}\norm{(-\Delta)^{-1}&\nabla\mathcal{P}(\rho_{0} u_t);\,(-\Delta)^{-1}\nabla(\rho_{0} u_t)}_{L^\infty(\R^d)}^4
\leqslant C \int_0^T\sigma^{\tfrac{d\beta}{4}}\norm{\rho_0 u_t}_{L^4(\R^d)}^{4\beta}\norm{\nabla(-\Delta)^{-1}(\rho_0 u_t)}_{L^2(\R^d)}^{4(1-\beta)}\nonumber\\
&\leqslant C_*\int_0^T\sigma^{\tfrac{d\beta}{4}} \norm{ \nabla  u_t}_{L^2(\R^d)}^{d\beta}\norm{ u_t}_{L^2(\R^d)}^{(4-d)\beta}
\norm{F,\nabla u}_{L^2(\R^d)}^{4(1-\beta)}\nonumber\\
&\leqslant C_*\left[\int_0^T \sigma \norm{ \nabla u_t}_{L^2(\R^d)}^4\right]^{\tfrac{d\beta }{4}}\left[\int_0^T \norm{ u_t}_{L^2(\R^d)}^{4}\right]^{(1-\tfrac{d}{4})\beta}\left[\int_0^T\norm{F,\;\nabla u}_{L^2(\R^d)}^{4}\right]^{1-\beta}
\end{align}
with 
\[
\beta=\dfrac{2d}{4+d}.
\]
Given that $\dot W^{1,\tfrac{d}{1-\alpha}}(\R^d)\hookrightarrow \dot \cC^{\alpha}(\R^d)$, we have: 
\begin{gather}
    \int_0^T\sigma^{r}  \norm{(-\Delta)^{-1}\nabla\mathcal{P}(\rho_{0} u_t);\,(-\Delta)^{-1}\nabla(\rho_{0} u_t)}_{\dot \cC^\alpha(\R^d)}^4
\leqslant C_* \int_0^T \sigma^{r} \norm{ u_t}_{L^{d/(1-\alpha)}(\R^d)}^4.
\end{gather}
For $d=2$, we have:
\begin{align*}
\int_0^T\sigma^\alpha \norm{ u_t}_{L^{2/(1-\alpha)}(\R^2)}^4& \leqslant C\int_0^T \sigma^\alpha \norm{\nabla  u_t}_{L^2(\R^2)}^{4\alpha}\norm{ u_t}_{L^2(\R^2)}^{4(1-\alpha)}\\
&\leqslant C \left[\int_0^T\sigma\norm{\nabla u_t}_{L^2(\R^2))}^4\right]^{\alpha}\left[\int_0^T\norm{ u_t}_{L^2(\R^2)}^{4}\right]^{1-\alpha}
\end{align*}
whereas for $d=3$, we only have for all $0<\alpha < 1/2$
\begin{align*}
\int_0^T\sigma^{2\alpha} \norm{ u_t}_{L^{3/(1-\alpha)}(\R^3)}^4&\leqslant  C\int_0^T \sigma^{2\alpha}\norm{\nabla u_t}_{L^2(\R^3)}^{2(1+2\alpha)} \norm{ u_t}_{L^2(\R^3)}^{2(1-2\alpha)}\\
&\leqslant C\sup_{[0,T]}\norm{ u_t}_{L^2(\R^3)}^{2(1-2\alpha)}\int_0^T \sigma^{2\alpha}\norm{\nabla u_t}_{L^2(\R^3)}^{2(1+2\alpha)}.
\end{align*}
For higher regularity, $1/2\leqslant \alpha<1$, we have:
\begin{align*}
\int_0^T\sigma^{\tfrac{2}{3}(5\alpha-1)} \norm{ u_t}_{L^{3/(1-\alpha)}(\R^3)}^4&\leqslant C  \int_0^T \sigma^{\tfrac{2}{3}(5\alpha-1)}\norm{\nabla u_t}_{L^4(\R^3)}^{\tfrac{16}{3}(\alpha-\tfrac{1}{2})} \norm{ \nabla u_t}_{L^2(\R^3)}^{\tfrac{4}{3}(5-4\alpha)}\\
&\leqslant C \left[\norm{\sqrt{\sigma }\nabla u_t}_{L^\infty((0,T),L^2(\R^3))}^2\right]^{\tfrac{2}{3}(2-\alpha)}\\
&\times \left[\int_0^T \sigma^{2}\norm{\nabla u_t}_{L^4(\R^3)}^{8/3}\right]^{2\alpha-1}\left[\int_0^T\norm{\nabla u_t}^2_{L^2(\R^3)}\right]^{2(1-\alpha)}.
\end{align*}
In conclusion, choosing $r$ as: 
\begin{gather}\label{c5.1}
r=\begin{cases}
                \max\left(1/3,\,2\alpha\right)\quad & \text{if}\quad d=2,\\
                \max\left(\dfrac{9}{14},\,2\alpha\right) \quad & \text{if}\quad d=3,\; \alpha\in (0,1/2)\\
                \dfrac{2}{3}(5\alpha-1) \quad &  \text{if }\quad d=3,\; \alpha\in [1/2,1).
            \end{cases}
\end{gather}
we have, owing to Young's inequality:
\begin{align}
    \int_0^T\sigma^r\norm{(-\Delta)^{-1}\nabla\mathcal{P}\left(\rho_{0}u_t\right),&\;(-\Delta)^{-1}\nabla\left(\rho_{0}u_t\right) }_{\cC^\alpha(\R^d)}^4\leqslant C_* \norm{F,\, \nabla u,\, u_t,\,\sqrt\sigma\nabla u_t}_{L^\infty((0,T),L^2(\R^d))}^4\nonumber\\
         &+C_*\left[\int_0^T \norm{F,\, \nabla u,\,u_t,\,  \nabla u_t}_{L^2(\R^d)}^2\right]^2+C_*\left[\int_0^T \sigma^{2}\norm{\nabla u_t}_{L^4(\R^d)}^{8/d}\right]^{d/2}.\label{c3.50}
\end{align}
The second term is discontinuous at the interface, and we have:
 \begin{align*}
 \left\|\nabla^2 (-\Delta)^{-1}\left[\dfrac{1}{2\widetilde\mu+\lambda}(-\Delta)^{-1}\dvg(\rho_0 u_t)\right]\right\|_{\cC^\alpha_{pw,\mathcal{C}_0}(\R^d)}&\leqslant  C\left(\norm{\lambda}_{\dot \cC^\alpha_{pw,\mathcal{C}_0}(\R^d)}+\ell_{\vph_0}^{-\tfrac{\alpha}{2}}\right)\norm{(-\Delta)^{-1}\nabla(\rho_0 u_t)}_{L^q(\R^d)}\\
 &+\left(1+\norm{\lambda}_{\dot \cC^\alpha_{pw,\mathcal{C}_0}(\R^d)}+\mathfrak{P}_{\mathcal{C}_0}\right)\norm{(-\Delta)^{-1}\nabla(\rho_0 u_t)}_{\cC^{\alpha}(\R^d)}.
 \end{align*}
 Finally, we interpolate the $L^q(\R^d)$ norm above between $L^2(\R^d)$ and $L^\infty(\R^d)$, and we recall \eqref{c3.50}  to obtain:
\begin{align}
    \int_0^T\sigma^r\left\|\frac{1}{\widetilde\mu}(-\Delta)^{-1}\nabla\mathcal{P}\left(  \rho_{0}u_t\right),\,\right.&\left.\,\nabla^2 (-\Delta)^{-1}\left[\dfrac{1}{2\widetilde\mu+\lambda}(-\Delta)^{-1}\dvg(\rho_0 u_t)\right]\right\|^4\nonumber\\
    &\leqslant  C_{**} \norm{F,\, \nabla u,\, u_t,\,\sqrt\sigma\nabla u_t}_{L^\infty((0,T),L^2(\R^d))}^4\nonumber\\
         &+C_{**}\left[\int_0^T \norm{F,\, \nabla u,\,u_t,\,  \nabla u_t}_{L^2(\R^d)}^2\right]^2+C_{**}\left[\int_0^T \sigma^{2}\norm{\nabla u_t}_{L^4(\R^d)}^{8/d}\right]^{d/2}.\label{c3.51}
\end{align}
\paragraph*{\textbf{Step 4: Closing of all estimates}}
    It follows from \eqref{c3.53},  \eqref{c3.52} and \eqref{c3.51} that: 
    \begin{align*}
        \int_0^T\sigma^r&\norm{\nabla u}_{\cC^\alpha_{pw,\mathcal{C}_0}(\R^d)}^4\leqslant C_{**} \norm{F,\, \nabla u,\, u_t,\,\sqrt\sigma\nabla u_t}_{L^\infty((0,T),L^2(\R^d))}^4+C_{**}\left[\int_0^T \norm{F,\, \nabla u,\,u_t,\,  \nabla u_t}_{L^2(\R^d)}^2\right]^2\nonumber\\
        &+C_{**}\left[\int_0^T \sigma^{2}\norm{\nabla u_t}_{L^4(\R^d)}^{8/d}\right]^{d/2}+C_{**}\int_0^T\sigma^r\left( \norm{\nabla u}_{L^q(\R^d)}^4+\norm{F}_{\cC^\alpha_{pw,\mathcal{C}_0}(\R^d)}^4\right)\\
        &+C\left( \mathfrak{P}_{\mathcal{C}_0} +\ell^{-\alpha}_{\vph_0}\right)^4\sup_{[0,T]}\left[\norm{\llbracket\mu\rrbracket}_{L^\infty(\mathcal{C}_0)}^4+\norm{\llbracket\mu\rrbracket,\, \llbracket \lambda\rrbracket}_{L^\infty(\mathcal{C}_0)}^4\left\|\dfrac{\scalar{\mu-\widetilde\mu}}{\scalar{\mu}}\right\|_{L^\infty(\mathcal{C}_0)}^4\right]\int_0^T\sigma^r\norm{\nabla u}_{L^\infty(\R^d)}^4\\ \vspace{0,3cm}
         &+C\sup_{[0,T]}\left[1+\norm{\lambda}_{\dot \cC^\alpha_{pw,\mathcal{C}_0}(\R^d)}+\left( \mathfrak{P}_{\mathcal{C}_0} +\ell^{-\alpha}_{\vph_0}\right)\norm{\llbracket\lambda\rrbracket}_{L^\infty(\mathcal{C}_0)}\right]^4\sup_{[0,T]}\norm{\mu-\widetilde\mu}_{\cC^\alpha_{pw,\mathcal{C}_0}(\R^d)}^4\int_0^T\sigma^r\norm{\nabla u}_{\cC^{\alpha}_{pw,\mathcal{C}_0}(\R^d)}^4.
    \end{align*}
    Therefore there exists a constant $[\mu]>0$ depending only on $d,\alpha$ and $\widetilde\mu$ such that if:
    \begin{align}
    \sup_{[0,T]}\left[1+\norm{\lambda}_{\dot \cC^\alpha_{pw,\mathcal{C}_0}(\R^d)}\right.&\left.+\left( \mathfrak{P}_{\mathcal{C}_0} +\ell^{-\alpha}_{\vph_0}\right)\norm{\llbracket\lambda\rrbracket}_{L^\infty(\mathcal{C}_0)}\right]\sup_{[0,T]}\norm{\mu-\widetilde\mu}_{\cC^\alpha_{pw,\mathcal{C}_0}(\R^d)}\nonumber\\
    &+\left( \mathfrak{P}_{\mathcal{C}_0} +\ell^{-\alpha}_{\vph_0}\right)\sup_{[0,T]}\left[\norm{\llbracket\mu\rrbracket}_{L^\infty(\mathcal{C}_0)}+\norm{\llbracket\mu\rrbracket,\, \llbracket \lambda\rrbracket}_{L^\infty(\mathcal{C}_0)}\left\|\dfrac{\scalar{\mu-\widetilde\mu}}{\scalar{\mu}}\right\|_{L^\infty(\mathcal{C}_0)}\right]\leqslant [\mu],
    \end{align}
    then the last two terms can be absorbed in the left hand side and interpolating the $L^q(\R^d)$ norm of the velocity gradient between $L^2(\R^d)$ and $L^\infty(\R^d)$ followed by Young's inequality, we find \eqref{eq3.49}.
\enddem
This completes the first part of the proof of \cref{Theorem 3.1} concerning the a priori estimates. We now move to the second part in the proof of \cref{Theorem 3.1}, focusing on the homotopy argument.

For $\theta\in [0,1]$, we define:
\[
 \mu^\theta= (1-\theta)\widetilde\mu +\theta \mu, \quad \text{ and } \quad \lambda^\theta= (1-\theta)\widetilde\lambda +\theta \lambda
\]
and we consider \eqref{ep1.6} where the viscosities $\mu,\,\lambda$ are replaced by $\mu^\theta$ and $\lambda^\theta$ respectively :
\begin{gather}\label{eq3.19}
    \begin{cases}
        \rho_0u_t -\dvg(2\mu^\theta\D u+\lambda^\theta\dvg uI_d)= \dvg F,\\
        u_{|t=0}=u_0.
    \end{cases}
\end{gather}
In order to ease the reading, we set $Y_T^\theta:=Y_T(\mu^\theta_0, \lambda^\theta_0)$. The principle of homotopy argument is to look for $\theta \in [0,1]$ such that \eqref{eq3.19} is uniquely solvable in $X_T$. Therefore, we define the following set:
\begin{gather}\label{c4.33}
\mathscr{E}:=\big\{\theta \in [0,1]\colon \forall (F,u_0)\in Y_T^\theta,\, \exists\,! \;u\in X_T \text{ solution of }\eqref{eq3.19}\big\}.
\end{gather}
It is obvious that \cref{Theorem 3.1} follows as soon as $1\in \mathscr{E}$. This will be obtained by proving that $\mathscr E$ is at the same time a non-empty, open and closed set of $[0,1]$. 
The proof is structured in two subsections: In the first one (see \cref{proof1}), we prove that $0\in \mathscr E$. In the second one (see \cref{proof2}), we show that $\mathscr E$ is both an open and a closed subset of $[0,1]$.  

\subsubsection{\textbf{Proof of $\mathscr{E}\neq \emptyset$}}\label[section]{proof1}

\paragraph{\textbf{Sketch of the proof}}

The aim is  to prove that $0\in \mathscr{E}$. 
When $\theta=0$\label{page}, the equations \eqref{eq3.19} reduces to:
\begin{gather}\label{eq3.37}
\begin{cases}
\rho_0 u_t -\widetilde\mu\Delta u-(\widetilde\mu+\widetilde\lambda) \nabla \dvg u=\dvg F,\\
u_{|t=0}= u_0
\end{cases}
\end{gather}
and the purpose is to prove that for all $(F,u_0)\in Y_T^0$ there exists a unique solution of \eqref{eq3.37} in $X_T$. 
 Solving \eqref{eq3.37} is not part of the classical theory for linear parabolic systems, as the initial density in front of $\dpt u$ is rough. Despite, one can add 
 $(\widetilde\rho-\rho_0) u_t$ on both sides of \eqref{eq3.37}, in order to obtain a heat equation with a right hand side that belongs to $L^2((0,T),L^2(\R^d))+ L^2((0,T),\dot H^{-1}(\R^d))$. As in \cite{danchin2020well}, we may require  small fluctuation of the density.
 In order to get around this difficulty, the strategy adopted is that, we consider a mollifying sequence of the initial density $\rho_0$, velocity $u_0$ and  $F$ respectively $\rho_{0,\delta}$ $u_{0,\delta}$, and $F_\delta$
and consider the equation: 
\begin{gather}\label{eq3.22}
\begin{cases}
\rho_{0,\delta} u_t -\widetilde\mu\Delta u-(\widetilde\mu+\widetilde\lambda) \nabla \dvg u=\dvg F_\delta,\\
u_{|t=0}= u_{0,\delta}.
\end{cases}
\end{gather}
The existence and uniqueness of global solution $u^\delta$ of the above system is based on, a change of unknowns, a study of the heat kernel and a homotopy argument. Taking the advantage that for any  $\delta>0$, $\rho_{0,\delta}$ is smooth enough and bounded away from vacuum, one has that $u$ solves \eqref{eq3.22} if and only if $v:= \rho_{0,\delta}u$ solves the following Cauchy problem: 
\begin{gather}\label{eq3.33}
\begin{cases}
v_t -\dvg \left(\dfrac{\widetilde\mu}{\rho_{0,\delta}}\D v\right)-\nabla\left(\dfrac{\widetilde\lambda}{\rho_{0,\delta}} \dvg v\right)&=\dvg F_\delta - \widetilde\mu\dvg \left( \rho_{0,\delta}^{-2} \left( \nabla \rho_{0,\delta}\otimes v+ v\otimes\nabla \rho_{0,\delta}\right)\right)\\
&-\widetilde\lambda\nabla\left(\rho_{0,\delta}^{-2} v\cdot \nabla \rho_{0,\delta} \right),\\
v_{|t=0}&= \rho_{0,\delta}u_{0,\delta}.
\end{cases}
\end{gather}
In order to ease computations, we consider the following problem for $\mu'_0,\lambda'_0,b_0=\mu'_0(x),\lambda'_0(x),b_0(x)$
\begin{gather}\label{eq3.23}
\begin{cases}
    v_t -\dvg (2\mu_0'\D v+\lambda_0'\dvg v I_d)= f+ \dvg ( b_0\cdot v),\\
    v_{|t=0}=v_0
\end{cases}
\end{gather}
and we have the following:
\begin{theo}\label[theo]{th2} Let $0< T< \infty$ and $f\in L^2((0,T)\times\R^d)$, $ 0<\mu_*\leqslant\mu_0',\,\lambda_0'\in W^{1,\infty}(\R^d)$, $b_0\in W^{1,\infty}(\R^d)$, and finally $v_0\in H^1(\R^d)$. Then there exists a unique solution 
    $v\in \cC([0,T], H^1(\R^d))$ with $v_t,\nabla^2 v\in L^2 ((0,T), L^2(\R^d))$ for the system \eqref{eq3.23}. Moreover, there exist  constants $\Lambda,\, K>0$ that depend on the norms of $\mu'_0,\,\lambda_0',\, b_0$ and the lower and upper bounds $\mu_*,\mu^*$ such that 
    \[
    \sup_{t\in [0,T]} e^{-2\Lambda t}\left\{\norm{v(t)}_{H^1(\R^d)}^2+\int_0^t \norm{v_t,\, \nabla v,\, \nabla^2 v}_{L^2(\R^d)}^2\right\}\leqslant K \left(\norm{v_0}^2_{H^1(\R^d)}+\int_0^T \norm{f(s)}_{L^2(\R^d)}^2ds\right).
    \]
\end{theo}
The proof of \cref{th2} is the purpose of the Step below, and it implies immediately:
\begin{lemm}\label[lemma]{lemA3}
    For all $\delta>0$, there exists a unique classical solution $u^\delta$ of the Cauchy problem  \eqref{eq3.22}. The sequence $(u^\delta)_\delta$ converges to $u\in X_T$ as $\delta \to 0$, the unique solution of  the Cauchy problem \eqref{eq3.37}.
\end{lemm}
All these discussions establish that $\mathscr{E}$ is a non-empty set. In the following, we present the proof of \cref{th2} and \cref{lemA3}.

\subparagraph{\textbf{Detailed proof of \cref{th2} and \cref{lemA3}}}
We start with the proof of \cref{th2}. 

\dem[\textbf{Proof of \cref{th2}}]
The proof of \cref{th2} is based on a homotopy argument. So, we define 
\[
{\mu_0'}^\theta= (1-\theta)\widetilde {\mu} + \theta \mu_0'\quad \text{and}\quad {\lambda_0'}^\theta= (1-\theta)\widetilde {\lambda} + \theta \lambda_0',
\]
and we consider the equation \eqref{eq3.23} where the viscosity $\mu_0',\lambda_0'$ are replaced by ${\mu_0'}^\theta,\, {\lambda_0'}^\theta$ respectively. Also, we define 
\[
X_T':= \{ v\in \cC([0,T], H^1(\R^d))\colon v_t,\, \nabla^2 v\in L^2((0,T)\times\R^d)\}
\]
and the map:
\begin{align*}
\Pi_\theta\colon \quad &X_T'\quad \to L^2((0,T)\times\R^d)\times H^1(\R^d)\\
& v\quad\quad\mapsto \left( v_t-\dvg (2{\mu_0'}^\theta\D v+{\lambda_0'}^\theta\dvg v I_d)-\dvg ( b_0\cdot v),\; v_{|t=0}\right)
\end{align*}
and finally we consider the set:
\[
\mathscr E':=\{\theta\in [0,1]\colon \Pi_\theta \text{ is one to one}\}.
\]
Obviously, the map $\Pi_\theta$ is well-defined and continuous, we will only have to show that $1\in \mathscr E'$ to conclude the 
proof of \cref{th2}.
To obtain $0\in \mathscr E'$, we observe that for $\theta=0$, the equation $\eqref{eq3.23}_\theta$ reduces to:
\begin{gather}\label{eq3.24}
v_t- \widetilde{\mu} \Delta v-(\widetilde{\mu}+\widetilde{\lambda})\nabla \dvg v= f+\dvg (b_0\cdot v).
\end{gather}
The existence of a unique solution $v$ of \eqref{eq3.24} is established in \cref{app2}, resulting in $0\in \mathscr E'$.
The remainder of the proof of \cref{th2} consists in  first deriving a priori estimates for the solution $v\in X'_T$ of \eqref{eq3.23}. Then, we prove that  $\mathscr E'$ is at the same time open and closed. 
\subparagraph{\textbf{A priori estimates}}
We multiply \eqref{eq3.23} by $v$ and integrate in time and space in order to obtain:
\begin{align*}
\dfrac{1}{2}\norm{v(t)}_{L^2(\R^d)}^2&+\int_0^t\int_{\R^d} 2\mu_0'\abs{\D v}^2+\int_0^t\int_{\R^d}\lambda_0'\abs{\dvg v}^2=\dfrac{1}{2}\norm{v_0}_{L^2(\R^d)}^2+\int_0^t\int_{\R^d} f v
-\int_0^t\int_{\R^d}b_0\cdot v \nabla v\\
&\leqslant \eta\int_0^t  \norm{\nabla v}_{L^2(\R^d)}^2+\dfrac{1}{2}\norm{v_0}_{L^2(\R^d)}^2+\int_0^t\norm{f}_{L^2(\R^d)}^2+ C_\eta \int_0^t(1+\norm{b_0}_{L^\infty(\R^d)}^2)\norm{v}_{L^2(\R^d)}^2)
\end{align*}
for some positive constant $\eta$ small related to the lower bound of the viscosity, 
in such way that the first term of the right hand side can be absorbed in the left hand side resulting in:
\begin{gather}\label{eq3.31}
\norm{v(t)}_{L^2(\R^d)}^2+\int_0^t\norm{\nabla v}^2_{L^2(\R^d)}\leqslant C_* 
\norm{v_0}_{L^2(\R^d)}^2+C_*\int_0^t\norm{f}_{L^2(\R^d)}^2+  C_*\int_0^t(1+\norm{b_0}_{L^\infty(\R^d)}^2)\norm{v}_{L^2(\R^d)}^2,
\end{gather}
where $C_*$ depends only on the bounds of the viscosity.  

Next, we use $\Delta v$ as a test function in  \eqref{eq3.23}.  The first term is:
\[
\int_{\R^d}v_t \Delta v=-\int_{\R^d} \nabla v_t \cdot \nabla v=-\dfrac{1}{2}\dfrac{d}{dt}\int_{\R^d}\abs{\nabla v}^2
\]
and the second one is:
\begin{align*}
    -2\int_{\R^d} \partial_k (\mu_0'\D ^{jk} v) \Delta v^j &-\int_{\R^d}\partial_j(\lambda_0'\dvg v)\Delta v^j
                     = -\int_{\R^d} \mu_0'\abs{\Delta v}^2 -\int_{\R^d} \mu_0' \partial_j \dvg v \Delta v^j-\int_{\R^d}\partial_k \mu_0' \D^{jk} v\Delta v^j\\
                     &-\int_{\R^d} \lambda_0' \abs{\nabla \dvg v}^2- \int_{\R^d}\dvg v \nabla \lambda_0'\cdot \nabla \dvg v\\
                     &=-\int_{\R^d} \mu_0'\abs{\Delta v}^2-\int_{\R^d}(\mu_0'+\lambda_0')\abs{\nabla \dvg v}^2+\int_{\R^d}\dvg v\partial_j \mu_0'\Delta v^j\\
                     &-\int_{\R^d} \dvg v\nabla \mu_0'\cdot \nabla \dvg v-\int_{\R^d}\partial_k \mu_0' \D^{jk} v\Delta v^j
                     - \int_{\R^d}\dvg v \nabla \lambda_0'\cdot \nabla \dvg v.
\end{align*}
Gathering all of theses computations, we have:
\begin{align*}
\dfrac{1}{2}\dfrac{d}{dt}\norm{\nabla v}_{L^2(\R^d)}^2&+\int_{\R^d} \mu_0'\abs{\Delta v}^2+\int_{\R^d}(\mu_0'+\lambda_0')\abs{\nabla \dvg v}^2=
\int_{\R^d} f^j \Delta v^j+\int_{\R^d} \dvg ( b_0\cdot v)\Delta v^j \\
+&\int_{\R^d}\int_{\R^d}\dvg v\partial_j \mu_0'\Delta v^j
                     -\int_{\R^d} \dvg v\nabla \mu_0'\cdot \nabla \dvg v-\int_{\R^d}\partial_k \mu_0' \D^{jk} v\Delta v^j
                     - \int_{\R^d}\dvg v \nabla \lambda_0'\cdot \nabla \dvg v.
\end{align*}
Thus, H\"older's and Young's inequalities yield:
\begin{gather}\label{eq3.32}
\norm{\nabla v(t)}_{L^2(\R^d)}^2+\int_0^t \norm{\Delta v}_{L^2(\R^d)}^2\leqslant C_* \norm{\nabla v_0}_{L^2(\R^d)}^2
+C_*\int_0^t \norm{f}_{L^2(\R^d)}^2
+C_* \int_0^t \norm{\nabla v}_{L^2(\R^d)}^2 \norm{\nabla\mu_0',\nabla\lambda_0',b_0,\nabla b_0}_{L^\infty(\R^d)}^2.
\end{gather}
Finally, we gather \eqref{eq3.31} and \eqref{eq3.32}, and we make use of Gr\"onwall's lemma and  \eqref{eq3.23} to obtain:
\begin{gather*}
\sup_{t\in [0,T]} e^{-\Lambda t}\left\{\norm{v(t)}_{H^1(\R^d)}^2+\int_0^t\norm{v_t,\,\nabla v,\, \nabla^2 v}_{L^2(\R^d)}^2\right\}\leqslant C_* \norm{v_0}_{H^1(\R^d)}^2
+C_*\int_0^T \norm{f}_{L^2(\R^d)}^2 
\end{gather*}
with
\[
\Lambda = C_*( 1+\norm{\nabla\mu_0',\nabla\lambda_0',b_0,\nabla b_0}_{L^\infty(\R^d)}^2).
\]
This ends the computations of the a priori estimates and we turn to the homotopy argument.
\subparagraph{\textbf{Closing of the homotopy argument}}
In this section, we will show that $\mathscr E'=[0,1]$. We have already showed that $\mathscr E'$ is not empty.  
 Also, we may notice that if $\Pi_\theta$ is invertible then, from the previous a priori estimates, $\Pi_\theta^{-1}$  is linear continuous with a norm that depends only on $\Lambda$ and time $T$. In particular, we have:
\[
K_T:=\sup_{\theta\in \mathscr E'}\norm{\Pi_\theta^{-1}}_{\mathcal{L}}<\infty.
\]
 Let us prove that $\mathscr E'$ is an open set. Let $\theta_0\in \mathscr E'$, then, $\Pi_{\theta_0}$ is a one to one map and  one can write,  for all $\theta\in [0,1]$, 
\[
\Pi_\theta= \Pi_{\theta_0}+ (\theta-\theta_0)(\Pi_1-\Pi_0)=\Pi_{\theta_0}\left( I+ (\theta-\theta_0)\Pi_{\theta_0}^{-1}(\Pi_1-\Pi_0)\right).
\]
Thus, all $\theta$ such that $2\abs{\theta-\theta_0}\norm{\Pi_{\theta_0}^{-1}(\Pi_1-\Pi_0)}_{\mathcal{L}(X_T')}<1$ belongs to $\mathscr E'$, therefore $\mathscr E'$ is an open set of $[0,1]$. We turn to the last part. Let $(\theta_n)_n$ be a sequence of element of $\mathscr E'$ that converges to some $\theta\in  [0,1]$.  As above, one can write,  for some $n_0$
\[
\Pi_\theta= \Pi_{\theta_{n_0}}\left( I+ (\theta-\theta_{n_0})\Pi_{\theta_{n_0}}^{-1}(\Pi_1-\Pi_0)\right).
\]
Since $K_T<\infty$, one can choose large $n_0$ such that 
\[
2\abs{\theta-\theta_{n_0}} K_T\norm{\Pi_1-\Pi_0}_{\mathcal{L}}<1
\]
thus $\Pi_{\theta}$ is one to one as composition of two one to one operators. Therefore $\theta \in \mathscr E'$ and consequently $\mathscr E'$ is a closed set of $[0,1]$. One then concludes that $1\in \mathscr E'$ and \cref{th2} then follows.
\enddem

\dem[\textbf{Proof of \cref{lemA3}}]
One deduces from \cref{th2}, \eqref{eq3.33} and then $\eqref{eq3.22}$ that the Cauchy problem \eqref{eq3.22} admits a unique solution $u^\delta\in \cC([0,T], H^1(\R^d))$ with $u_t^\delta, \nabla ^2u^\delta\in L^2((0,T)\times\R^d)$. Moreover, by differentiating  \eqref{eq3.22} with respect to space, and noting that $\rho_0^\delta$ and $u_0^\delta$ are smooth, we obtain that for all $s\geqslant 2$, $
u^\delta \in \cC([0,T], H^s(\R^d)),$ and $u_t\,,\nabla ^2 u^\delta \in L^2((0,T), H^{s-1}(\R^d))$.  In the same manner, we differentiate
\eqref{eq3.22} with respect to time, in order to obtain the equation satisfied by $u_t^\delta $. Of course, $F_t^\delta$ is regular enough and the regularity $u_t^\delta \in L^2((0,T), H^{s-1}(\R^d))$ is adequate for the differentiating to make sense. Thus, applying again \cref{th2} one deduces from the previous discussion that $u_t^\delta \in \cC ([0,T], H^s(\R^d))$
and $u_{tt}^\delta,\nabla^2 u_t^\delta \in L^2((0,T), H^{s-1}(\R^d))$. This closes the existence of classical solution 
of the Cauchy problem \eqref{eq3.22}.

  The sequence  $(u^\delta)_\delta$  verifies the estimates for $u$, $u_t$, (in \eqref{ep1.5}) $u_{tt}$ (in \eqref{eq3.34})  and the $L^p$ estimates in \cref{Theorem 3.1} for $u$ and $u_t$, uniformly with respect to the regularisation parameter $\delta>0$ with the exception of the piecewise H\"older estimate.  Hence, one can let $\delta$ go to zero in order to obtain a solution of the Cauchy problem \eqref{eq3.37}. Since $(u^\delta)_\delta $ is bounded in $L^\infty((0,T), H^1(\R^d))$ and its time derivative $(u_t^\delta)_\delta$ is bounded in $L^2((0,T), H^1(\R^d))$, it then turns out  that $u\in \cC([0,T], H^1(\R^d))$. On the other hand, given that 
$u_t\in L^\infty((0,T),L^2(\R^d))$ and  $\sigma u_{tt}\in L^2((0,T), L^2(\R^d))$ it immediately follows  $u_t\in \cC((0,T], L^2(\R^d))$. It only remains the continuity of $u_t$ at the initial time.   To achieve this, we first observe that for all $\tau>0$:
\begin{gather*}
       \dfrac{1}{2}\int_{\R^d}\rho_{0,\delta}\abs{ u_t^\delta(\tau)}^{2}+\int_{0}^{\tau}
\widetilde\mu\int_{\R^d}  \abs{\nabla  u_t^\delta}_{L^2(\R^d)}^{2}+(\widetilde\mu+\widetilde\lambda)\int_{0}^{\tau}\int_{\R^d} (\dvg u_t^\delta)^{2}
=\dfrac{1}{2}\int_{\R^d}\rho_{0,\delta} \abs{{u_t^\delta}_{|t=0}}^{2}- \int_0^{\tau}\int_{\R^d} \dpt F^{jk}_\delta\partial_k({u^\delta})^j_t,
\end{gather*}
and hence:
\[
\dfrac{1}{2}\int_{\R^d}\rho_{0,\delta}\abs{ u_t^\delta(\tau)}^{2}
\leqslant \dfrac{1}{2}\int_{\R^d}\rho_{0,\delta} \abs{{u_t^\delta}_{|t=0}}^{2}+ C_* \int_0^{\tau}\norm{\dpt F}_{L^2(\R^d)}^2.
\]
Using the compact embedding of  $H^1(\R^d)$ into $L^2_\loc (\R^d)$ and the Aubin-Lions Lemma, we obtain $u^\delta_t(\tau)\to u_t(\tau)$ a.e in $\R^d$. Consequently, applying Fatou's Lemma yields:
\[
\dfrac{1}{2}\int_{\R^d}\rho_{0}\abs{ u_t(\tau)}^{2}
\leqslant \dfrac{1}{2}\int_{\R^d}\rho_{0} \abs{{u_t}_{|t=0}}^{2}+ C_* \int_0^{\tau}\norm{\dpt F}_{L^2(\R^d)}^2,
\]
from which we deduce:
\[
\limsup_{\tau\to 0} \norm{\sqrt{\rho_0} u_t(\tau)}_{L^2(\R^d)}\leqslant \norm{\sqrt{\rho_0} {u_t}_{|t=0}}_{L^2(\R^d)}.
\]
Since
\[
\norm{\sqrt{\rho_0} {u_t}_{|t=0}}_{L^2(\R^d)}\leqslant \liminf_{\tau\to 0} \norm{\sqrt{\rho_0} u_t(\tau)}_{L^2(\R^d)},
\]
 we infer that $\sqrt{\rho_0} u_t(\tau)$ convergences  to $\sqrt{\rho_0} {u_t}_{|t=0}$ strongly in $L^2(\R^d)$ and consequently we obtain the continuity  of $u_t$ at the initial time given that   $\rho_0$ is bounded away 
from vacuum. 

This solution is unique: the difference $\delta u$ of two solutions $u_{1}$ and $u_2$  of \eqref{eq3.37} verifies 
\[
\begin{cases}
\rho_0 \delta u_t-\widetilde\mu\Delta \delta u-(\widetilde\mu+\widetilde\lambda)\nabla\delta u=0,\\
\delta u_{|t=0}=0
\end{cases}
\]
and testing the above equation with $\delta$ leads to $\delta u=0$. 

From \eqref{eq3.37}, it is obvious to derive the following expression for $\nabla u$:
\[
\nabla u= -\frac{1}{\widetilde\mu}(-\Delta)^{-1}\nabla\left[  \mathcal{P}\left(  \rho_{0}u_t\right)  -\mathcal{P}\dvg F  \right]  \nonumber\\
  -\dfrac{1}{2\widetilde \mu+\widetilde\lambda}(-\Delta)^{-1}\nabla\left[ \mathcal{Q}(
\rho_{0}u_t)  -\mathcal{Q}  \dvg F  \right].
\]
where $\mathcal{P}$ is the Leray projector onto the space of divergence-free vector fields and $\mathcal{Q}$ is the projector onto the space of curl-free vector fields. 
Finally, we can easily show that \eqref{eq3.49} also holds for $\nabla u$.  
\enddem

This concludes the proof that $\mathscr{E}$ is a nonempty set. The only remaining task is to prove that  $\mathscr{E}$ is both an open and  closed subset of $[0,1]$.   This will be addressed in the following step.


\subsubsection{\textbf{Remainder of the proof of \cref{Theorem 3.1}}}\label[section]{proof2}
We now proceed to the final step in the proof of \cref{Theorem 3.1}, which is to show that $\mathscr{E}$ is both an open and closed subset of $[0,1]$. 
\dem 
We begin by 
proving that $\mathscr{E}$ is an open subset of $[0,1]$. So let $\theta_0\in \mathscr E$, the goal is to find some $\eta>0$ such that 
\begin{gather}\label{eq3.53}
(\theta_0-\eta,\theta_0+\eta)\cap [0,1]\subset \mathscr E.
\end{gather}
For this purpose, let $\theta\in [0,1]$, $(F,u_0)\in Y_T^\theta$ and consider the Cauchy problem \eqref{eq3.19} which can be written as:
\begin{gather}\label{eq3.52}
\begin{cases}
\rho_0 u_t- \dvg(2\mu^{\theta_0}\D u+\lambda^{\theta_0}\dvg u I_d)=\dvg F+(\theta-\theta_0)\dvg \{2(\mu-\widetilde\mu)\D u+(\lambda-\widetilde\lambda)\dvg u I_d\},\\
u_{|t=0}=u_0.
\end{cases}
\end{gather}
We will apply a fixed-point argument to solve the above Cauchy problem with a small $\eta$ that does  not depend on $F$ nor $u_0$. In fact, the proof shows that the constant $\eta$ only depends on the norm of the viscosity and initial density as required in \eqref{ep1.2}, as well on the interface characteristics. The solution will be constructed in  $X_{T,u_0}$ which is defined as: 
\[
X_{T,u_0}:=\{ v\in X_T\colon v_{|t=0}=u_0\},
\]
a complete subset of the Banach space $X_T$. So for $v\in X_{T,u_0}$, we consider the following equation:
\begin{gather}\label{eq3.40}
\begin{cases}
\rho_0 u_t- \dvg(2\mu^{\theta_0}\D u+\lambda^{\theta_0}\dvg u I_d)=\dvg F+(\theta-\theta_0)\dvg \{2(\mu-\widetilde\mu)\D v+(\lambda-\widetilde\lambda)\dvg v I_d\},\\
u_{|t=0}=u_0.
\end{cases}
\end{gather}
The regularity of $v$ and $\mu,\lambda$ ensures that the right hand side written as $\dvg F(v)$ with 
\begin{gather}\label{eq3.38}
    F(v):= F+ (\theta-\theta_0)\{2(\mu-\widetilde\mu)\D v+(\lambda-\widetilde\lambda)\dvg v I_d\}
\end{gather}
verifies $(F(v), u_0)\in Y_T^{\theta_0}$. The last condition in the definition of $Y_T^{\theta_0}$ (see \eqref{c7.6}) is valid since $v_{|t=0}=u_0$.  Indeed, we have
\begin{align}
\dvg \Big(2\mu^{\theta_0}_0\D u_0&+\lambda^{\theta_0}_0 \dvg u_0I_d+F(v)_{|t=0}\Big)\notag\\
&=\dvg \Big(2\big(\mu^{\theta_0}_0+(\theta-\theta_0)(\mu_0-\widetilde\mu)\big)\D u_0+\big(\lambda_0^{\theta_0}+(\theta-\theta_0)(\lambda_0-\widetilde\lambda)\big)\dvg u_0 I_d+ F_{|t=0}\Big)\notag\\
&=\dvg \Big(2\mu^{\theta}_0\D u_0+\lambda^{\theta}_0 \dvg u_0I_d+F_{|t=0}\Big)\in L^2(\R^d)\; \text{ since }\; (F,u_0)\in Y^\theta_T.\label{c7.5}
\end{align}

Now, owing to the fact that $\theta_0\in \mathscr E$, the underlying equation admits a unique solution 
$\Phi(v)\in X_{T,u_0}$ that verifies all the estimates in \cref{Theorem 3.1}, where $F$ is replaced by $F(v)$. 

In the following we are looking for some small $\theta-\theta_0$ and some $R_T$ which depends on the initial data such that the condition $\norm{v}_{X_T}\leqslant R_T$
implies the same for $\Phi(v)$.   As a standard notation, from now on $A\lesssim B$ means $A\leqslant C B$ with a  constant $C$ depending on the  bounds of $\mu,\lambda,\rho_0$, interface characteristics, and on all of the norms of $\mu,\lambda$ as required in \eqref{ep1.2}.
\paragraph{\textbf{Step 1: Estimates for $\Phi(v)$ and $\partial_t \Phi(v)$}}
Since $\Phi(v)$ solves the Cauchy problem \eqref{eq3.40}, $(F(v),u_0)\in Y_T^{\theta_0}$, and $\theta_0\in \mathscr E$, 
one has from the a priori estimate the following estimate for $\Phi(v)$:
   \begin{align}
     \sup_{[0,T]}\norm{\Phi(v),\,\partial_t \Phi(v),\, \nabla \Phi(v)}_{L^2(\R^d)}^{2}&+\int_{0}^{T}
\norm{\nabla \Phi(v),\, \nabla \dpt \Phi(v),\, \dpt\Phi(v)}_{L^2(\R^d)}^{2}\lesssim \norm{u_0}_{H^1(\R^d)}^2 +\norm{\partial_t \Phi(v)_{|t=0}}_{L^2(\R^d)}^2\nonumber\\
&+ \sup_{[0,T]}\norm{F(v)}_{L^2(\R^d)}^2+\int_0^T \norm{\dpt F(v)}_{L^{2}(\R^d)}^2 +\int_0^T\norm{F(v)}_{L^2\cap L^6(\R^d))}^2.
\end{align}
On the one hand, from \eqref{eq3.38}, we infer that:
\begin{align}
\sup_{[0,T]}\norm{F(v)-F}_{ L^2(\R^d)}^2+\int_0^T \norm{F(v)-F}_{ L^2(\R^d)}^2 &\lesssim  \abs{\theta-\theta_0}^2\left\{\sup_{[0,T]}\norm{\nabla v}_{ L^2(\R^d)}^2+\int_0^T\norm{\nabla v}_{ L^2(\R^d)}^2\right\}\nonumber\\
&\lesssim \abs{\theta-\theta_0}^2 \norm{v}_{X_T}^2\label{eq3.42}
\end{align}
and 
\[
\int_0^T\norm{F(v)-F}_{L^6(\R^d))}^2
\lesssim \abs{\theta-\theta_0}^2 \norm{v}_{X_T}^2.
\]
On the other hand, from \eqref{eq3.38} one can express the time derivative of $F(v)$  as:
\begin{gather}\label{eq3.47}
\dpt F(v)=F_t+(\theta-\theta_0)\{2(\mu-\widetilde\mu)\D v_t+(\lambda-\widetilde\lambda)\dvg v_t I_d
+2\mu_t\D v+\lambda_t\dvg v I_d\}
\end{gather}
and consequently, it holds:
\begin{gather}\label{eq3.43}
\int_0^T \norm{\dpt F(v)-F_t}_{ L^2(\R^d)}^2\lesssim \abs{\theta-\theta_0}^2 \left\{\norm{\nabla v}_{L^4((0,T)\times\R^d)}^2+\int_0^T\norm{\nabla v_t}_{L^2(\R^d)}^2\right\}
\lesssim \abs{\theta-\theta_0}^2  \norm{v}_{X_T}^2.
\end{gather}
Gathering all of these estimates, one has the following:
\begin{gather}\label{eq3.41}
     \sup_{[0,T]}\norm{\Phi(v),\,\partial_t \Phi(v),\, \nabla \Phi(v)}_{L^2(\R^d)}^{2}+\int_{0}^{T}
\norm{\nabla \Phi(v),\, \nabla \dpt \Phi(v),\, \dpt\Phi(v)}_{L^2(\R^d)}^{2}
\lesssim \norm{(F,u_0)}_{Y_T^{\theta}}^2 +\abs{\theta-\theta_0}^2  \norm{v}_{X_T}^2.
\end{gather}
Also, one has the following estimate for $\nabla \Phi(v)$:
\begin{gather}\label{eq3.45}
\norm{\nabla\Phi(v)}_{L^\infty((0,T),L^6(\R^d))} \lesssim \norm{F(v)}_{L^\infty((0,T),L^6(\R^d))}+ \sup_{[0,T]}\norm{F(v), \dpt\Phi(v), \nabla \Phi(v)}_{L^2(\R^d)}.
\end{gather}

Owing to the fact that the estimate \eqref{eq3.41} provides a bound for $\dpt \Phi(v)$ in $L^\infty((0,T), L^2(\R^d))$ and for $\Phi(v)$ in $L^\infty((0,T), H^1(\R^d))$ and owing to the estimate \eqref{eq3.42}, we only need to estimate the $L^\infty((0,T),L^6(\R^d))$ norm of $F(v)$. In fact, this can be easily obtained from \eqref{eq3.38} as:
\[
\norm{F(v)-F}_{L^\infty((0,T),L^6(\R^d))} \lesssim \abs{\theta-\theta_0} \norm{\nabla v}_{L^\infty((0,T),L^6(\R^d))}
\lesssim \abs{\theta-\theta_0} \norm{v}_{X_T}.
\]
Gathering all of theses estimates, one has: 
\begin{gather}\label{eq3.44}
\norm{\nabla\Phi(v)}_{L^\infty((0,T),L^6(\R^d))} \lesssim  \norm{(F,u_0)}_{Y_T^{\theta}}+ \abs{\theta-\theta_0} \norm{v}_{X_T}.
\end{gather}
We turn to the $L^{16}((0,T), L^8(\R^d))$ estimate for $\nabla \Phi(v)$. For this purpose, we recall from \eqref{eq3.36} that:
\[
\int_0^T\norm{\nabla \Phi(v)}_{L^8(\R^d)}^{16} \lesssim  \int_0^T\norm{ F(v)}_{L^8(\R^d)}^{16}+ \norm{F(v), \dpt \Phi(v), \nabla \Phi(v)}_{L^2\cap L^\infty((0,T), L^2(\R^d))}^{16}+\norm{\nabla \dpt \Phi(v)}_{L^2((0,T)\times \R^d)}^{16}.
\]
We have already obtained  bound for the last  two terms above in \eqref{eq3.41} and \eqref{eq3.42} and we have only to estimate the first one. From \eqref{eq3.38}, we have:
\[
\int_0^T\norm{F(v)-F}_{L^8(\R^d))}^{16} \lesssim \abs{\theta-\theta_0}^{16} \int_0^T\norm{\nabla v}_{L^8(\R^d))}^{16}
\lesssim \abs{\theta-\theta_0}^p \norm{v}_{X_T}^{16}
\]
and hence:
\begin{gather}\label{eq4.7}
    \int_0^T\norm{\nabla \Phi(v)}_{L^8(\R^d)}^{16} \lesssim   \norm{(F,u_0)}_{Y_T^{\theta}}^p+\abs{\theta-\theta_0}^{16}  \norm{v}_{X_T}^{16}.
\end{gather}
\paragraph{\textbf{Step 2: Estimates for $\partial_{tt}\Phi(v)$ and $L^4$-norm estimate for $\nabla \dpt \Phi(v)$}}
As in the previous step, one first  recalls that $\Phi(v)$ verifies:
\begin{multline*}
     \int_0^T\norm{ \sqrt\sigma \partial_{tt}\Phi(v),\; \sigma \nabla \partial_{tt}\Phi(v)}_{L^2(\R^d)}^2+\norm{\sqrt\sigma \dpt\nabla \Phi(v),\; \sigma \partial_{tt}\Phi(v)}^2_{L^\infty((0,T),L^2(\R^d))}\\
     \lesssim \int_0^{\sigma(T)}\norm{\dpt \nabla\Phi(v) }^2_{L^2(\R^d)}+\int_0^{T} \norm{\partial_tF(v),\, \sigma\partial_{tt}F(v)}_{L^2(\R^d)}^2
    +  \norm{\nabla \Phi(v)}_{L^\infty((0,T),L^4(\R^d))}^2.
\end{multline*}
All terms that appear in the right hand side above are controlled in \eqref{eq3.41}, \eqref{eq3.44}, \eqref{eq3.42} except the $L^2((0,T)\times \R^d)$ norm of $\sigma \partial_{tt} F(v)$ which is:
\[
\int_0^T\sigma^2\norm{\partial_{tt}F(v)}_{L^2(\R^d)}^2\leqslant \int_0^T\sigma^2\norm{F_{tt}}_{L^2(\R^d)}^2+\int_0^T\sigma^2\norm{\partial_{tt}F(v)-F_{tt}}_{L^2(\R^d)}^2.
\]
In order to link the last term to the norm of $v$, we first express $\partial_{tt}F(v)-F_{tt}$, from \eqref{eq3.47} as
\[
\partial_{tt}F(v)-F_{tt}=(\theta-\theta_0)\{2(\mu-\widetilde\mu)\D v_{tt}+(\lambda-\widetilde\lambda)\dvg v_{tt} I_d
+4\mu_t\D v_t+2\lambda_t\dvg v_t I_d+ 2\mu_{tt}\D v+\lambda_{tt}\dvg vI_d\}
\]
and then
\begin{align*}
    \int_0^T\sigma^2\norm{\partial_{tt}F(v)-F_{tt}}_{L^2(\R^d)}^2&\lesssim \abs{\theta-\theta_0}^2 \norm{v}_{X_T}^2
\end{align*}
from which, we deduce that 
\begin{align}\label{eq4.6}
\int_0^T\norm{ \sqrt\sigma \partial_{tt}\Phi(v),\; \sigma \nabla \partial_{tt}\Phi(v)}_{L^2(\R^d)}^2&+\norm{\sqrt\sigma \dpt\nabla \Phi(v),\; \sigma \partial_{tt}\Phi(v)}^2_{L^\infty((0,T),L^2(\R^d))}\nonumber\\
&\lesssim \norm{(F,u_0)}_{Y_T^{\theta}}^2+ \abs{\theta-\theta_0}^2\norm{v}_{X_T}^2 .
\end{align}
This ends the computations of the higher-order Hoff estimates. Now, we may turn to the $L^{8/d}((0,T),L^4(\R^d))$ norm estimate for $\sigma^{\tfrac{d}{4}}\nabla\dpt \Phi(v)$. Thanks to 
\eqref{eq3.48}, H\"older's and  Young's inequality, we have:
\begin{multline*}
\int_0^T\sigma^{2}\norm{\nabla \dpt \Phi(v)}_{L^4(\R^d)}^{8/d}\lesssim \int_0^T\sigma^{2}\norm{\dpt F(v)}_{L^4(\R^d)}^{8/d}+\left[\int_0^T\norm{\nabla\Phi(v)}_{L^8(\R^d)}^{16}\right]^{1/(2d)}+\norm{\sigma\partial_{tt}\Phi(v)}_{L^\infty((0,T),L^2(\R^d))}^{8/d}\\
       +\norm{\nabla \Phi(v)}_{L^4((0,T)\times\R^d)}^{8/d}+\left[\int_0^T\norm{\nabla \dpt \Phi(v),\, \dpt F(v)}_{L^2(\R^d)}^2\right]^{4/d}.
\end{multline*}
The estimates \eqref{eq4.6} and \eqref{eq4.7}, \eqref{eq3.44} and \eqref{eq3.43} control all of the terms that appear in the left hand side above, except the $L^{8/d}((0,T), L^4(\R^d))$-norm of $\sigma^{\tfrac{d}{4}}\dpt F(v)$. Starting from \eqref{eq3.47}, one has:
\begin{align*}
\int_0^T\sigma^{2}\norm{\dpt F(v)-F_t}_{L^4(\R^d)}^{8/d}&\lesssim\abs{\theta-\theta_0}^{8/d}\left\{\int_0^T\sigma^{2}\norm{\nabla v_t}_{L^4(\R^d)}^{8/d}+\left[\int_0^T\norm{\nabla v}_{L^8(\R^d)}^{16}\right]^{1/(2d)}\right\}\\
&\lesssim \abs{\theta-\theta_0}^{8/d}  \norm{v}_{X_T}^{8/d}.
\end{align*}
Gathering all of the estimates one obtains the following:
\[
\int_0^T\sigma^{2}\norm{\dpt \nabla \Phi(v)}_{L^4(\R^d)}^{8/d}\lesssim  \norm{(F,u_0)}_{Y_T^{\theta}}^{8/d}+ \abs{\theta-\theta_0}^{8/d}  \norm{v}_{X_T}^{8/d}.
\]
We now turn to the last step:
\paragraph{\textbf{Step 3: Piecewise H\"older-norm estimate for $\nabla \Phi(v)$}}
As above, thanks to \eqref{eq3.49},  the following H\"older estimate holds for $\nabla \Phi(v)$:
\begin{align*}
\int_0^T\sigma^{r}\norm{\nabla \Phi(v)}_{\cC^\alpha_{pw,\mathcal{C}_0}(\R^d)}^4 
         &\lesssim  \int_0^T\sigma^{r}\norm{F(v)}_{\cC^\alpha_{pw,\mathcal{C}_0}(\R^d)}^4+\left[\int_0^T \sigma^{4/d}\norm{\dpt \nabla \Phi(v)}_{L^4(\R^d)}^{8/d}\right]^{d/2} \\
         &+ \norm{F(v),\,  \nabla \Phi(v),\, \dpt \Phi(v),\,\sqrt\sigma\nabla \dpt \Phi(v)}_{L^\infty((0,T),L^2(\R^d))}^4\nonumber\\
         &+\left[\int_0^T \norm{F(v),\, \nabla \Phi(v),\,\dpt \Phi(v),\,  \dpt\nabla \Phi(v)}_{L^2(\R^d)}^2\right]^2.
\end{align*}
According to the above computations, one only needs to estimate the first term of the right hand side above and from \eqref{eq3.38} one has:
\[
\int_0^T\sigma^{r}\norm{ F(v)-F}_{\cC^\alpha_{pw,\mathcal{C}_0}(\R^d)}^4
\lesssim \abs{\theta-\theta_0}^4\norm{v}_{X_T}^4.
\]

Gathering all of these estimates, one has the following:
\begin{gather}\label{c3.35}
\norm{\Phi(v)}_{X_T}\leqslant C\left( R_0+ \abs{\theta-\theta_0} \norm{v}_{X_T} \right)
\end{gather}
where $R_0$ depends on the norm of $(F, u_0)\in Y_T^\theta$ and the constant $C$ that depends not only on the bounds of the  density and viscosity, on the interface characteristics, and all of the norm of the viscosity $\mu,\, \lambda$ and initial density $\rho_0$ as required in \eqref{ep1.2}. So, choosing large $R_T$, for instance $R_T=2 CR_0$ and $\abs{\theta-\theta_0}\leqslant \min\{1, 1/(2 C)\}$, one obtains $\norm{\Phi(v)}_{X_T}\leqslant R_T$ as soon as $\norm{v}_{X_T}\leqslant R_T$.  We thus define a map $\Phi$ from the close ball $B_{X_{T,u_0}}(0, R_T)$ of $X_{T,u_0}$ to itself.

\paragraph{\textbf{Step 4: Construction of the approximate sequence and proof of its convergence}}
We construct the sequence of approximate solutions $\big(u^n\big)_{n\in \N}$, inductively, as follows: 
\begin{gather*}
    u_{n+1}= \Phi(u_n), \; n\geqslant 0,
\end{gather*}
with 
\begin{gather}\label{c7.3}
\begin{cases}
\rho_0 \dpt u_{n+1}- \dvg(2\mu^{\theta_0}\D u_{n+1}+\lambda^{\theta_0}\dvg u_{n+1} I_d)&=\dvg F+(\theta-\theta_0)\dvg \big(2(\mu-\widetilde\mu)\D u_n\\
&+(\lambda-\widetilde\lambda)\dvg u_n I_d\big),\\
{u_{n+1}}_{|t=0}&=u_0.
\end{cases}
\end{gather}
It follows from the conclusion of Step 3 that the sequence  $(u_n)_n$  is uniformly bounded in $X_T$:
\begin{gather}\label{c7.2}
\sup_{n\in \N} \norm{ u_n}_{X_T}\leqslant R_T.
\end{gather}
In the following lines, we will prove that $(u_n)_n$ converges strongly in $L^\infty((0,T), L^2(\R^d))\cap L^2((0,T), H^1(\R^d))$ to some $u\in X_{T,u_0}$ which satisfies \eqref{eq3.52}.

To achieve this, we take the differences in \eqref{c7.3} and  obtain: 
\[
\begin{cases}
    \rho_0\dpt w_{n+1}-\dvg\Big(2\mu^{\theta_0}\D (w_{n+1})+\lambda^{\theta_0}\dvg (w_{n+1}) I_d\Big)= (\theta-\theta_0)\dvg\Big( 2(\mu-\widetilde\mu)\D (w_n)+(\lambda-\widetilde\lambda)\dvg (w_n) I_d\Big),\\
    {w_{n+1}}_{|t=0}=0,
\end{cases}
\]
where $w_{k+1}=u_{k+1}-u_k$, with $k\geqslant 1$. By multiplying the equation above by $w_{n+1}$, integrating in space, and applying Korn's inequality (see \eqref{c7.1}), along with Young's and H\"older's inequalities, we obtain:
\[
\dfrac{1}{2}\dfrac{d}{dt}\norm{\sqrt{\rho}_0 w_{n+1}}_{L^2(\R^d)}^2+\dfrac{\mu_*}{2}\norm{\nabla w_{n+1}}_{L^2(\R^d)}^2
\leqslant C\abs{\theta-\theta_0}^2\norm{\nabla w_n}_{L^2(\R^d)}^2.
\]
Time integration, along with the condition ${w_{n+1}}_{|t=0}=0$ yields:
\[
\sup_{[0,T]} \norm{w_{n+1}}_{L^2(\R^d)}^2+\int_0^T\norm{\nabla w_{n+1}}_{L^2(\R^d)}^2\leqslant C \abs{\theta-\theta_0}^2\Bigg(\sup_{[0,T]} \norm{w_{n}}_{L^2(\R^d)}^2+\int_0^T\norm{\nabla w_{n}}_{L^2(\R^d)}^2\Bigg),
\]
and therefore:
\begin{align*}
\sup_{[0,T]} \norm{u_{n+1}-u_n}_{L^2(\R^d)}^2+\int_0^T\norm{\nabla \big(u_{n+1}-u_n\big)}_{L^2(\R^d)}^2&\leqslant C^n \abs{\theta-\theta_0}^{2n}\Bigg(\sup_{[0,T]} \norm{w_{1}}_{L^2(\R^d)}^2+\int_0^T\norm{\nabla w_{1}}_{L^2(\R^d)}^2\Bigg)\\
&\leqslant 2C^n \abs{\theta-\theta_0}^{2n} R_T^2.
\end{align*}
By choosing 
\begin{gather}\label{c3.36}
\abs{\theta-\theta_0}^2< \frac{1}{2C},
\end{gather}
we end up with the strong convergence of $(u_n)_n$ to a limit $u$ in $L^\infty((0,T),L^2(\R^d))\cap L^2((0,T), H^1(\R^d))$. 
It is  a simple matter to check from \eqref{c7.3} that $u$ is a weak solution of 
\begin{gather}\label{c7.4}
\begin{cases}
\rho_0 \dpt u- \dvg(2\mu^{\theta_0}\D u+\lambda^{\theta_0}\dvg u I_d)=\dvg F+(\theta-\theta_0)\dvg \big(2(\mu-\widetilde\mu)\D u+(\lambda-\widetilde\lambda)\dvg u I_d\big),\\
u_{|t=0}=u_0,
\end{cases}
\end{gather}
and consequently, 
of \eqref{eq3.19}. Then, by virtue of lower semi-continuity of norms, 
we deduce form the uniform bound \eqref{c7.2}:
\[
\norm{w}_{X_T}\leqslant R_T.
\]
Thus $u$ has all the regularity required in \eqref{ep1.3}, and it only remains to establish time continuity:
\[
u\in \cC([0,T], H^1(\R^d)),\, \text{ and }\, \dpt u\in \cC([0,T],  L^2(\R^d)).
\]
The first one follows from the fact that $\dpt u\in L^2((0,T), H^1(\R^d))$. Moreover, since $\sqrt\sigma u_{tt}\in L^2((0,T)\times \R^d)$, it transpires that $u_t\in \cC((0,T], L^2(\R^d))$. Therefore, it only remains to prove the continuity of $ u_t$ at the initial time. To achieve this, we observe that the sequence  $(u_n)_n$ satisfies the following estimate (see \eqref{eq4.1}):
\begin{align*}
   \dfrac{1}{2}\int_{\R^d}\rho_{0}\abs{ \dpt u_n(\tau)}^{2}&\leqslant \dfrac{1}{2}\int_{\R^d}\rho_{0} \abs{{\dpt u_n}_{|t=0}}^{2} + C\Bigg(\int_0^\tau \norm{ F_t}_{L^{2}(\R^d)}^2 +  \int_0^\tau\norm{\mu_t^\theta,\,\lambda_t^\theta}_{L^4(\R^d)}^2\norm{\nabla u_n}_{L^4(\R^d)}^2\Bigg)\\
   &\leqslant \dfrac{1}{2}\int_{\R^d}\rho_{0} \abs{{\dpt u_n}_{|t=0}}^{2}+ \tau CR_T^2+C\int_0^\tau \norm{ F_t}_{L^{2}(\R^d)}^2.
\end{align*}
For all $n\geqslant 1$, since ${u_{n+1}}_{|t=0}={u_n}_{|t=0}=u_0$, it follows from  \eqref{c7.4}-\eqref{c7.5} that:
\[
\rho_0 {\dpt u_{n+1}}_{|t=0}=\dvg \Big(2\mu_0^\theta\D u_0+\lambda_0^\theta\dvg u_0+ F_{|t=0}\Big)={\rho_0u_t}_{|t=0}.
\]
Next, using the compact embedding $H^1(\R^d)$ into $L^2_\loc (\R^d)$ and Aubins-Lions theorem, we obtain $w^{n}_t(\tau)\to u_t(\tau)$ a.e in $\R^d$. Therefore, Fatou's Lemma  yields:
\[
 \dfrac{1}{2}\int_{\R^d}\rho_{0}\abs{ u_t(\tau)}^{2}\leqslant \dfrac{1}{2}\int_{\R^d}\rho_{0} \abs{{u_t}_{|t=0}}^{2}+ \tau CR_T^2+C\int_0^\tau \norm{ F_t}_{L^{2}(\R^d)}^2
\]
from which we deduce 
\[
\limsup_{\tau\to 0}\norm{\sqrt{\rho_{0}} u_t(\tau)}_{L^2(\R^d)}\leqslant \norm{\sqrt{\rho_{0}} {u_t}_{|t=0}}_{L^2(\R^d)}.
\]
This yields the continuity of $u_t$ at the initial time, and thereby completing the proof that $u\in X_{T,u_0}$.

Given that the constant $C$ appearing in \eqref{c3.35} and \eqref{c3.36} does not depend on $u_0$ nor $F$, it makes the condition \eqref{eq3.53} true, for all $(F,u_0)\in Y_T^{\theta}$
and one deduces that $\mathscr E$ is an open set of $[0,1]$. 

We turn to the proof of the fact that $\mathscr E$ is a closed set of $[0,1]$. 
Let $(\theta_n)_n$ a sequence of element of $\mathscr E$ that converges to some $\theta_0\in [0,1]$. Our aim is to prove that $\theta_0\in \mathscr E$. In other words, the following Cauchy problem 
\begin{gather}
\begin{cases}
\rho_0 u_t- \dvg(2\mu^{\theta_n}\D u+\lambda^{\theta_n}\dvg u I_d)=\dvg F+(\theta_0-\theta_n)\dvg \{2(\mu-\widetilde\mu)\D u+(\lambda-\widetilde\lambda)\dvg u I_d\},\\
u_{|t=0}=u_0
\end{cases}
\end{gather}
admits a unique solution in $X_T$ for all $(F,u_0)\in Y_T^{\theta_0}$. 
By following the computations  in Step 1-Step 4, and given that the norms of $\mu^\theta$ and $\lambda^\theta$ are uniformly bounded with respect to $\theta$, we can choose a sufficiently large $n$ independently of $F$, $u_0$. Hereby, $\mathscr E$ is a closed set of $[0,1]$.

In conclusion,  $\mathscr E$ is at same time a non-empty, open and closed set of $[0,1]$, since $[0,1]$ is a connected set, $\mathscr E=[0,1]$, $1\in \mathscr E$ and \cref{Theorem 3.1} finally holds.
\enddem
\section{\texorpdfstring{Study of the full nonlinear system \protect{\eqref{c4.39}}}{}}\label[section]{nonlinearised}
\subsection{Preliminary}
This section is devoted to the solvability of the full nonlinear two-fluid model in Lagrangian coordinates \eqref{c4.39}. The initial data $(c_0, \rho_0,u_0)$ is such that the compatibility condition holds true, that is:
\begin{gather}\label{c4.34}
\dvg(2\mu(\rho_0,c_0)\D u_0+(\lambda(\rho_0,c_0)\dvg u_0- P(\rho_0,c_0)) I_d)\in L^2(\R^d).
\end{gather}
Our method relies on a fixed-point argument, thus we will begin by setting some tools. We recall that such solution $u$ verifies (for simplicity, we write $u$ instead of $\overline{u}$):
\begin{gather}\label{c4.35}
    \rho_0 u_t-\dvg\{2\mu(J^{-1}_{ u}\rho_0,c_0)\D  u+\lambda(J^{-1}_{ u} \rho_0,c_0))\dvg u I_d\}=\dvg I (u,u)
\end{gather}
where 
\[
I(u,u)= I_1(u,u)+I_2(u,u)+ I_3(u)
\]
and 
\begin{gather}\label{c7.9}
\begin{cases}
I_1(v,w)&= 2\mu(J^{-1}_{v}\rho_0,c_0)\{\adj(D\X_{v})\D_{A_v} w-\D w\},\\
I_2(v,w)&=\lambda(J^{-1}_{v}\rho_0,c_0)\{\adj(D\X_{v})\dvg_{A_v} w-\dvg  w\},\\
I_3(v)&=  \adj(D\X_{v})\{P(J^{-1}_{v}\rho_0,c_0)-\widetilde P\}I_d.
\end{cases}
\end{gather}
According to the previous section, we will construct $u$ in a closed subset of the vector space $X_T$ as defined
in  \eqref{ep1.3}
for some small time $T$. For some $v\in X_T$ with $v_{|t=0}= u_0$, we consider the linear Cauchy problem:
\begin{gather}\label{eq3.1}
    \begin{cases}
        \rho_0 u_t-\dvg\{2\mu(J^{-1}_{ v}\rho_0,c_0)\D  u+\lambda(J^{-1}_{ v} \rho_0,c_0)\dvg \ u Id\}=\dvg I (v,v),\\
        u_{|t=0}=u_0.
    \end{cases}
\end{gather}
The validity of the last condition in the definition of $Y_T:=Y_T(\mu(\rho_0,c_0),\lambda(\rho_0,c_0))$ (see \eqref{c7.6}) results from the fact that $v_{|t=0}=u_0$, ${\big(J_v^{-1}\big)}_{|t=0}=1$, $I(v,v)_{|t=0}=\big(P(\rho_0,c_0)-\widetilde P\big) I_d$ and \eqref{c4.34}. In the following section, we set up tools for the construction of the solution to \eqref{c4.35}. 
 
\subsection{Implementation of the fixed-point theorem}
In this section, we will justify the well-posedness  of a contracting map whose fixed-point is the solution of the nonlinear system \eqref{c4.35}. To this end, we fix some $v$ in $X_{T,u_0}=\{ u\in X_T\colon u_{|t=0}=u_0\}$.  The first part of the following section will consist in proving that the regularity of $v\in X_T$ is sufficient to ensure that $I(v,v)$ belongs to $Y_T$. From there, we will be able to deduce the existence of a solution for the Cauchy problem \eqref{eq3.1} using \cref{Theorem 3.1}. This leads to the following result.
\begin{prop}\label[prop]{prop41}
    There are $T^*>0$ and $R>0$ such that for all $T\in (0,T_*)$ the  map: 
\begin{align}
    \Phi_T\colon \; X_{T,u_0}(R)&\longrightarrow  X_{T,u_0}(R)\nonumber\\
                       v&\longmapsto  \Phi_T(v):= u,\label{c4.44}
\end{align}
where $u$ is the unique solution of the Cauchy problem \eqref{eq3.1}, and $X_{T,u_0}(R):=\{ u\in X_{T,u_0}\colon \norm{u}_{X_T}\leqslant R\}$, is well-defined.
\end{prop}

The final step of this section is dedicated to proving that for some small $T\in (0,T_*)$, the map $\Phi_T$ admits a fixed point, see \cref{prop42} below.
 
\dem[Proof of \cref{prop41}]
To ease readability, for $T>0$, and a normed space $\mathcal{X}$, the space $L^p((0,T), \mathcal{X})$ will be denoted $L^p_T\mathcal{X}$. Also, we will only present the main steps, technical computations are carried out in \cref{esti1} below. 

First of all, let us point out that since $\sigma^{r/4} \nabla v\in L^4_T\cC^\alpha_{pw,\mathcal{C}_0}(\R^d)$, with  $0<r<8/3$, one has:
\[
c=\int_0^T\norm{\nabla v}_{\cC^\alpha_{pw,\mathcal{C}_0}(\R^d)}\leqslant  T^{\tfrac{3-r}{4}}\norm{v}_{X_T}.
\]
Suppose that  $\norm{v}_{X_T}\leqslant R$ for some $R>0$ to be fixed later. One chooses $T\leqslant 1$ small in such a way that $2c<1$. Thus, the 
constants involved in \cref{lem3} and \cref{lem4} below are assumed to be fixed.  Furthermore,  the constant in the notation $\lesssim$  in all subsequent steps, may depend on $c<1/2$, the bounds of $\rho_0,\mu(\rho_0,c_0),\lambda(\rho_0,c_0)$, the piecewise H\"older norm of $\mu(\rho_0,c_0)$, $\lambda(\rho_0,c_0)$ and the interface regularity.

The following bounds for $I(v,v)$ will be proved in  \cref{esti1-step1} (we freely use $T\leqslant 1$):
\begin{gather}\label{epq45}
\sup_{[0,T]}\norm{I(v,v)}_{ L^2(\R^d)} +\norm{I(v,v)}_{ L^\infty_TL^6(\R^d))} \lesssim \norm{P(\rho_0,c_0)-\widetilde P}_{L^2\cap L^6(\R^d)}+T^{\tfrac{3-r}{4}} R(1+R),
\end{gather}
\begin{gather}\label{epq49}
\norm{I(v,v)}_{ L^{16}_TL^8(\R^d)}\lesssim  T^{1/16}\norm{P(\rho_0,c_0)-\widetilde P}_{L^8(\R^d)}+T^{\tfrac{3-r}{4}} R(1+R),
\end{gather}
\begin{gather}\label{epq53}
\norm{\sigma^{r/4}I(v,v)}_{L^4_T\cC^\alpha_{pw,\mathcal{C}_0}(\R^d))}\lesssim T^{1/4} \norm{P(\rho_0,c_0)-\widetilde P}_{\cC^\alpha_{pw,\mathcal{C}_0}(\R^d)}+T^{\tfrac{3-r}{4}}R(1+ R).
\end{gather}
On the other hand, the following estimates for $\dpt I(v,v)$ will be proved in \cref{esti1-step2} (we freely use $T^{\tfrac{3-r}{4}} R\leqslant 1$): 
\begin{gather}\label{epq47}
\norm{\dpt I(v,v)}_{L^2_TL^2(\R^d)}\lesssim  T^{1/2}R(R+1)+ T^{\tfrac{3-r}{4}}R^2,
\end{gather}
\begin{gather}\label{epq50}
\norm{\sigma^{\tfrac{d}{4}} \dpt I(v,v)}_{L^{8/d}_TL^4(\R^d)}\lesssim   T^{\tfrac{d}{8}}R(R+1)+T^{\tfrac{3-r}{4}}R^2,
\end{gather}
and  finally, we will derive in \cref{esti1-step3} the following estimates for $\partial_{tt} I(v,v)$:
\begin{gather}\label{epq48}
\norm{\sigma\partial_{tt} I(v,v)}_{L^2_TL^2(\R^d)}\lesssim T^{1/2}R( R+ 1)+ T^{\tfrac{1}{2}-\tfrac{d}{8}}R^2(1+ T^{\tfrac{d}{8}} R)+ T^{\tfrac{3-r}{4}}R^2.
\end{gather}
We obtain that $I(v,v)$ satisfies all conditions required in the definition of  $Y_T$. Now, we shall ensure that the viscosity  $\mu(J_{v}^{-1}\rho_0,c_0)$ and $\lambda(J_{v}^{-1}\rho_0,c_0)$ have the regularities  required in \eqref{ep1.2}, and the dynamic viscosity $\mu(J_{v}^{-1}\rho_0,c_0)$ remains  
close to $\widetilde\mu$. We  prove in \cref{esti1-step4}:
\begin{gather}\label{c4.40}
\begin{cases}
    \norm{\mu(J^{-1}_v\rho_0,c_0)-\widetilde\mu }_{L^\infty_T\cC^\alpha_{pw,\mathcal{C}_0}(\R^d)}&\leqslant \norm{\mu(\rho_0,c_0)-\widetilde\mu}_{ \cC^\alpha_{pw,\mathcal{C}_0}(\R^d)}+ \norm{\partial_\rho\mu}_{L^\infty}T^{\tfrac{3-r}{4}}R,\\
    \norm{\llbracket\mu(J^{-1}_v\rho_0,c_0)\rrbracket}_{L^\infty_TL^\infty(\mathcal{C}_0)}&\leqslant \norm{\llbracket\mu(\rho_0,c_0)\rrbracket}_{L^\infty(\mathcal{C}_0)}+ \norm{\partial_\rho\mu}_{L^\infty}T^{\tfrac{3-r}{4}}R,
\end{cases}
\end{gather}
with similar bounds for $\lambda(J^{-1}_v\rho_0,c_0)$. Assuming \eqref{c3.37}-\eqref{ep3.2}, we  choose $T$ small such that 
\begin{align}
     \sup_{[0,T]}\left[1\right.&\left.+\norm{\lambda (J_{v}^{-1}\rho_0,c_0)}_{\dot \cC^\alpha_{pw,\mathcal{C}_0}(\R^d)}+\left( \mathfrak{P}_{\mathcal{C}_0} +\ell^{-\alpha}_{\vph_0}\right)\llbracket\lambda(J_{v}^{-1}\rho_0,c_0)\rrbracket_{L^\infty(\mathcal{C}_0)}\right]\sup_{[0,T]}\norm{\mu(J_{v}^{-1}\rho_0,c_0)-\widetilde\mu}_{\cC^\alpha_{pw,\mathcal{C}_0}(\R^d)}\nonumber\\
    &+\left( \mathfrak{P}_{\mathcal{C}_0} +\ell^{-\alpha}_{\vph_0}\right)\sup_{[0,T]}\norm{\llbracket\mu(J_v^{-1}\rho_0,c_0)\rrbracket}_{L^\infty(\mathcal{C}_0)}\nonumber\\
    &+\sup_{[0,T]}\left[\norm{\llbracket\mu(J_v^{-1}\rho_0,c_0)\rrbracket,\, \llbracket \lambda(J_v^{-1}\rho_0,c_0)\rrbracket}_{L^\infty(\mathcal{C}_0)}\left\|1-\dfrac{\widetilde\mu}{\scalar{\mu(J_v^{-1}\rho_0,c_0)}}\right\|_{L^\infty(\mathcal{C}_0)}\right]\leqslant \dfrac{[\mu]}{2},\nonumber
    \end{align}
and 
\[
\mu_{*,0}/4\leqslant \mu(J_{v}^{-1}\rho_0,c_0)\quad \text{and}\quad \lambda(J_{v}^{-1}\rho_0,c_0),\; \mu(J_{v}^{-1}\rho_0,c_0) \leqslant 4\max\left\{\sup_{\R^d}\mu(\rho_0,c_0),\; \sup_{\R^d} \lambda(\rho_0,c_0)\right\}.
\]

We  prove in \cref{esti1-step5} below, the following estimates for the time derivative of the viscosity:
\begin{gather}\label{c4.41}
\norm{\dpt \mu(J^{-1}_v\rho_0,c_0),\,\dpt \lambda(J^{-1}_v\rho_0,c_0)}_{L^\infty_TL^4(\R^d)}\lesssim R,\quad\quad\; \norm{\dpt \mu(J^{-1}_v\rho_0,c_0),\, \dpt \lambda(J^{-1}_v\rho_0,c_0)}_{L^{16}_TL^8(\R^d)}\lesssim R,
\end{gather}
and 
\begin{gather}\label{c4.42}
\norm{\dpt \mu(J^{-1}_v\rho_0,c_0),\,\dpt \lambda(J^{-1}_v\rho_0,c_0)}_{L^1_TL^\infty(\R^d)}+
\norm{\sqrt{\sigma}\dpt \mu(J^{-1}_v\rho_0,c_0),\, \sqrt{\sigma}\dpt \lambda(J^{-1}_v\rho_0,c_0)}_{L^2_TL^\infty(\R^d)}\lesssim T^{\tfrac{3-r}{4}}R.
\end{gather}
Finally, the $L^2_TL^4(\R^d)$-norm of the second time derivative of the viscosity is obtained in \cref{esti1-step6}, from which we deduce:
\begin{gather}\label{c4.43}
\norm{\sigma\partial_{tt}\mu(J^{-1}_v\rho_0,c_0),\, \sigma\partial_{tt}\lambda(J^{-1}_v\rho_0,c_0)}_{L^{2}_TL^4(\R^d)}\lesssim 
T^{\tfrac{1}{2}-\tfrac{d}{8}}R(1+R).
\end{gather}
Hence, the viscosities $\mu$ and $\lambda$ verify all the conditions stated in \eqref{ep1.2} and in particular \cref{Theorem 3.1} provides us with a unique solution 
$u\in X_{T,u_0}$ of the Cauchy problem  \eqref{eq3.1}. 
\paragraph{\textbf{Estimates for $u$ and $u_t$}} 
  \begin{align}
     \sup_{[0,T]}\norm{u,\,u_t,\, \nabla u}_{L^2(\R^d)}^{2}&+\norm{\nabla u}_{L^\infty_TL^6(\R^d)}^2+\norm{\nabla u}_{L^{16}_TL^8(\R^d)}^2+\int_{0}^{T}
\norm{\nabla u,\, \nabla u_t,\, u_t}_{L^2(\R^d)}^{2}\nonumber\\
&\lesssim \left(\norm{(I(v,v),u_0)}_{Y_T}^2+TR^2\norm{I(v,v)}_{L^\infty_T (L^2\cap L^6)(\R^d)}^2\right) \exp\left(C_*T^{\tfrac{3-r}{4}}R+C_*TR^2\right)\label{epq46}
\end{align}
with the estimates of $I(v,v)$ given in \eqref{epq45}-\eqref{epq48}. 
\paragraph{\textbf{Estimates for $u_{tt}$}}
\begin{multline}\label{epq51}
    \int_0^T\norm{ \sqrt\sigma u_{tt},\; \sigma \nabla u_{tt}}_{L^2(\R^d)}^2+\norm{\sqrt\sigma \nabla u_t,\; \sigma u_{tt}}^2_{L^\infty_TL^2(\R^d)}\lesssim \left[ \int_0^{\sigma(T)}\norm{\nabla u_t}^2_{L^2(\R^d)}\right.\\
    \left.+\int_0^{T} \norm{\dpt I(v,v),\, \sigma \partial_{tt} I(v,v)}_{L^2(\R^d)}^2+  (T R^2+T^{1-\tfrac{d}{4}}R^2(1+R)^2)\norm{\nabla u}_{L^\infty_TL^4(\R^d)}^2\right]\times\exp\left( C_* T^{\tfrac{3-r}{4}}R\right).
\end{multline}
The norms of $u$ involved in this estimate are controlled in \eqref{epq46}, while  those of $I(v,v)$ are given in \eqref{epq45}-\eqref{epq48}.
\paragraph{\textbf{$L^{8/d}_TL^4(\R^d)$-norm estimate  for $\sqrt{\sigma}\nabla u_t$} }
\begin{multline}\label{epq54}
\int_0^T\sigma^{2}\norm{\nabla u_t}_{L^4(\R^d)}^{8/d}\lesssim \int_0^T\sigma^{2}\norm{\dpt I(v,v)}_{L^4(\R^d)}^{8/d}+T^{1-\tfrac{1}{d}}R^{8/d}\norm{\nabla u}_{L^{16}_TL^8(\R^d)}^{8/d}+ \norm{\sigma u_{tt}}_{L^\infty_TL^2(\R^d)}^{8/d}\\
        +\left[\int_0^T\norm{\nabla u_t,\, \dpt I(v,v)}_{L^2(\R^d)}^{2}\right]^{4/d}+T^{4/d}R^{8/d}\norm{\nabla u}_{L^\infty_TL^4(\R^d)}^{8/d}.
\end{multline}
The norms of $I(v,v)$ are given in \eqref{epq50} and \eqref{epq47}, while the norms of $u$ are controlled in  \eqref{epq46} and \eqref{epq51}.
\paragraph{ \textbf{$L^4_T\cC^\alpha_{pw,\mathcal{C}_0}(\R^d)$ estimate for $\sigma^{r/4}\nabla u$}} 
\begin{align*}
\int_0^T\sigma^{r}\norm{\nabla u}_{\cC^\alpha_{pw,\mathcal{C}_0}(\R^d)}^4 
         &\lesssim  \int_0^T\sigma^{r}\norm{I(v,v)}_{\cC^\alpha_{pw,\mathcal{C}_0}(\R^d)}^4 +\norm{I(v,v),\,  \nabla u,\, u_t,\,\sqrt\sigma\nabla u_t}_{L^\infty_TL^2(\R^d)}^4\nonumber\\
         &+\left[\int_0^T \norm{I(v,v),\, \nabla u,\,u_t,\,  \nabla u_t}_{L^2(\R^d)}^2\right]^2+\left[\int_0^T \sigma^{2}\norm{\nabla u_t}_{L^4(\R^d)}^{8/d}\right]^{d/2}.
\end{align*}
The norms of $I(v,v)$ are given in \eqref{epq45} and \eqref{epq53} whereas the norms of $u$ involved are controlled in  \eqref{epq46}, \eqref{epq51} and \eqref{epq54}.

By combining all these computations, we choose $R$ related to the norm of the initial data and a small $T^*$ such that for all $0<T\leqslant T_*$, 
\[
\norm{u}_{X_T}\leqslant R.
\]
This concludes the proof of \cref{prop41}. 
\enddem

The final step of this section is dedicated to constructing the fixed point $u\in X_{T,u_0}$ of $\Phi_T$. 
\begin{prop}\label[prop]{prop42}
    There exists $T\in (0, T_*)$ and $u\in X_{T,u_0}$ such that 
    \[
u=\Phi_T(u).
\]
\end{prop}

This fixed point $u$ satisfies the nonlinear equations \eqref{c4.35}, thereby completing the proof of \cref{remaL}.

\dem[Proof of \cref{prop42}] 
To prove \cref{prop42}, we begin by considering the following sequence:
\begin{gather*}
    u_{n+1}= \Phi_T(u_n), \; n\geqslant 0,
\end{gather*}
with 
\begin{gather}\label{c7.7}
    \begin{cases}
        \rho_0 \dpt u_{n+1}-\dvg\{2\mu(J^{-1}_{u_n}\rho_0,c_0)\D  u_{n+1}+\lambda(J^{-1}_{ u_n} \rho_0,c_0)\dvg  u_{n+1} I_d\}=\dvg \big(I (u_n,u_n)\big),\\
        {u_{n+1}}_{|t=0}=u_0.
    \end{cases}
\end{gather}
From \cref{prop41}, it follows that the sequence $(u_n)_n$  is uniformly bounded in $X_T$:
\begin{gather}\label{c7.11}
\sup_{n\in \N} \norm{ u_n}_{X_T}\leqslant R.
\end{gather}
We take the differences in \eqref{c7.7} and  obtain: 
\begin{gather}\label{c7.8}
\begin{cases}
    \rho_0\dpt w_{n+1}&-\dvg\Big(2\mu(J_{u_n}^{-1}\rho_0,c_0)\D w_{n+1}+\lambda(J_{u_n}^{-1}\rho_0,c_0)\dvg w_{n+1} I_d\Big)\\
    &= \dvg \big(\delta \mathcal I\big)+ \dvg \big(I(u_n,u_n)-I(u_{n-1},u_{n-1})\big),\\
    {w_{n+1}}_{|t=0}&=0,
\end{cases}
\end{gather}
where $w_{k+1}=u_{k+1}-u_k$, $k\geqslant 1$, and 
\begin{gather*}
    \delta \mathcal I=2\big((\mu(J^{-1}_{u_n}\rho_0,c_0)-\mu(J^{-1}_{u_{n-1}}\rho_0,c_0)\big)\D u_{n-1}+
    \big(\lambda(J^{-1}_{u_n}\rho_0,c_0)-\lambda(J^{-1}_{u_{n-1}}\rho_0,c_0)\big)\dvg u_{n-1} I_d.
\end{gather*}
An energy estimate in \eqref{c7.8} yields:
\begin{gather}\label{c7.10}
\sup_{[0,T]}\norm{ w_{n+1}}_{L^2(\R^d)}^2+\int_0^T\norm{\nabla w_{n+1}}_{L^2(\R^d)}^2\lesssim \int_0^T\norm*{\delta \mathcal{I},\,I(u_n,u_n)-I(u_{n-1},u_{n-1})}_{L^2(\R^d)}^2.
\end{gather}
To estimate the right-hand side above in terms of the norms of $u_n$ and $u_{n-1}$, we rely on the following result from \cite[Lemma A.4]{danchinLagrangianapproach}.
\begin{lemm}\label[lemma]{lem5}
Let $v_1,\, v_2$ be two vectors fields  verifying  
\[
c=\max_j\int_0^T\norm{D v_j}_{L^\infty(\R^d)}<1.
\]
There exists a constant $k=k(c)$ such that the following estimates hold true: for all $p\in [2,\infty]$
\begin{gather*}
\begin{cases}
    \norm{A_{v_2}(t)-A_{v_1}(t)}_{L^p(\R^d)}&\leqslant k \norm{D (v_2-v_1)}_{L^1_tL^p(\R^d)},\\
    \norm{\adj(D\X_{v_2}(t))-\adj(D\X_{v_1}(t))}_{L^p(\R^d)}&\leqslant k \norm{D (v_2-v_1)}_{L^1_tL^p(\R^d)},\\
    \norm{J^{\pm 1}_{v_2}(t)-J^{\pm 1}_{v_1}(t)}_{L^p(\R^d)}&\leqslant k \norm{D (v_2-v_1)}_{L^1_tL^p(\R^d)}.
\end{cases}     
    \end{gather*}
\end{lemm}
\cref{lem5} immediately implies:
\begin{align*}
\int_0^T\norm{\delta\mathcal{I}(t)}_{L^2(\R^d)}^2dt&\lesssim \int_0^T \norm{J_{u_n}^{-1}(t)-J_{u_{n-1}}^{-1}(t)}^2_{L^2(\R^d)}\norm{D u_{n-1}(t)}_{L^\infty(\R^d)}^2dt\\
&\lesssim \bigg(\int_0^T t\norm{D u_{n-1}(t)}_{L^\infty(\R^d)}^2dt\bigg)\bigg(\int_0^T\norm{ D w_n}_{L^2(\R^d)}^2\bigg)\\
&\lesssim T^{\frac{3-r}{2}} R^2 \int_0^T\norm{ D w_n}_{L^2(\R^d)}^2.
\end{align*}
Similarly, by expressing  (see the definitions of $I_1$ and $I_2$ in \eqref{c7.9})
\begin{align*}
    I_3\big(u_n\big)-I_3(u_{n-1}\big)&=\{\adj (D \X_{u_n}) - \adj (D \X_{u_{n-1}})\}\left(P(J_{u_n}^{-1}\rho_0,c_0)-\widetilde P\right)\\
                 &+ \adj (D \X_{u_{n-1}})\{P(J_{u_n}^{-1}\rho_0,c_0)-P(J_{u_{n-1}}^{-1}\rho_0,c_0)\},
\end{align*}
and
\begin{align*}
    I_1\big(u_n,u_n\big)-I_1\big(u_{n-1},u_{n-1}\big)&= 2 \{ \mu(J^{-1}_{u_n}\rho_0,c_0)-\mu(J^{-1}_{u_{n-1}}\rho_0,c_0)\}\{\adj (D \X_{u_n})\D _{A_{u_n}} {u_n}-\D {u_n}\}\nonumber\\
                     &+2\mu(J^{-1}_{u_{n-1}}\rho_0,c_0)\{\adj(D \X_{u_n})\D _{A_{u_n}} w_n-\D w_n\}\nonumber\\
                     &+2\mu(J^{-1}_{u_{n-1}}\rho_0,c_0)\{\adj(D \X_{u_n})-\adj(D \X_{u_{n-1}})\}\D_{A_{u_n}} u_{n-1}\nonumber\\
                     &+2\mu(J^{-1}_{u_{n-1}}\rho_0,c_0)\adj(D \X_{u_{n-1}})\{\D_{A_{u_n}}-\D_{A_{u_{n-1}}}\} u_{n-1},
\end{align*}
we infer 
\[
\int_0^T\norm{I_3\big(u_n\big)-I_3\big(u_{n-1}\big)}_{L^{2}(\R^d)}^2 \lesssim T^2\int_0^T \norm{\nabla  w_n}_{L^{2}(\R^d)}^2,
\]
and 
\[
\int_0^T\norm{I_1\big(u_n,u_n\big)-I_1\big(u_{n-1},u_{n-1}\big)}_{L^{2}(\R^d)}^2\lesssim T^{\frac{3-r}{2}} R^2 \int_0^T\norm{ \nabla w_n}_{L^2(\R^d)}^2.
\]
A similar estimate also holds for  $I_2\big(u_n,u_n\big)-I_2\big(u_{n-1},u_{n-1}\big)$. It therefore follows from \eqref{c7.10} that:
\begin{gather*}
\sup_{[0,T]}\norm{w_{n+1}}_{L^2(\R^d)}^2+\int_0^T\norm{\nabla w_{n+1}}_{L^2(\R^d)}^2\leqslant C T^{\frac{3-r}{2}}\big(1+R^2\big)\Bigg(\sup_{[0,T]}\norm{ w_{n}}_{L^2(\R^d)}^2+\int_0^T\norm{\nabla w_{n}}_{L^2(\R^d)}^2\Bigg).
\end{gather*}
We then choose $T<T_*$ small such that
\[
C T^{\frac{3-r}{2}}\big(1+R^2\big)\leqslant \dfrac{1}{2},
\]
which leads to
\[
\sup_{[0,T]}\norm{w_{n+1}}_{L^2(\R^d)}^2+\int_0^T\norm{\nabla w_{n+1}}_{L^2(\R^d)}^2\leqslant \dfrac{1}{2}\Bigg(\sup_{[0,T]}\norm{ w_{n}}_{L^2(\R^d)}^2+\int_0^T\norm{\nabla w_{n}}_{L^2(\R^d)}^2\Bigg),
\]
and hence
\begin{align*}
\sup_{[0,T]}\norm{w_{n+1}}_{L^2(\R^d)}^2+\int_0^T\norm{\nabla w_{n+1}}_{L^2(\R^d)}^2&\leqslant \dfrac{1}{2^n}\Bigg(\sup_{[0,T]}\norm{w_1}_{L^2(\R^d)}^2+\int_0^T\norm{\nabla \big( w_{1}\big)}_{L^2(\R^d)}^2\Bigg)\\
&\leqslant \dfrac{R^2}{2^{n-1}}.
\end{align*}
This implies the strong convergence of the sequence $(u_n)_n$ in $L^\infty((0,T), L^2(\R^d))\cap L^2((0,T), H^1(\R^d))$ to some limit $u$. Furthermore, by virtue of lower semi-continuity of norms, we deduce from the uniform bound  \eqref{c7.11}:
\[
\norm{u}_{X_T}\leqslant R.
\]
This enables us to apply \cref{lem5}  with $v_1= u_n$ and $v_2=u$:
\begin{gather*}
\begin{cases}
    \norm{A_{u_n}-A_{u}}_{L^\infty_TL^2(\R^d)}&\leqslant k \norm{D (u_n-u)}_{L^1_TL^2(\R^d)},\\
    \norm{\adj(D\X_{u_n})-\adj(D\X_{u})}_{L^\infty_TL^2(\R^d)}&\leqslant k \norm{D (u_n-u)}_{L^1_TL^2(\R^d)},\\
    \norm{J^{\pm 1}_{u_n}-J^{\pm 1}_{u}}_{L^\infty_TL^2(\R^d)}&\leqslant k \norm{D (u_n-u)}_{L^1_TL^2(\R^d)},
\end{cases}     
    \end{gather*}
from which we deduce the strong convergence in  $L^\infty_T L^2(\R^d)$ of the sequences $\big( I_d-A_{u_n}\big)_n$, $\big( I_d-\adj(D\X_{u_n})\big)_n$ and $\big(J^{\pm 1}_{u_n}-1\big)_n$   to $I_d-A_{u}$, $I_d-\adj(D\X_{u})$ , and $J^{\pm 1}_{u}-1$, respectively. In fact, these strong convergences hold in $L^\infty((0,T), L^p(\R^d))$ for all $p\in [2, \infty)$, given that the sequences are uniformly bounded in $L^\infty((0,T)\times \R^d)$. Additionally, the bound in \eqref{c7.11} indicates that the sequence  $(\nabla u_n)_n$ is bounded in $L^{16}((0,T), L^8(\R^d))$, which implies that $(\nabla u_n)_n$ converges strongly to $\nabla u$ in  $L^p((0,T)\times \R^d)$ for all $p\in [2,8)$. We are then able to pass to the limit in \eqref{c7.7} and conclude that the limit $u$  is a weak solution of \eqref{c4.35}. The remaining step is to establish the following time continuity:
\[
u\in \cC([0,T], H^1(\R^d)),\,\text{ and }\, u_t\in \cC([0,T], L^2(\R^d)).
\]
This property can be proved by following the arguments in Step 4 of \cref{proof2} and using \eqref{epq47}. This completes the proof of \cref{prop42}.
\enddem


\subsection{Estimates}\label[section]{esti1}
In this section,, we derive the estimates for $I(v,v)$, $\mu(\rho_0J^{-1}_v,c_0)$ and $\lambda(\rho_0J^{-1}_v,c_0)$ used in the proof of \cref{prop41} above. We first recall some estimates for the flow of a Lipschitz vector field. We refer to \cite[Lemma A.3]{danchinLagrangianapproach} for the proof.
\begin{lemm}\label[lemma]{lem3}
    Let $v$ be a vector field with $\nabla v\in L^1((0,T), L^\infty(\R^d))$. Suppose that 
    \[
    c=\int_0^T\norm{\nabla v}_{L^\infty(\R^d)}<1.
    \]
    There exists a constant $k=k(c)$ such that the following estimates hold true:
   \[ 
   \begin{cases}
   \norm{I_d-\adj(D \X_v(t));\; I_d-A_v(t)}_{L^\infty(\R^d)}&\leqslant k \norm{\nabla v}_{L^1((0,t), L^\infty(\R^d))},\\
   \norm{J_v^{\pm 1}(t)-1}_{L^\infty(\R^d)}&\leqslant k \norm{\nabla v}_{L^1((0,t), L^\infty(\R^d))}.
   \end{cases}
   \]
   Moreover, one has the following:
   \[
   \begin{cases}
   \norm{\adj( D\X_v(t)) D_{A_v(t)} w-D w }_{L^p(\R^d)} &\leqslant k  \norm{\nabla v}_{L^1((0,t), L^\infty(\R^d))}\norm{D w}_{L^p(\R^d)},\\
   \norm{\adj(D \X_v(t))\dvg_{A_v(t)} w-\dvg w I_d}_{L^p(\R^d)}&\leqslant k  \norm{\nabla v}_{L^1((0,t), L^\infty(\R^d))}\norm{D w}_{L^p(\R^d)}.
   \end{cases}
   \]
\end{lemm}
In the spirit of the above result, and owing to the fact that $\cC^\alpha_{pw,\mathcal{C}_0}(\R^d)$ is an algebra, 
one has the following estimates.
\newline
\begin{lemm}\label[lemma]{lem4}
    Let $v$ be a vector field fulfilling  $\nabla v\in L^1((0,T), \cC^\alpha_{pw,\mathcal{C}_0}(\R^d))$. Suppose that 
    \begin{gather}\label{eq3.65}
    c=\int_0^T\norm{\nabla v}_{\cC^\alpha_{pw,\mathcal{C}_0}(\R^d)}<1.
    \end{gather}
    There exists a constant $k=k(c)$ such that the following estimates holds true:
  \[
  \begin{cases}
  \norm{I_d-\adj(D \X_v(t));\; I_d-A_v(t)}_{\cC^\alpha_{pw,\mathcal{C}_0}(\R^d)}&\leqslant k \norm{\nabla v}_{L^1((0,t), \cC^\alpha_{pw,\mathcal{C}_0}(\R^d))},\\
   \norm{J_v^{\pm 1}(t)-1}_{\cC^\alpha_{pw,\mathcal{C}_0}(\R^d)}&\leqslant k \norm{\nabla v}_{L^1((0,t), \cC^\alpha_{pw,\mathcal{C}_0}(\R^d))}.
   \end{cases}
   \]
   Moreover, one has the following:
   \[
   \begin{cases}
   \norm{\adj( D\X_v(t)) \D_{A_v(t)} w-\D w }_{\cC^\alpha_{pw,\mathcal{C}_0}(\R^d)} &\leqslant k  \norm{\nabla v}_{L^1((0,t), \cC^\alpha_{pw,\mathcal{C}_0}(\R^d))}\norm{D w}_{\cC^\alpha_{pw,\mathcal{C}_0}(\R^d)},\\
   \norm{\adj(D \X_v(t))\dvg_{A_v(t)} w-\dvg w I_d}_{\cC^\alpha_{pw,\mathcal{C}_0}(\R^d)}&\leqslant k  \norm{\nabla v}_{L^1((0,t), \cC^\alpha_{pw,\mathcal{C}_0}(\R^d))}\norm{D w}_{\cC^\alpha_{pw,\mathcal{C}_0}(\R^d)}.
   \end{cases}
   \]
\end{lemm}
Notice that we will need to estimate the first and the second time derivative of $I(v,v)$ as it is required in the definition of $Y_T$.  From \eqref{eq3.4} we have $\adj(D \X_v ) \cdot D \X_v= \det (D\X_v) I_d$ and  
$\dpt D \X_v = D v$  and consequently:
\begin{lemm}\label[lemma]{lemexp}
The following expressions hold true:
     \begin{gather}\label{eq3.2}
\dpt J^{-1}_v(t,y)=-J_v^{-1}(t,y)\dvg_{A_v} (v),
\end{gather}
\begin{gather}\label{eq3.6}
\dpt A_v=-A_v \cdot D v \cdot A_v,
\end{gather}
\begin{gather}\label{eq3..3}
\dpt \adj(D \X_v )=(D v\colon \adj (D\X_v)) A_v -\adj(D \X_v ) \cdot D v\cdot A_v.
\end{gather}
Next, taking the time derivative of the above quantities, one has:
\begin{gather}\label{eq3.58}
J_v \partial_{tt} J^{-1}_v=(\dvg_{A_v} (v))^2- D v_t\colon A_v+D v\colon (A_v \cdot D v \cdot A_v ),
\end{gather}
\begin{gather}
    \partial_{tt} A_v=2A_v \cdot (D v \cdot A_v)^2- A_v \cdot D v_t \cdot A_v,
\end{gather}
\begin{align}
\partial_{tt} \adj(D \X_v )&=(D v_t\colon \adj (D\X_v) )A_v+J_v(D v\colon A_v)^2A_v-(D v\colon (\adj (D\X_v)\cdot D v\cdot A_v) )A_v\nonumber\\
&-2(D v\colon \adj (D\X_v) )A_v\cdot D v\cdot A_v+2J_vA_v \cdot (D v\cdot A_v)^2-\adj(D \X_v ) \cdot D v_t\cdot A_v.\label{eq3.59}
\end{align}
\end{lemm}
We now proceed to the estimates for $I(v,v)$, $\mu(\rho_0J^{-1}_v,c_0)$ and $\lambda(\rho_0J^{-1}_v,c_0)$.
\subsubsection{Lebesgue and H\"older norm estimates for $I(v,v)$}\label[section]{esti1-step1}
We first recall 
\begin{gather}\label{c4.46}
I(v,v)= I_1(v,v)+I_2(v,v)+ I_3(v)
\end{gather}
where
\begin{gather}\label{c4.45}
\begin{cases}
I_1(v,v)&= 2\mu(J^{-1}_{v}\rho_0,c_0)\{\adj(D\X_{v})\D_{A_v} v-\D v\},\\
I_2(v,v)&=\lambda(J^{-1}_{v}\rho_0,c_0)\{\adj(D\X_{v})\dvg_{A_v} v-\dvg  v\},\\
I_3(v)&=  \adj(D\X_{v})\{P(J^{-1}_{v}\rho_0,c_0)-\widetilde P\}I_d.
\end{cases}
\end{gather}
From \cref{lem3}, we have the following  estimates for the first two terms:
\[
\norm{I_1(v,v);I_2(v,v)}_{L^m_TL^{q}(\R^d)}\lesssim 
\norm{\nabla  v}_{L^1_TL^\infty(\R^d)}\norm{\nabla v}_{L^m_T L^{q}(\R^d)}
\]
for all $2\leqslant m,q\leqslant \infty$. Regarding the H\"older estimates, one has from \cref{lem4} the following:
\begin{gather*}
     \norm{\sigma^{r/4}(I_1(v,v);I_2(v,v))}_{L^4_T\cC^\alpha_{pw,\mathcal{C}_0}(\R^d)}
    \lesssim \norm{\nabla  v}_{L^1_T\cC^\alpha_{pw,\mathcal{C}_0}(\R^d)}\norm{\sigma^{r/4}\nabla  v}_{L^4_T\cC^\alpha_{pw,\mathcal{C}_0}(\R^d)}.
\end{gather*}
It now remains the norm of $I_3(v)$. We may first write:
\[
I_3(v)=\adj(D\X_{v})\{P(J^{-1}_{v}\rho_0,c_0) -P(\rho_0,c_0)+P(\rho_0,c_0)-\widetilde P\}
\]
and thus, we deduce from \cref{lem3}, the following:
\[
\norm{I_3(v)}_{L^m_T L^{q}(\R^d)}\lesssim T^{1/m}\norm{\nabla v}_{L^1_TL^{q}(\R^d)} + T^{1/m}\norm{P(\rho_0,c_0)-\widetilde P}_{L^{q}(\R^d)}
\]
and finally, from \cref{lem4}, we have the following H\"older estimate for $I_3(v)$:
\[
\norm{\sigma^{r/4}I_3(v)}_{L^4_T\cC^\alpha_{pw,\mathcal{C}_0}(\R^d)}\lesssim   T^{1/4}\norm{\nabla  v}_{L^1_T\cC^\alpha_{pw,\mathcal{C}_0}(\R^d))}+  T^{1/4} \norm{P(\rho_0,c_0)-\widetilde P}_{\cC^\alpha_{pw,\mathcal{C}_0}(\R^d)}.
\]
We gather all of these estimates in order to obtain: 
\[
\norm{I(v,v)}_{L^m_TL^q(\R^d)}\lesssim \norm{\nabla  v}_{L^1_TL^\infty(\R^d)}\norm{\nabla v}_{L^m_T L^{q}(\R^d)}
+T^{1/m}\norm{\nabla v}_{L^1_TL^{q}(\R^d)} + T^{1/m}\norm{P(\rho_0,c_0)-\widetilde P}_{L^{q}(\R^d)}
\]
and 
\begin{multline*}
\norm{\sigma^{r/4}I(v,v)}_{L^4_T\cC_{pw,\mathcal{C}_0}^{\alpha}(\R^d)}\lesssim  \norm{\nabla  v}_{L^1_T\cC^\alpha_{pw,\mathcal{C}_0}(\R^d)}\norm{\sigma^{r/4}\nabla  v}_{L^4_T\cC^\alpha_{pw,\mathcal{C}_0}(\R^d)}
+T^{1/4}\norm{\nabla  v}_{L^1_T\cC^\alpha_{pw,\mathcal{C}_0}(\R^d))}\\+  T^{1/4} \norm{P(\rho_0,c_0)-\widetilde P}_{\cC^\alpha_{pw,\mathcal{C}_0}(\R^d)}.
\end{multline*}
This proves \eqref{epq45}, \eqref{epq49}, and \eqref{epq53}, and we now move on to estimating the time derivative of $I(v,v)$.
\subsubsection{Lebesgue norm estimates for $\dpt I(v,v)$}\label[section]{esti1-step2}
In this step, we will estimate the $L^2_TL^2(\R^d)$ of $\dpt I(v,v)$ and the $L^{8/d}_TL^4(\R^d)$ norm of $\sigma^{\tfrac{d}{4}} \dpt I(v,v)$. Let us begin  by computing  the time derivative of the third term $I_3(v)$. We have:
\begin{gather}\label{eq3.55}
\dpt I_3(v)=\dpt \adj(D\X_v)\{P(J_v^{-1} \rho_0,c_0)-\widetilde P\}+\rho_0 \partial_\rho  P(J_v^{-1} \rho_0,c_0)\adj(D\X_v)\dpt J^{-1}_v
\end{gather}
and thus one obtains from \eqref{eq3.2} and \eqref{eq3..3} the following: 
\begin{gather}\label{eq3..8}
\norm{\dpt I_3(v)}_{L^2_TL^2(\R^d)}\lesssim \norm{\nabla v}_{L^2_TL^2(\R^d)}\left(\norm{\nabla v}_{L^1_TL^\infty(\R^d)}+1\right).
\end{gather}
On the other hand, the following estimate holds true for $\dpt I_3(v)$:
\begin{gather*}
\norm{\sigma^{\tfrac{d}{4}} \dpt I_3(v)}_{L^{8/d}_TL^4(\R^d)}\lesssim \norm{\nabla v}_{L^{8/d}_TL^4(\R^d)}\left(\norm{\nabla v}_{L^1_T L^\infty(\R^d)}+1\right).
\end{gather*}
We turn to the estimate of the norms of $\dpt I_1(v,v)$ and $\dpt I_2(v,v)$. First, we begin by writing:
\begin{align}
    \dpt I_1(v,v)&=2\rho_0\partial_\rho \mu(J_v^{-1}\rho_0,c_0)\dpt J_v^{-1}\{\adj(D\X_{v})\D_{A_v} v-\D v\}\nonumber\\
                &+2\mu(J^{-1}_{v}\rho_0,c_0)\{\adj( D\X_{v})\D_{A_v}   v_t-\D  v_t\}\nonumber\\
                &+2\mu(J^{-1}_{v}\rho_0,c_0)\dpt\adj( D\X_{v})\D_{A_v}  v\nonumber\\
                &+2\mu(J^{-1}_{v}\rho_0,c_0)\adj (D\X_{v})\left(D v\cdot \dpt A_v+ \dpt A^{T}_v\cdot\nabla v\right).\label{eq3.5}
\end{align}
Using \cref{lem3} and the expression of the time derivative of $J_v^{-1}$ \eqref{eq3.2}, the following bounds hold for the first term of the expression above:
\begin{align*}
    \norm{2\rho_0\partial_\rho \mu(J_v^{-1}\rho_0,c_0)\dpt J_v^{-1}\{\adj(D\X_{v})\D_{A_v} v-\D v\}}_{L^2_TL^2(\R^d)}\lesssim \norm{\nabla v}_{L^1_TL^\infty(\R^d)}\norm{\nabla v}_{L^4_TL^4(\R^d)}^2
\end{align*}
and 
\begin{align*}
    \norm{\sigma^{\tfrac{d}{4}}\rho_0\partial_\rho \mu(J_v^{-1}\rho_0,c_0)\dpt J_v^{-1}\{\adj(D\X_{v})\D_{A_v} v-\D v\}}_{L^{8/d}_T L^4(\R^d)}\lesssim \norm{\nabla v}_{L^1_TL^\infty(\R^d)}\norm{\nabla v}_{L^{16/d}_TL^8(\R^d)}^2.
\end{align*}
The same argument helps us to control the second term of the right hand side of \eqref{eq3.5} as follows:
\[
 \norm{2\mu(J^{-1}_{v}\rho_0,c_0)\{\adj( D\X_{v})\D_{A_v} v_t-\D   v_t\}}_{L^2_TL^2( \R^d)}
 \leqslant \norm{\nabla v}_{L^1_TL^\infty(\R^d)}\norm{\nabla v_t}_{L^2_TL^2(\R^d)}
\]
and
\[
 \norm{\sigma^{\tfrac{d}{4}}\mu(J^{-1}_{v}\rho_0,c_0)\{\adj( D\X_{v})\D_{A_v} v_t-\D   v_t\}}_{L^{8/d}_TL^4( \R^d)}
 \leqslant \norm{\nabla v}_{L^1_TL^\infty(\R^d)}\norm{\sigma^{\tfrac{d}{4}}\nabla v_t}_{L^{8/d}_TL^4(\R^d)}.
\]
Next, from  the expression of the time derivative of $A_v$, \eqref{eq3.6}, we have: 
\begin{align*}
D v\cdot \dpt A_v+ \dpt A^{T}_v\cdot\nabla v&=-(D v\cdot A_v)^2- (A_v^T\cdot\nabla v)^2
\end{align*}
thus, the others terms of \eqref{eq3.5} are controlled as follows:
\begin{align*}
    \norm{\mu(J^{-1}_{v}\rho_0,c_0)\dpt\adj( D\X_{v})\D_{A_v}  v;\;\, &\mu(J^{-1}_{v}\rho_0,c_0)\adj (D\X_{v})\{D v\cdot \dpt A_v+ \dpt A^{T}_v\cdot\nabla v\}}_{L^2_TL^2(\R^d)}\lesssim\norm{\nabla v}_{L^4_TL^4(\R^d)}^2
\end{align*}
and 
\begin{align*}
    \norm{\sigma^{\tfrac{d}{4}}\mu(J^{-1}_{v}\rho_0,c_0)\dpt\adj( D\X_{v})\D_{A_v}  v;\; &\sigma^{\tfrac{d}{4}}\mu(J^{-1}_{v}\rho_0,c_0)\adj (D\X_{v})\{D v\cdot \dpt A_v+ \dpt A^{T}_v\cdot\nabla v\}}_{L^{8/d}_TL^4(\R^d)}\lesssim\norm{\nabla v}_{L^{16/d}_TL^8 (\R^d)}^2.
\end{align*}
Gathering all of these computations and noticing that the estimate for $\dpt I_1(v,v)$  and $\sigma^{\tfrac{d}{4}} I_1(v,v)$ also holds for $\dpt I_2(v,v)$  and $\sigma^{\tfrac{d}{4}} I_2(v,v)$ respectively, one has:
\begin{gather*}
    \norm{\dpt I(v,v)}_{L^2_TL^2(\R^d)}\lesssim \left(\norm{\nabla v}_{L^2_TL^2(\R^d)}+\norm{\nabla v}_{L^4_TL^4(\R^d)}^2\right)\left(\norm{\nabla v}_{L^1_TL^\infty(\R^d)}+1\right)
    +\norm{\nabla v}_{L^1_TL^\infty(\R^d)}\norm{\nabla v_t}_{L^2_TL^2(\R^d)}
\end{gather*}
and
\begin{multline*}
    \norm{\sigma^{\tfrac{d}{4}} \dpt I(v,v)}_{L^{8/d}_TL^4(\R^d)}\lesssim \left(\norm{\nabla v}_{L^{8/d}_TL^4(\R^d)}+\norm{\nabla v}_{L^{16/d}_TL^8(\R^d)}^2\right)\left(\norm{\nabla v}_{L^1_T L^\infty(\R^d)}+1\right)\\
    +\norm{\nabla v}_{L^1_TL^\infty(\R^d)}\norm{\sigma^{\tfrac{d}{4}}\nabla v_t}_{L^{8/d}_TL^4(\R^d)}.
\end{multline*}
This leads to \eqref{epq47} and \eqref{epq50}.
\subsubsection{$L^2_TL^2(\R^d)$ norm estimate for $\sigma \partial_{tt} I(v,v)$}\label[section]{esti1-step3}
As above, we may begin by computing the second time derivative of $I(v,v)$. First, from \eqref{eq3.55}, we have the following expression for $\partial_{tt} I_3(v)$:
\[
\partial_{tt} I_3(v,v)= \partial_{tt}I_{3,1}+ \partial_{tt}I_{3,2}
\]
where
\begin{gather}\label{eq3.56}
\partial_{tt}I_{3,1}=\partial_{tt} \adj(D\X_v)\{P(J_v^{-1} \rho_0,c_0)-\widetilde P\}+ \rho_0 \partial_\rho P(J^{-1}_v\rho_0,c_0)\dpt J^{-1}_v \dpt \adj(D\X_v)
\end{gather}
and 
\begin{align}
\partial_{tt}I_{3,2}&=\rho_0^2 \partial_\rho^2 P(J_v^{-1}\rho_0,c_0) (\dpt J^{-1}_v)^2\adj(D \X_v)+\rho_0 \partial_\rho P(J_v^{-1}\rho_0,c_0)\dpt J_v^{-1}\dpt \adj (D \X_v)\nonumber\\
&+\rho_0 \partial_\rho P(J_v^{-1}\rho_0,c_0)\partial_{tt} J_v^{-1} \adj (D \X_v).\label{eq3.57}
\end{align}
Using the expression of the time derivative of $J^{-1}_v$, \eqref{eq3.2} and of $\adj(D\X_v)$, \eqref{eq3..3}, the last term of \eqref{eq3.56} and the first two terms of the right hand side of \eqref{eq3.57} can be estimated as follows:
\begin{gather*}
\norm{\sigma\rho_0 \partial_\rho P(J^{-1}\rho_0,c_0)\dpt J^{-1}_v \dpt \adj(D\X_v),\;\sigma\rho_0^2 \partial_\rho^2 P(J_v^{-1}\rho_0,c_0) (\dpt J^{-1}_v)^2\adj(D \X_v)}_{L^2_TL^2(\R^d)}
\lesssim \norm{\nabla v}_{L^4_TL^4(\R^d)}^2.
\end{gather*}
Owing to the expression of the second time derivative of $J^{-1}_v$ \eqref{eq3.58} and $\adj (D\X_v)$ \eqref{eq3.59}, the remainder terms of $\partial_{tt} I_3(v)$ are estimated as follows:
\begin{align*}
\norm{\sigma \partial_{tt} \adj(D\X_v)\{P(J_v^{-1} \rho_0,c_0)-\widetilde P\},&\; \sigma \rho_0 \partial_\rho P(J_v^{-1}\rho_0,c_0)\partial_{tt} J_v^{-1} \adj (D \X_v)}_{L^2_TL^2(\R^d)}\lesssim \norm{\sigma \nabla v_t}_{L^2_TL^2(\R^d)}+\norm{\nabla v}_{L^4_TL^4(\R^d)}^2.
\end{align*}
Gathering all of these estimates, one has:
\begin{gather}\label{eq3.70}
\norm{\sigma \partial_{tt} I_3(v)}_{L^2_TL^2(\R^d)}\lesssim  \norm{\sigma \nabla v_t}_{L^2_TL^2(\R^d)}+\norm{\nabla v}_{L^4_TL^4(\R^d)}^2,
\end{gather}
and we turn to the $L^2_TL^2(\R^d)$ norm estimate of $\sigma \partial_{tt} I_1(v,v)$. From \eqref{eq3.5}, we express: 
\[
\partial_{tt} I_1(v,v)=\sum_{k}^{4}\partial_{tt} I_{1,k}(v,v)
\]
where:
\begin{align*}
\partial_{tt}I_{1,1}(v,v)&=2(\rho_0)^2\partial_\rho^2\mu(J_v^{-1}\rho_0,c_0)(\dpt J_v^{-1})^2\{\adj(D\X_{v})\D_{A_v} v-\D v\}\\
       &+2\rho_0\partial_\rho\mu(J_v^{-1}\rho_0,c_0)\partial_{tt} J_v^{-1}\{\adj(D\X_{v})\D_{A_v} v-\D v\}\\
       &+2\rho_0\partial_\rho\mu(J_v^{-1}\rho_0,c_0)\dpt J_v^{-1}\{\adj(D\X_{v})\D_{A_v} v_t-\D v_t\}\\
       &+2\rho_0\partial_\rho\mu(J_v^{-1}\rho_0,c_0)\dpt J_v^{-1}\dpt \adj(D\X_{v})\D_{A_v} v\\
       &+2\rho_0\partial_\rho\mu(J_v^{-1}\rho_0,c_0)\dpt J_v^{-1}\adj(D\X_{v})\{D v\cdot  \partial_t A_v + \partial_t A_v^T\cdot \nabla v\},
\end{align*}
\begin{align*}
   \partial_{tt} I_{1,2}(v,v)&=2\rho_0\partial_\rho\mu(J^{-1}_{v}\rho_0,c_0)\dpt J^{-1}_{v} \{\adj( D\X_{v})\D_{A_v}   v_t-\D  v_t\}\\
                &+2\mu(J^{-1}_{v}\rho_0,c_0)\{\adj( D\X_{v})\D_{A_v}   v_{tt}-\D  v_{tt}\}\\
                &+2\mu(J^{-1}_{v}\rho_0,c_0)\dpt \adj( D\X_{v})\D_{A_v}   v_t\\
                &+2\mu(J^{-1}_{v}\rho_0,c_0)\adj( D\X_{v})\{\D v_t\cdot \dpt A_v  +\dpt A_v^T\cdot \nabla v_t\},
\end{align*}
\begin{align*}
    \partial_{tt} I_{1,3}(v,v)&=2\rho_0\partial_\rho\mu(J^{-1}_{v}\rho_0,c_0)\dpt J^{-1}_{v} \dpt\adj( D\X_{v})\D_{A_v}  v+
    2\mu(J^{-1}_{v}\rho_0,c_0)\partial_{tt}\adj( D\X_{v})\D_{A_v}  v\\
    &+2\mu(J^{-1}_{v}\rho_0,c_0)\dpt\adj( D\X_{v})\{\D_{A_v}  v_t+ D v \cdot\dpt A_v+\dpt A_v^T\cdot \nabla v\}
\end{align*}
and finally, 
\begin{align*}
    \partial_{tt} I_{1,4}(v,v)&=2\rho_0\partial_\rho\mu(J^{-1}_{v}\rho_0,c_0)\dpt J^{-1}_v\adj (D\X_{v})\left(D v\cdot \dpt A_v+ \dpt A^{T}_v\cdot\nabla v\right)\\
                              &+2\mu(J^{-1}_{v}\rho_0,c_0)\dpt \adj (D\X_{v})\left(D v\cdot \dpt A_v+ \dpt A^{T}_v\cdot\nabla v\right)\\
                              &+2\mu(J^{-1}_{v}\rho_0,c_0)\adj (D\X_{v})\left(D v_t\cdot \dpt A_v+ \dpt A^{T}_v\cdot\nabla v_t+D v\cdot \partial_{tt} A_v+ \partial_{tt} A^{T}_v\cdot\nabla v\right).
\end{align*}
Combining \eqref{eq3.2}, \eqref{eq3.6}, \eqref{eq3..3}, \eqref{eq3.58} and \eqref{eq3.59} one has:
\begin{align*}
\norm{\sigma\partial_{tt} I_1(v,v)}_{L^2_TL^2(\R^d)}&\lesssim  \norm{\nabla v}_{L^1_TL^\infty(\R^d)}\left(\norm{\nabla v}_{L^6_TL^6(\R^d)}^3
+\norm{\sigma^{\tfrac{d}{4}}\nabla v_t}_{L^{8/d}_T L^4(\R^d)}\norm{\nabla v}_{L^{8/(4-d)}_TL^4(\R^d)}+\norm{\sigma\nabla v_{tt}}_{L^2_TL^2(\R^d)}\right)\\
         &+\norm{\sigma^{\tfrac{d}{4}}\nabla v_t}_{L^{8/d}_T L^4(\R^d)}\norm{\nabla v}_{L^{8/(4-d)}_T L^4(\R^d)}+\norm{\nabla v}_{L^6_TL^6(\R^d)}^3.
\end{align*}
Finally, the same estimate holds true for the $L^2_TL^2(\R^d)$ norm of $\sigma \partial_{tt} I_2(v,v)$ and gathering \eqref{eq3.70} we obtain \eqref{epq48}. We turn to the estimation of the norms of the viscosity.
\subsubsection{Lebesgue and H\"older norm estimate for viscosity}\label[section]{esti1-step4}
This step makes more sense to the estimate of the lower, upper bound and the H\"older regularity of the viscosity $\mu$ and $\lambda$. With the help of  
\cref{lem3}, we have the following estimates for $\mu$ 
\[
\begin{cases}
\norm{\mu(J^{-1}_v\rho_0,c_0)-\widetilde\mu }_{L^\infty_TL^\infty(\R^d)}&\leqslant \norm{\mu(\rho_0,c_0)-\widetilde\mu}_{L^\infty(\R^d)}+ \norm{\partial_\rho\mu}_{L^\infty}\norm{\nabla v}_{L^1_T L^\infty(\R^d)},\\
\norm{\llbracket \mu(J^{-1}_v\rho_0,c_0)\rrbracket}_{L^\infty_T\mathcal{C}_0}&\leqslant \norm{\llbracket\mu(\rho_0,c_0)\rrbracket}_{L^\infty_T\mathcal{C}_0}+ \norm{\partial_\rho\mu}_{L^\infty}\norm{\nabla v}_{L^1_T L^\infty(\R^d)}
\end{cases}
\]
and a similar estimates hold for $\lambda(J^{-1}_v\rho_0,c_0)$. Also, one has the following H\"older estimate:
\[
\norm{\mu(J^{-1}_v\rho_0,c_0)}_{L^\infty_T\dot \cC^\alpha_{pw,\mathcal{C}_0}(\R^d))}\leqslant \norm{\mu(\rho_0,c_0)}_{\dot \cC^\alpha_{pw,\mathcal{C}_0}(\R^d)}+ \norm{\partial_\rho \mu}_{L^\infty}\norm{\nabla v}_{L^1_T \cC^\alpha_{pw,\mathcal{C}_0}(\R^d)}.
\]
The same estimate also holds for $\lambda(J^{-1}_v\rho_0,c_0)$ and we obtain \eqref{c4.40}. We now turn to the estimate for the time derivative of the viscosity.
\subsubsection{Estimates for time derivative of the viscosity}\label[section]{esti1-step5}
We estimate the norm of the time derivative 
of $\mu(J^{-1}_v\rho_0,c_0)$ and $\lambda(J^{-1}_v\rho_0,c_0)$ as required in \eqref{ep1.2}. First, we express:
\begin{gather}\label{eq3.60}
\dpt \mu(J^{-1}_v\rho_0,c_0)= -\rho_0 \partial_\rho\mu (J^{-1}_v\rho_0,c_0)J_v^{-1}\dvg_{A_v} (v),
\end{gather}
and it transpires that:
\[
\norm{\dpt \mu(J^{-1}_v\rho_0,c_0)}_{L^m_TL^q(\R^d)}\lesssim \norm{\nabla v}_{L^m_TL^q(\R^d)},
\quad
\norm{\sqrt{\sigma}\dpt \mu(J^{-1}_v\rho_0,c_0)}_{L^2_TL^\infty(\R^d)}\lesssim \norm{\sqrt\sigma\nabla v}_{L^2_TL^\infty(\R^d)}.
\]
A similar estimates hold true for the time derivative of $\lambda(J^{-1}_v\rho_0,c_0)$ and implicitly \eqref{c4.41}-\eqref{c4.42} follow. We turn to the estimates of the norm of $\partial_{tt}\mu$ and $\partial_{tt}\lambda$ as required in \eqref{ep1.2}.
\subsubsection{Estimates for the second time derivative of the viscosity}\label[section]{esti1-step6}
From \eqref{eq3.60}, it is obvious that the second time derivative of $\mu(J^{-1}_v\rho_0,c_0)$
reads:
\begin{align*}
\partial_{tt}\mu(J^{-1}_v\rho_0,c_0)&= \rho_0^2 \partial^2_\rho\mu(J^{-1}_v\rho_0,c_0)(J_v^{-1}\dvg_{A_v} (v))^2+\rho_0\partial_\rho \mu(J^{-1}_v\rho_0,c_0)J_v^{-1}(\dvg_{A_v} (v))^2\\
&-\rho_0 \partial_\rho\mu(J^{-1}_v\rho_0,c_0)J_v^{-1}\dvg_{A_v} v_t-\rho_0 \partial_\rho\mu(J^{-1}_v\rho_0,c_0)J_v^{-1}\dpt A_v^T\colon\nabla v,
\end{align*}
hence:
\[
\norm{\sigma\partial_{tt}\mu(J^{-1}_v\rho_0,c_0)}_{L^{2}_TL^4( \R^d)}\lesssim 
\norm{\nabla v}_{L^{4}_TL^8 (\R^d)}^2+\norm{\sigma\nabla v_t}_{L^{2}_TL^4( \R^d)}.
\]
A similar estimate holds true for the second time derivative of $\lambda(J^{-1}_v\rho_0,c_0)$. This leads to \eqref{c4.43} and concludes the estimates used in \cref{prop41}.

\subsection*{Acknowledgment}
  I would like to express my gratitude to my PhD advisors, Cosmin Burtea and David Gérard-Varet, for their insightful discussions and thorough review of this work. I also sincerely thank the referees for their careful reading and valuable comments and suggestions on the manuscript.

\subsection*{Funding}
This project has received funding
from the European Union’s Horizon 2020 research and innovation
program under the Marie Skłodowska-Curie grant agreement No 945332.
I am grateful for the support of the SingFlows project grant (ANR-18- CE40-0027) of the French National Research Agency (ANR).
This work has been partially supported by the project CRISIS (ANR-20-CE40-0020-01), operated by the French National Research Agency (ANR). 

\appendix
\section{Quantitative H\"older estimates}\label[appendix]{app1}
This section contains two results concerning the regularity of even-order Riesz transform of discontinuous functions. The first one is a recent contribution by Gancedo and Garc\'ia-Ju\'arez \cite{gancedo2021quantitative}. Roughly speaking, they established that these operators are continuous on the space of piecewise H\"older continuous functions, provided that the  discontinuity surfaces 
enjoy certain regularities.  The second result, obtained in Hoff \cite{hoff2002dynamics} for $d=2$, concerns the extension of discontinuous functions.

The quantitative H\"older estimates for even-order Riesz operator obtained by Gancedo and Garc\'ia-Ju\'arez reads as follows: 
\begin{lemm}[\cite{gancedo2021quantitative,gancedo2021global}]\label[lemma]{lemA2}
   Assume the hypotheses in  \cref{notabene} hold for $\mathcal{C}$ and consider a singular integral operator $\mathcal{T}$:
   \begin{gather}
    \mathcal{T}(g)(x)=p.v.\int_{\R^d} K(x-y) g(y)dy
   \end{gather}
    where the kernel $K$ is:
   \[
     K(x)=\dfrac{P_{2l}(x)}{\abs{x}^{d+2l}}.
   \]
Here, $P_{2l}$ is a homogeneous polynomial in $\R^d$ of degree $2l$, $l\in \N\setminus\{0\}$.
\begin{itemize}
    \item \textbf{$L^\infty(\R^d)$-Estimate.}
        Given $\alpha'\in (0,1)$ and  $q>d/\alpha'$, there exists a constant $C= C(d,q,\alpha')$ such that for all
    $g\in \cC^{\alpha'}_{pw,\mathcal C}(\R^d)$ and $w\in \cC^{\alpha'}_{pw,\mathcal C}(\R^d)\cap L^q(\R^d)$:
     \begin{align}
    \norm{\mathcal T(g w)}_{L^\infty(\R^d)} &\leqslant C\left(\norm{g}_{\dot \cC^{\alpha'}_{pw,\mathcal C}(\R^d)}+\ell_{\vph}^{-\tfrac{d}{q}}\norm{g}_{L^\infty(\R^d)} \right)\norm{w}_{L^{q}(\R^d)} +C\norm{g}_{L^\infty(\R^d)}\norm{w}_{\cC^{\alpha'}_{pw,\mathcal C}(\R^d)}.\label{c3.31}
\end{align}
\item  \textbf{$\dot \cC^\alpha_{pw,\mathcal C}(\R^d)$-Estimate.} There exists a constant $C=C(d,\alpha)$ such that for all $w\in \cC^\alpha_{pw,\mathcal C}(\R^d)$:
    \begin{gather}\label{ez1.96}
 \norm{\mathcal{T} (w)}_{\dot \cC^\alpha_{pw,\mathcal C}(\R^d)}\leqslant  C\norm{w}_{\dot \cC^\alpha_{pw,\mathcal C}(\R^d)}+ C\left(1+\abs{\mathcal{C}}\right) \norm{w}_{L^\infty(\R^d)}
\mathfrak{P}\big(\norm{\mathcal{C}}_{\text{Lip}}+ \abs*{\mathcal{C}}_*\big)\left|\mathcal{\mathcal{C}}\right|_{\dot \cC^{1+\alpha}},
\end{gather}
where $\mathfrak{P}$ is a polynomial depending on $P_{2l}$.
\end{itemize}
\end{lemm}

For $d=2$, \cref{lemA2} summarizes  the computations  in Step 4 and Step 5 of \cite{gancedo2021global}. It should be noted that the proof uses a geometric lemma derived in \cite{bertozzi1993global}  for $d=2$. This extra cancellation property of even-order Riesz operators is also valid in higher dimensions. We refer to \cite[Section 3]{mateu2009extra}, and more precisely to the estimate of the term in equation (12). Let us now move on to the second result, which concerns 
the extension of piecewise H\"older continuous functions.
\begin{lemm}\label[lemma]{exten}
    Assume the hypotheses in  \cref{notabene} hold  for $\mathcal{C}$  and let $q\in [1,\infty]$. For $g\in L^q(\R^d)\cap \cC^\alpha_{pw,\mathcal C}(\R^d)$ there exists $g^e\in \cC^\alpha(\R^d)$ such that $g^e=g$ on $D$ (or $D^c$)
    and:
    \begin{gather}\label{c3.16}
        \begin{cases}
            \norm{g^e}_{L^\infty(\R^d)}                      &\leqslant C \norm{g}_{L^\infty(\R^d)},\\
            \norm{g-g^e}_{L^q(\R^d)}                         &\leqslant C \ell_{\vph}^{1/q}\norm{\llbracket g\rrbracket}_{L^q(\mathcal{C})},\\
            \norm{g^e}_{\dot \cC^\alpha(\R^d)}+\norm{g-g^e}_{\dot \cC^\alpha_{pw,\mathcal C}(\R^d)}&\leqslant C \norm{g}_{\dot\cC^\alpha_{pw,\mathcal C}(\R^d)}+ C\norm{\llbracket g\rrbracket}_{L^\infty (\mathcal{C})}\ell^{-\alpha}_{\vph}.
        \end{cases}
    \end{gather}
    The constant $C$ above is independent from $q,\,\alpha$ and $\mathcal{C}$.
\end{lemm}
For $d=2$, the lemma above does not differ from Lemma 5.2 in \cite{hoff2002dynamics}, except for the dependence on the interface regularity.  The proof for $d=3$ follows almost the same lines.

\dem[Proof of \cref{exten}]
Given that $\mathcal{C}\in \cC^{1+\alpha}$, for each $x=(x',x_d)\in \mathcal{C}$ \footnote{For $y=(y_1,\dots,y_d)\in \R^d$ we use the notation $y=(y',y_d)$, where $y'=(y_1,\dots, y_{d-1})$.} there exist a neighborhood  $W'$ of $x'$ in $\R^{d-1}$ and a function $z\colon W'\mapsto (x_d-\ell, x_d+\ell)\in \cC^{1+\alpha}$ such that:
\[
\mathcal{C}\cap W= \{ (y',z(y'))\colon y'\in W' \}, \quad\text{ where }\quad W=W'\times (x_d-\ell, x_d+\ell).
\]
For all $y\in \R^d$, we have:
\[
\abs{\nabla \vph (y)}\geqslant  \abs{\nabla\vph (x)}-\abs{\nabla \vph(x)-\nabla\vph(y)}\geqslant \abs{\nabla\vph (x)}-\abs{x-y}^\alpha\norm{\nabla\vph}_{\dot \cC^\alpha}.
\]
Therefore, for $y$ such that:
\[
\abs{x-y}^\alpha \norm{\nabla \vph}_{\dot \cC^\alpha} \leqslant \dfrac{1}{2}\norm{\nabla\vph}_{\text{inf}},
\]
it follows that 
\[
\abs{\nabla \vph (y)} \geqslant \dfrac{1}{2}\norm{\nabla\vph}_{\text{inf}}.
\]
In consequence, we can consider the local chart $(W',z)$ provided by the implicit function theorem,  and  $W'$ can be taken as 
the $(d-1)$-dimensional ball centered at $x'$ with radius $\ell=C_d\ell_{\vph}$. Furthermore, we might have:
\begin{gather}\label{c3.19}
\norm{\nabla z}_{L^\infty(W')}\leqslant 2.
\end{gather}

We now cover $\mathcal{C}$ with finite open bounded sets $W_j$, $j=2,\dots, N$  and we set $W_0=D$ and $W_1=\R^d\setminus \overline D$. For each $j=2,3,\dots,N$, we define the projection:
\begin{align*}\vspace{0,2cm}
   \Xi_j \colon W_j\;&\; \to W_j\cap \mathcal{C} \\ 
     x=(x',x_d)&\;\;\mapsto \Xi_j(x)=(x',z_j(x'))   
\end{align*} 
and the local extension $g^e_j$:
\begin{gather}\label{c3.22}
g^e_j(x)=
    \begin{cases}
        g(x) &\text{if} \quad x\in W_j\cap \overline{D},\\
        g(x)+\llbracket g\left(\Xi_j(x)\right)\rrbracket &\text{if} \quad x\in W_j \cap (\R^d\setminus\overline{D}).
    \end{cases}
\end{gather}
Clearly, $g^e_j$ is $\alpha$-H\"older continuous on $W_j\cap \overline{D}$ and $W_j \cap (\R^d\setminus D)$. However, 
it is even more than that: $g_j^e$ is $\alpha$-H\"older in the whole $W_j$ with: 
\begin{gather}\label{c3.21}
\begin{cases}
     \norm{g_j^e-g}_{L^\infty(W_j)}&\leqslant  \norm{\llbracket g\rrbracket}_{L^\infty(W_j\cap \mathcal{C})},\\
    \norm{g_j^e}_{\dot \cC^\alpha(W_j)}&\leqslant  C\norm{g}_{\dot \cC^\alpha_{pw,\mathcal C}(W_j)}.
\end{cases}
\end{gather}
Deriving $\eqref{c3.21}_1$, and obtaining the estimate of the piecewise H\"older norm  for $g_j^e$ is straightforward from \eqref{c3.22}. Next, we consider $x\in W_j\cap \overline{D}$ and $y\in W_j \cap (\R^d\setminus \overline D)$ and we express:
\begin{align}
\abs{g^e_j(x)-g^e_j(y)}&= \abs{g(x)-g^+(y',z_j(y'))+g^-(y',z_j(y'))-g(y)}\nonumber\\
                       &\leqslant \norm{g}_{\dot \cC^\alpha_{pw,\mathcal C}(W_j)}\left(\abs{x'-y'}^\alpha+\abs{x_d-z_j(y')}^\alpha+\abs{y_d-z_j(y')}^\alpha\right).\label{c3.20}
\end{align}
Recalling  $y_d\leqslant z_j(y')$ and  $z_j(x')\leqslant  x_d$, we have:
\begin{itemize}
    \item If $z_j(y')\leqslant x_d$, then $y_d\leqslant z_j(y')\leqslant  x_d $ and hence:
    \[
    \abs{x_d-z_j(y')}+ \abs{y_d-z_j(y')}\leqslant \abs{x_d-y_d}.
    \]
    \item If $z_j(y')> x_d$, then $z_j(x')\leqslant x_d\leqslant z_j(y')$ and using \eqref{c3.19} we find:
    \[
   \abs{x_d-z_j(x')}+ \abs{x_d-z_j(y')}\leqslant \abs{z_j(x')-z_j(y')}\leqslant 2 \abs{x'-y'}.
    \]
    Whence:
    \[
    \abs{y_d-z_j(y')}\leqslant \abs{y_d-x_d}+\abs{x_d-z_j(y')}\leqslant 3\abs{x-y}.
    \]
\end{itemize}
Summing up, \eqref{c3.20} becomes:
\[
\abs{g^e_j(x)-g^e_j(y)}\leqslant C \norm{g}_{\dot \cC^\alpha_{pw,\mathcal C}(W_j)}\abs{x-y}^\alpha
\]
and $\eqref{c3.21}_2$ follows.

The global extension $g^e$ reads: 
\begin{gather}\label{c3.17}
g^e(x)= g(x)\left[ \psi_{0}(x)+ \psi_{1}(x)\right]+\sum_{j=2}^{N} g^e_j(x) \psi_j(x),
\end{gather}
where $(\psi_j)_j$ is a partition of unity subordinate to the cover $(W_j)_{j}$ of $\R^d$. Given that $\psi_j$ is compactly supported in $W_j$,  $g\psi_0$, $g\psi_1$ and $g^e_j\psi_j$, $j\geqslant 2$, extend over $\R^d$ into $\cC^\alpha(\R^d)$ functions. Therefore,  $g^e$ itself belongs to $\cC^\alpha(\R^d)$, and  $\eqref{c3.16}_{1,2}$ immediately holds for $q=\infty$ from \eqref{c3.21}. From \eqref{c3.17}, we observe: 
\begin{gather}\label{c3.23}
g^e(x)=g(x)+ \sum_{j=2}^{N} \left(g^e_j(x)-g(x)\right) \psi_j(x),
\end{gather}
and whence (recall \eqref{c3.22} and the local finiteness of the cover):
\[
 g^e(x)=g(x)\quad\text{for}\quad x\in \overline {D} \quad\text{and}\quad \norm{g^e-g}_{\dot \cC^\alpha_{pw,\mathcal C}(\R^d)}\leqslant C\norm{g}_{\dot \cC^\alpha_{pw,\mathcal C}(\R^d)}+ C\norm{\llbracket g\rrbracket}_{L^\infty(\mathcal{C})}\sup_{j\geqslant 2}\norm{\psi_j}_{\dot\cC^\alpha(\R^d)}.
\]

Difference in \eqref{c3.17} yields:
\begin{align}
2\big(g^e(x)-g^e(y)\big)&= \big(g(x)-g(y)\big)\big(\psi_0(x)+\psi_1(x)+\psi_0(y)+\psi_1(y)\big)\notag\\
&+ \big(g(x)+g(y)\big)\sum_{j=0}^{j=1}\left(\psi_j(x)-\psi_j(y)\right)+\sum_{j =2}^N \left( g^e_j(x)-g_j^e(y)\right)\big(\psi_j(x)+\psi_j(y)\big)\notag\\
             & +\sum_{j= 2}^N\big( g_j^e(x)+g_j^e(y)\big)\left(\psi_j(x)-\psi_j(y)\right)\notag\\
             &=\big(g(x)-g(y)\big)\big(\psi_0(x)+\psi_1(x)+\psi_0(y)+\psi_1(y)\big)\notag\\
             &+\sum_{j= 2}^N \left( g^e_j(x)-g_j^e(y)\right)\big(\psi_j(x)+\psi_j(y)\big)+\sum_{j=2}^N \big((g^e_j-g)(x)+(g^e_j-g)(y)\big)\left(\psi_j(x)-\psi_j(y)\right)\notag\\
             &+ \big( g(x)+g(y)\big)\sum_{j=0}^{N}\left(\psi_j(x)-\psi_j(y)\right),\label{c7.12}
\end{align}
where we notice that the last term  vanishes. 
This leads to
\[
\norm{g^e}_{\dot \cC^\alpha(\R^d)}\leqslant C \norm{g}_{\dot\cC^\alpha_{pw,\mathcal C}(\R^d)}+ C\norm{\llbracket g\rrbracket}_{L^\infty (\mathcal{C})}\sup_{j}\norm{\psi_j}_{\dot \cC^\alpha}
\]
and  $\eqref{c3.16}_3$  follows since we have:
\[
\sup_{j}\norm{\psi_j}_{\dot \cC^\alpha(\R^d)}\leqslant C  \ell^{-\alpha}.
\]
We may observe that when $x\in D$ and $y\in \R^d\setminus \overline{D}$, it follows that $\psi_1(x)=0\;\text{ and }\; \psi_0(y)=0$. Next, we choose $\widetilde x\in \mathcal{C}$ (which implies $\psi_0(\widetilde x)=\psi_1(\widetilde x)=0$, since $\psi_0$ and $\psi_1$ vanish on $\R^d\setminus D$ and $\overline{D}$, respectively) such that 
\[
\abs{x-\widetilde x},\, \abs{y-\widetilde x} \leqslant  \abs{x-y}.
\]
Consequently, the first term of the right hand side of  \eqref{c7.12} shall be 
expressed as follows:
\begin{align*}
\big(g(x)-g(y)\big)\big(\psi_0(x)+\psi_1(x)+\psi_0(y)+\psi_1(y)\big)&=\big(g(x)-g(y)\big)\big(\psi_0(x)+\psi_1(y)\big)\\
&=\big(g(x)-g^+(\widetilde x)+g^- (\widetilde x)-g(y)\big) \big(\psi_0(x)+\psi_0(y)\big)\\
&+\big( g^+(\widetilde x)-g^- (\widetilde x)\big) \big(\psi_0(x)-\psi_0(\widetilde x)+\psi_1(y)-\psi_1(\widetilde x)\big).
\end{align*}
We now proceed with the proof of $\eqref{c3.16}_2$ for $q\in [1,\infty)$. From \eqref{c3.23} and H\"older's inequality, we deduce:
\begin{align*}
    \norm{g^e-g}_{L^q(\R^d)}^q&=\int_{D^c} \left|\sum_{j=1}^{N}\left(g^e_j(x)-g(x)\right) \psi_j(x)\right|^qdx\\
                              &\leqslant \int_{D^c}\left(\sum_{j=1}^{N}\psi_j(x)\right)^{q-1}\left(\sum_{j=1}^{N}\abs{g^e_j(x)-g(x)}^q \psi_j(x)\right)\\
                              &\leqslant \sum_{j=1}^{N}\int_{D^c\cap W_j}\abs{g^e_j(x)-g(x)}^q \psi_j(x)dx\\
                              &\leqslant 2\ell\sum_{j=1}^{N}\int_{W_j'}\abs{\llbracket g(\Xi_j(x'))\rrbracket}^q dx' \quad (\text{recall} \;\;\eqref{c3.22})\\
                              &\leqslant C \ell \sum_{j=1}^{N}\int_{W_j\cap \mathcal{C}}\abs{\llbracket g\rrbracket}^qd\sigma
\end{align*}
and $\eqref{c3.16}_2$ follows. This concludes the proof of \cref{exten}.
\section{Study of linear non homogeneous heat equation}\label[appendix]{app2}
This section  addresses  the well-posedness  of the following system of PDE:
\begin{gather}\label{eq3.25}
\begin{cases}
v_t- \widetilde\mu \Delta v-(\widetilde\mu+\widetilde\lambda)\nabla \dvg v= f+\dvg (b_0\cdot v),\\
v_{|t=0}=v_0.
\end{cases}
\end{gather}
The main result is presented in the following theorem.
\begin{theo}\label[theo]{th3}
    Let  $T\in (0, \infty)$, $v_0\in H^1(\R^d)$, $b_0\in W^{1,\infty}(\R^d)$ and $f\in L^2((0,T)\times \R^d)$.  Then, the Cauchy problem \eqref{eq3.25} admits a unique solution $v\in \cC([0,T], H^1(\R^d))$. Moreover  $v_t,\, \nabla^2 v\in L^2 ((0,T)\times\R^d)$.
\end{theo}
\subsection{Preliminary}
 We define
\[
\mathcal{A}_k=
\begin{cases}
    \mathcal{P} \quad \text{if}\quad k=\widetilde\mu,\\
        \mathcal{Q} \quad \text{if}\quad k=2\widetilde\mu+\widetilde\lambda
\end{cases}
\]
where $\mathcal{P}$ (resp. $\mathcal{Q}$) is the projector onto the space of divergence-free (resp. curl-free) vector fields.  
Applying $\mathcal{A}_k$ to \eqref{eq3.25}, we infer that the divergence-free and curl-free part of any regular\footnote{It is sufficient that $\dpt v\in L^2((0,T)\times \R^d)$ and $v\in L^2((0,T), H^2(\R^d))$.} solution $v$ of \eqref{eq3.25} verify:
\begin{gather}\label{c3.12}
\begin{cases}
\dpt \mathcal{A}_k v-k \Delta \mathcal{A}_kv =\mathcal{A}_k f+\mathcal{A}_k\dvg (b_0\cdot v),\\
\mathcal{A}_kv_{|t=0}=\mathcal{A}_kv_0.
\end{cases}
\end{gather}
Conversely, since $\mathcal{P}+ \mathcal{Q}=I$, any  regular vector field $v$ that fulfills \eqref{c3.12} is a solution of \eqref{eq3.25}.  

We shall now consider the following Cauchy problem
\begin{gather}\label{c2.17}
\begin{cases}
\dpt \mathcal{A}_kv-k\Delta \mathcal{A}_kv =\mathcal{A}_k f+\mathcal{A}_k\dvg (b_0\,\colon w),\\
\mathcal{A}_kv_{|t=0}=\mathcal{A}_kv_0,
\end{cases}
\end{gather}
where $w\in \cC([0,T], H^1(\R^d))$ is a given vector field. The solution of this  linear non-homogeneous heat equation is given by 
the Duhamel formula:
\begin{gather}\label{eq3.28}
\mathcal{A}_kv(t) =\mathcal{T}(k t) \mathcal{A}_kv_0+\int_{0}^t  \mathcal{T}(k (t-s))\mathcal{A}_kf(s)ds+\int_{0}^t  \mathcal{T}(k (t-s))\mathcal{A}_k\dvg (b_0\,\colon w)(s)ds,
\end{gather}
where $(\mathcal{T}(t))_t$ is the heat semi-group. It follows  $v\in \cC([0,T], H^1(\R^d))$ with the following estimate:
\begin{gather}\label{eq3.26}
\sup_{[0,T]}\norm{v}_{H^1(\R^d)} \leqslant \norm{v_0}_{H^1(\R^d)}+\sum_k\int_0^T \norm{\mathcal{T}(k (t-s))\mathcal{A}_k\left[f(s)+\dvg (b_0\,\colon w)\right](s)}_{H^1(\R^d)}ds.
\end{gather}
Additionally, based on  \cite[Lemma 2.2]{danchin2020well}, we infer that $\dpt v,\, \nabla^2 v\in L^2((0,T)\times \R^d)$ and that there exists a constant $C_k=C_k(k)>0$ such that we have:
\begin{gather*}
\int_0^T\norm{\dpt \mathcal{A}_kv,\nabla^2 \mathcal{A}_kv}_{L^2(\R^d)}^2 \leqslant C_k \norm{\mathcal{A}_kv_0}_{H^1(\R^d)}^2+C_k\int_0^T\left\|\int_0^t \Delta\mathcal{T}(k (t-s))\mathcal{A}_k\left[f +\dvg (b_0\,\colon w)\right](s)ds\right\|_{L^2(\R^d)}^2dt.
\end{gather*}
This therefore defines a map 
\begin{align*}
\Phi\colon\,\; \cC([0,T], H^1(\R^d))&\quad\to\quad \cC([0,T], H^1(\R^d))\\ 
                              w     & \quad\mapsto \quad \Phi(w)
\end{align*}
where  $\Phi(w)$ is the unique solution of \eqref{c2.17}, which is given by \eqref{eq3.28} and fulfills \eqref{eq3.26}. In the following we will set some tools in order to apply the Banach's fixed-point theorem. 
\subsection{\textbf{Tools}}
We begin by the following Lemma:
\begin{lemm}\label[lemma]{lemC1}
    There exists a constant $C=C\left(\norm{b}_{W^{1,\infty}(\R^d)}\right)$ such that:
    \[
    \forall\; \vph \in H^1(\R^d),\quad\, \forall\, t>0\quad  \norm{\mathcal{T}(k t)\mathcal{A}_k\dvg (b_0\,\colon \vph)}_{H^1(\R^d)}\leqslant C (k t)^{-\tfrac{1}{2}}\norm{\vph }_{H^1(\R^d)}.
    \]
\end{lemm}

\dem 
We first write:
\[
\mathcal{T}(k t)\mathcal{A}_k\dvg (b_0\,\colon \vph)= G_{k t}*\mathcal{A}_k \dvg (b_0\,\colon \vph)=\nabla G_{k t}*\mathcal{A}_k(b_0\,\colon \vph)
\]
where $G_t$ is the heat kernel. Thanks to the continuity of $\mathcal{A}_k$ on $L^2(\R^d)$ and  Young's inequality we obtain: 
\[
\norm{\mathcal{T}(k t)\mathcal{A}_k\dvg (b_0\,\colon \vph)}_{L^2(\R^d)}\leqslant \norm{\nabla G_{k t}}_{L^1(\R^d)}\norm{b_0}_{L^\infty(\R^d)}\norm{\vph}_{L^2(\R^d)},
\]
where
\[
\norm{\nabla G_{t}}_{L^1(\R^d)}=\dfrac{(4\pi t)^{-\tfrac{d}{2}}}{2t}(2\sqrt{t})^{d+1}\abs{\partial B(0,1)}\int_{0}^{\infty} r^d e^{-r^2}dr.
\]
Similarly, we obtain:
\begin{align*}
\norm{\nabla \mathcal{T}(k t)\mathcal{A}_k\dvg (b_0\,\colon \vph)}_{L^2(\R^d)}&\leqslant \norm{\nabla G_{k t}}_{L^1(\R^d)}\norm{\nabla (b_0\,\colon\vph)}_{L^2(\R^d)}\\
&\leqslant \norm{\nabla G_{k t}}_{L^1(\R^d)}\norm{b_0}_{W^{1,\infty}(\R^d)}\norm{\vph}_{H^1(\R^d)}.
\end{align*}
This concludes the proof of \cref{lemC1}.
\enddem
We now turn to the stability of equations \eqref{eq3.28} and to do this,  we set 
\[
R_{k,T}=\sup_{t\in [0,T]}\left\|\mathcal{T}(k t) \mathcal{A}v_0+\int_{0}^t  \mathcal{T}(k (t-s))\mathcal{A}_kf(s)ds\right\|_{H^1(\R^d)}\quad\text{and}\quad R_T=R_{\widetilde\mu,T}+R_{2\widetilde\mu+\widetilde\lambda,T}.
\]
\begin{lemm}\label[lemma]{lemB2} 
  There exist  constants $K>0$ and $\beta_0>0$  such that for all $\beta\geqslant \beta_0$ the condition 
  \begin{gather}\label{eq3.29}
  \forall\, t\in [0,T],\quad \norm{w(t)}_{H^1(\R^d)} \leqslant e^{\beta t} K R_T
  \end{gather}
  implies the same for $\Phi(w)$, for all $w\in \cC([0,T], H^1(\R^d))$.
\end{lemm}
\dem 
From \eqref{eq3.28}, \cref{lemC1} and assumption \eqref{eq3.29} on $w$, we have:
\begin{align*}
\norm{\mathcal{A}_k\Phi(w)(t)}_{H^1(\R^d)}  &\leqslant R_{k,T} +\int_0^t\norm{\mathcal{T}(k (t-s))\mathcal{A}_k\dvg (b_0\,\colon w)(s)}_{H^1(\R^d)}ds\\
&\leqslant R_{k,T}+C k^{-\tfrac{1}{2}}KR_T\int_0^t\dfrac{e^{\beta s}}{\sqrt{t-s}}ds.
\end{align*}
For the condition \eqref{eq3.29} to be satisfied  for $\Phi(w)$, it is sufficient that 
\begin{gather}\label{eq3.30}
1+2C{\widetilde\mu}^{-\tfrac{1}{2}}K\int_0^t\dfrac{e^{\beta s}}{\sqrt{t-s}}ds \leqslant K e^{\beta t}.
\end{gather}
Owing to a change of variable, we have:
\[
    e^{-\beta t}\int_0^t \dfrac{e^{\beta s}}{\sqrt{t-s}}ds=\dfrac{1}{\sqrt{\beta}}\int_{0}^{\beta t} \dfrac{e^{-s}}{\sqrt{s}}ds
    \leqslant\dfrac{1}{\sqrt{\beta}}\int_{0}^{\infty} \dfrac{e^{-s}}{\sqrt{s}}ds
    = \sqrt{\dfrac{2}{\beta}} \int_{0}^{\infty} e^{-\tfrac{r^2}{2}}ds
    =\sqrt{\dfrac{\pi}{\beta}}
\]
hence,
\[
1+2C{\widetilde\mu}^{-\tfrac{1}{2}}K\int_0^t\dfrac{e^{\beta s}}{\sqrt{t-s}}ds\leqslant 1+2CKe^{\beta t} \sqrt{\dfrac{\pi}{\widetilde\mu\beta}}.
\]
We  choose $K$ as 
\[
K=\dfrac{1}{4C}\sqrt{\dfrac{\beta_0 \widetilde\mu}{\pi}}
\]
with $\beta_0$ large such that $K\geqslant 2$, and the condition \eqref{eq3.30} becomes 
\[
1+ Ck^{-\tfrac{1}{2}}K\int_0^t\dfrac{e^{\beta s}}{\sqrt{t-s}}ds\leqslant 1+\dfrac{e^{\beta t}}{2}\sqrt{\dfrac{\beta_0}{\beta}}\leqslant  K e^{\beta t}.
\]
This concludes the proof of \cref{lemB2}.
\enddem
\begin{lemm}\label[lemma]{lemB3}
    There exists  $\beta_0'>0$  such that
    for all $\beta\geqslant \beta'_0$ and for all $v,\,w\in \cC([0,T],H^1(\R^d))$, one has:
    \begin{gather}\label{c3.13}
     \sup_{t\in [0,T]} e^{-\beta t}\norm{\Phi(v)(t)-\Phi(w)(t)}_{H^1(\R^d)}\leqslant \dfrac{1}{2} \sup_{t\in [0,T]} e^{-\beta t}\norm{v(t)-w(t)}_{H^1(\R^d)}.
    \end{gather}
\end{lemm}
\dem 
Starting from \eqref{eq3.28} and \cref{lemC1}, one has: 
\begin{align*}
\norm{\mathcal{A}_k\Phi(v)(t)-\mathcal{A}_k\Phi(w)(t)}_{H^1(\R^d)}&= \left\|\int_{0}^t  \mathcal{T}(k (t-s))\mathcal{A}_k\dvg (b_0\,\colon (v-w))(s)ds\right\|_{H^1(\R^d)}\\
&\leqslant C k^{-\tfrac{1}{2}}\int_0^t\dfrac{1}{\sqrt{t-s}} \norm{v(s)-w(s)}_{H^1(\R^d)} ds\\
&\leqslant Ck^{-\tfrac{1}{2}}\sup_{s\in [0,T]} e^{-\beta s}\norm{v(s)-w(s)}_{H^1(\R^d)} \int_0^t \dfrac{e^{\beta s}}{\sqrt{t-s}}ds.
\end{align*}
It then turns out that, 
\[
\sup_{t\in [0,T]}e^{-\beta t}\norm{\mathcal{A}_k\Phi(v)(t)-\mathcal{A}_k\Phi(w)(t)}_{H^1(\R^d)} \leqslant Ck^{-\tfrac{1}{2}}\sqrt{\dfrac{\pi}{\beta}}
\sup_{t\in [0,T]} e^{-\beta t}\norm{v(t)-w(t)}_{H^1(\R^d)}
\]
and \cref{lemB3} follows by choosing $\beta'_0>0$ large enough.
\enddem
\subsection{\textbf{Final step}}
We now turn to the application of the Banach's fixed-point theorem in the following vector space 
\[
E:=\left\{ v\in \cC([0,T], H^1(\R^d))\colon \sup_{t\in [0,T]}e^{-\beta t}\norm{ v(t)}_{H^1(\R^d)} <\infty \right\}
\]
endowed with the norm:
\[
\forall\, v\in E,\quad \norm{v}_{E}=\sup_{t\in [0,T]}e^{-\beta t}\norm{ v(t)}_{H^1(\R^d)}.
\]
Above, $\beta$ is a large number such that \eqref{eq3.29}-\eqref{c3.13} hold true. It is obvious that $(E, \norm{\,\cdot\,}_E)$ is a Banach space and, in particular the closed ball $E_T= B(0, KR_T)$ is complete with respect to the norm $\norm{\,\cdot\,}_E$. On the first hand, \cref{lemB2} ensures that $\Phi$ maps $E_T$ into itself and the \cref{lemB3} ensures that $\Phi$ is a contracting map. Hence, by Banach's  fixed-point theorem, $\Phi$ admits a unique fixed-point $v$
 in $E_T$ that solves the equation \eqref{eq3.25}. Furthermore,  $v$ verifies 
 \[
 \mathcal{A}_kv =\mathcal{T}(k t) \mathcal{A}v_0+\int_{0}^t  \mathcal{T}(k (t-s))\mathcal{A}_kf(s)ds+\int_{0}^t  \mathcal{T}(k (t-s))\mathcal{A}_k\dvg (b_0\cdot v)(s)ds
 \]
and  \cite[Lemma 2.2]{danchin2020well} implies $v_t, \nabla^2 v\in L^2((0,T)\times \R^d)$.

 \vspace{0,3cm}
 \paragraph{\textbf{Data availability}}
 There are no data sets generated or analyzed in this study, therefore data sharing does not apply.

\vspace{0,3cm}
\paragraph{\textbf{Conflict of interest}} The author declares that they have no conflict of interest.
{\small 
\bibliographystyle{acm}
\bibliography{JMFM-V3/Biblio-JMFM-V3}
}
\end{document}